%% file: 2d-isotropy.tex
\documentclass[nopreprintline,12pt]{elsarticle}
\biboptions{sort&compress}
\usepackage[nottoc]{tocbibind}

\usepackage{amsthm}

\usepackage{amsmath,amssymb,url,tikz,enumitem,color,hyperref,xspace,comment}

\input{headeforpics.tex}
\usetikzlibrary{matrix,arrows}
\tikzstyle{twocell}=[-implies,double equal sign distance]
\tikzstyle{shorttwocell}=[-implies,double equal sign distance,shorten >=7pt, shorten <=7pt]

\numberwithin{equation}{section}
\newtheorem{theorem}[equation]{Theorem}
\newtheorem{lemma}[equation]{Lemma}
\newtheorem{proposition}[equation]{Proposition}
\newtheorem{corollary}[equation]{Corollary}
\theoremstyle{definition}
\newdefinition{definition}[equation]{Definition}
\newtheorem{example}[equation]{Example}
\newdefinition{remark}[equation]{Remark}

\newcommand{\cat}[1]{\ensuremath{\mathbf{#1}}}
\newcommand{\CC}{\cat{C}}
\newcommand{\DD}{\cat{D}}
\newcommand{\EE}{\cat{E}}
\newcommand{\Grp}{\cat{Grp}}
\newcommand{\Mon}{\cat{Mon}}
\newcommand{\Set}{\cat{Set}}
\newcommand{\presh}[1][C]{[\cat{#1}\op,\Set]}
\newcommand{\op}{\ensuremath{{}^{\mathrm{op}}}}
\newcommand{\rev}{\ensuremath{{}^{\mathrm{rev}}}}
\newcommand{\id}[1][]{\ensuremath{\mathrm{id}_{#1}}}
\newcommand{\ie}{\text{i.e.,}\xspace}
\newcommand{\eg}{\text{e.g.}\xspace}
\DeclareMathOperator{\ob}{ob}
\DeclareMathOperator{\nat}{nat}
\DeclareMathOperator{\Aut}{Aut}

\newcommand{\MM}{\mathcal{M}} 
\newcommand{\Mps}{\mathcal{M}_{ps}} 
\newcommand{\ZZ}{\mathcal{Z}} 
\newcommand{\Zps}{\mathcal{Z}_{ps}} 
\newcommand{\LL}{\mathcal{L}} 
\newcommand{\NN}{\mathcal{N}} 
\newcommand{\Pic}[1]{\mathrm{Pic}(\cat{#1})}
\newcommand{\Picb}{\mathrm{Pic}}
\newcommand{\lslice}[2]{#1/\!/\cat{#2}}
\newcommand{\pslice}[2]{#1/\cat{#2}}

\begin{document}

\begin{frontmatter}

\title{Inner autoequivalences in general and those of monoidal categories in particular}

\author{Pieter Hofstra\fnref{fn1}} 
\fntext[fn1]{Pieter Hofstra passed away shortly after giving his final comments on this manuscript.} 
\author{Martti Karvonen\corref{cor1}}
\ead{martti.karvonen@uottawa.ca}
\address{University of Ottawa, Department of Mathematics and Statistics, Ottawa, Canada}
\cortext[cor1]{Corresponding author}

\begin{abstract}
We develop a general theory of (extended) inner autoequivalences of objects of any $2$-category, generalizing the theory of isotropy groups to the $2$-categorical setting. 
We show how dense subcategories let one compute isotropy in the presence of binary coproducts, 
unifying various known one-dimensional results and providing tractable computational tools in the two-dimensional setting. 
In particular, we show that the isotropy $2$-group of a monoidal category coincides with its \emph{Picard $2$-group}, \ie the $2$-group on its weakly invertible objects. 
\end{abstract}


\end{frontmatter}

\tableofcontents

\section{Introduction}\label{sec:intro}

An inner automorphism of a group $G$ is an automorphism $G\to G$ that is of the form $g\cdot (-)\cdot g^{-1}$, \ie one that equals conjugation by an element of $G$. 
For an inner automorphism of $G$, there might be several elements of $G$ that induce it---for instance, all central elements induce the identity on $G$. 
However, a fixed $g\in G$\/ induces more automorphisms than just one on $G$: 
in fact, for any homomorphism $\phi\colon G\to H$\/ it induces an inner automorphism of $H$ by conjugation with $\phi(g)$. 
This family of automorphisms constitutes what~\cite{bergman:inner} calls an \emph{extended inner automorphism}: formally,
 this can be defined as a natural automorphism of the projection functor $G/\cat{Grp}\to \cat{Grp}$\/ (where $G/\cat{Grp}$\/ is the 
 co-slice category and where the projection takes a homomorphism $G \to H$\/ to $H$). 
The main result in loc.\ cit.\ is then that every such extended inner automorphism is 
induced in the above manner by a unique element of $G$, so that there is a natural isomorphism
\[ \Aut(G/\cat{Grp}\to \cat{Grp}) \cong G.\]
 This motivates the following definition.

\begin{definition} Let $\CC$ be a category and $X$ an object of $\CC$. 
Then the (covariant) \emph{isotropy group of $\CC$ at $X$} is the group $\ZZ(X)$ of natural automorphisms of the projection functor $P_X\colon X/\CC\to \CC$. 
\end{definition}

Explicitly, an automorphism 
\begin{equation}\label{eq:alpha} 
\alpha=(\alpha_f)_{f\colon X \to A} \in \ZZ(X)=_{\mathrm{def}} \Aut(P_X\colon X/\CC\to \CC) \end{equation}
 consists of an automorphism $\alpha_f\colon A\to A$ for each $f\colon X\to A$ such that the square
    \[\begin{tikzpicture}
    \matrix (m) [matrix of math nodes,row sep=2em,column sep=4em,minimum width=2em]
    {
     X & X\\
     A & A  \\};
    \path[->]
    (m-1-1) edge node [left] {$f$} (m-2-1)
           edge node [above] {$\alpha_{\id[A]}$} (m-1-2)
    (m-1-2) edge node [right] {$f$} (m-2-2)
    (m-2-1) edge node [below] {$\alpha_{f}$}  (m-2-2);
  \end{tikzpicture}\]
commutes. (In the terminology of universal algebra, $\alpha_f$\/ \emph{extends} $\alpha_{\id[A]}$.) 
Moreover, the naturality of $\alpha_f$\/ then amounts to requiring that for each $g\colon A \to B$, the square
    \[\begin{tikzpicture}
    \matrix (m) [matrix of math nodes,row sep=2em,column sep=4em,minimum width=2em]
    {
     A & A\\
     B & B  \\};
    \path[->]
    (m-1-1) edge node [left] {$g$} (m-2-1)
           edge node [above] {$\alpha_f$} (m-1-2)
    (m-1-2) edge node [right] {$g$} (m-2-2)
    (m-2-1) edge node [below] {$\alpha_{gf}$}  (m-2-2);
  \end{tikzpicture}\]
commutes for any $g\colon A\to B$. 

A morphism $x\colon X\to Y$ induces a homomorphism $\ZZ(X)\to \ZZ(Y)$ as follows: 
first, note that that $x$ induces a functor $x^*\colon Y/\CC\to X/\CC$ fitting into a strictly commuting triangle
    \[\begin{tikzpicture}
    \node (x) at (0,1) {$Y/\CC$};
    \node (a) at (-1,-1) {$X/\CC$};
    \node (b) at (1,-1) {$\CC$};
    \draw[->] (x) to node[left] {$x^*$} (a);
    \draw[->] (x) to node[right] {$P_Y$} (b);
    \draw[->] (a) to node[below] {$P_X$} (b);
    \end{tikzpicture}\]
so that given $\alpha\in\ZZ(X)$ we can define $\ZZ(x)\alpha$ by whiskering along $x^*$, \ie $\ZZ(x)\alpha:=\alpha x^*$. 
In concrete terms, $\ZZ(x)\alpha$ is defined for $f\colon Y\to A$ by $(\ZZ(x)\alpha)=\alpha_{fx}$. Consequently, $\ZZ$ is functorial in $X$; as such, it can be thought
of as rectifying the failure of the assignment $X \mapsto \Aut(X)$\/ to be functorial. Indeed, there is, for each $X$, a comparison homomorphism
\[ \theta_X\colon \ZZ(X) \to \Aut(X). \]

In this viewpoint, one thinks of an arbitrary automorphism $h\colon X\to X$\/ as an ``abstract inner automorphism'' 
if it can be extended to such a family $\alpha$ with $\alpha_{\id} = h$; equivalently, these are the automorphisms in the image of
the comparison map $\theta_X$.  
Modulo size issues (discussed in the next section), we therefore have a functor 
\[ \ZZ\colon \cat{C}\to\cat{Grp}\] 
called the (covariant) \emph{isotropy group} of $\cat{C}$. This functor (or rather,
the contravariant version) was studied in a topos-theoretic context
in~\cite{funketal:crossedtoposes}, and in the context of (essentially) algebraic theories
 in~\cite{hofstraetal:isotropyofalgtheories,hofstraetal:picard}.

In particular, in~\cite{hofstraetal:picard} it is shown that for the $1$-category 
$\cat{StrMonCat}$\/ of strict monoidal categories and strict monoidal functors, 
the isotropy group evaluated at a monoidal category $\CC\in \cat{StrMonCat}$\/ is 
isomorphic to the \emph{strict Picard group} of $\CC$: 
this is the group
of strictly invertible objects of $\CC$; here, an object $X$ of $\CC$\/ is called strictly invertible 
when there exists another object $Y$\/ such that $X\otimes Y=I=Y\otimes X$, 
where $I$ is the tensor unit of $\cat{C}$.

In this work we introduce and investigate a two-dimensional version of the theory. We believe this is motivated in part 
by its intrinsic appeal; however, the example of the Picard group provides a concrete incentive to develop the generalization
to $2$-categories. Indeed, given a non-strict monoidal category, such as the category of modules over a ring, or the category
of vector bundles over a space, we cannot apply directly the aforementioned result of~\cite{hofstraetal:picard}; instead,
we would have to pass to a category of isomorphism classes of objects. Of course, this is similar in spirit to what is done traditionally,
since the Picard group of a ring is usually defined to be the group of isomorphism classes of invertible modules. However,
from the modern perspective of higher categories one associates instead with a (non-strict) monoidal category its Picard $2$-group,
that is, the monoidal groupoid of weakly invertible objects and isomorphisms between them. (For general discussion of
$2$-groups and related structures, see~\cite{baezlauda:2groups}.)

Thus we aim to study, for a $2$-category $\cat{C}$ and an object $X$ of $\CC$, the $2$-group
 of ``extended inner autoequivalences'' of $\CC$ by taking the pseudonatural autoequivalences 
 of the projection  $X/\cat{C}\to\cat{C}$.  Somewhat more explicitly, an object $\alpha$\/ of this $2$-group comprises, 
 for each $1$-cell $f\colon X \to A$,  an equivalence $\alpha_f\colon A \to A$\/ rather than an isomorphism; 
moreover the naturality squares of $\alpha$\/ are now required to commute up to 
 coherent invertible $2$-cells (that are now part of the data of $\alpha$) rather than up to equality.

It turns out that moving to the two-dimensional setting is not as straightforward as one might hope.
The first and most elementary obstacle is the following. For a $2$-category $\cat{C}$, there are various notions of (co-)slice: there is the strict slice, the pseudo-slice,
and the lax slice. Once we choose one of those, we may consider the projection $2$-functor $P_X\colon X/\cat{C} \to \cat{C}$, and then we have another choice to make: 
are we to take $2$-natural, pseudonatural, or lax natural autoequivalences of $P_X$? As it turns out, several of the possible combinations fail. For example,
if we were to consider the pseudo-slice $X/\cat{C}$, and define $\ZZ(X)$\/ to be the $2$-group of pseudonatural autoequivalences of $P_X$, then
while we can still define an action of $\ZZ$\/ on $1$-cells via whiskering (as in the one-dimensional case), $\ZZ$\/  fails to be a $2$-functor, 
because there in general is no well-defined action on $2$-cells. This particular problem vanishes when one instead considers the lax slice and lax natural 
autoequivalences of the projection $X/\!/\cat{C} \to \cat{C}$, 
but then there is another technical hurdle, namely that lax natural transformations do not form the $2$-cells of a bicategory (or even tricategory),
which prevents us from straightforwardly deducing the desired functoriality. 
We therefore have to exercise some caution when setting up the general theory. (See Section~\ref{sec:2disotropy} for details.)
As it turns out, however, there is an important situation in which the aforementioned problem goes away: when the ambient $2$-category has binary coproducts, it makes no difference whether one defines 2-isotropy in terms of the lax slice or in terms of the pseudo-slice. 

A different obstacle, and perhaps a more important one, is that most prior computations 
of isotropy groups proceed by first identifying the isotropy group with a suitable group of definable automorphisms\footnote{The idea that inner automorphisms are related to the definable ones is already present in~\cite{freyd:algebravaluedfunctors}, which inspired the title of this work.} and then computing the latter by intricate syntactic arguments that are not that straightforward to generalize to two dimensions.
Moreover, even with a suitable logical framework for reasoning syntactically about 
(certain classes of) $2$-categories, carrying out syntactic arguments similar to those in~\cite{hofstraetal:isotropyofalgtheories,hofstraetal:picard,parker:IsotropyofGrothendieckToposes,parker:presheaves} seems unwieldy. 
To overcome this issue, we introduce a novel technique for computing isotropy groups, 
relying on the existence of binary coproducts and on a convenient choice of a dense subcategory. 
This technique is already interesting in the one-dimensional case, as it lets us compute various known isotropy groups with high efficiency, 
systematizes known results and sheds light on prior proofs.
 For instance, the aforementioned result of Bergman~\cite{bergman:inner} revolves around considering the 
 coproduct of a group with the free groups on one and two generators. 
 That these two groups appear in the argument is no accident, for they form a dense subcategory of the category of groups. 
 Finally, this technique readily generalizes to the two-dimensional setting, and makes establishing some of the coherence conditions tractable
  by boiling the situation down to a small amount of data.

The plan of the paper is as follows. 
We begin in Section~\ref{sec:1d} by explaining the ``coproducts-plus-density" technique with numerous examples in the one-dimensional case.  
After that, we will briefly review some required background concerning $2$-categories in Section~\ref{sec:2dbackground} 
before defining isotropy $2$-groups and observing some basic properties in Section~\ref{sec:2disotropy}. 
In Section~\ref{sec:2disotropywithcoprods} we observe how the situation simplifies in the presence of binary coproducts,
 and in Section~\ref{sec:2ddensity} we simplify it further in the presence of a dense subcategory. 
 We then reap the payoff of these tools, by computing various isotropy $2$-groups of interest:
  \begin{itemize}
    \item For the $2$-category of groupoids, we show that the isotropy $2$-groups vanish in Theorem~\ref{thm:groupoidshaveno2disotropy}. 
    Thus groupoids have no nontrivial inner autoequivalences in our sense of the term. This should be compared with the results in~\cite{garner:innerofgroupoids}.
    \item For the $2$-category of indexed categories on $\CC$, \ie pseudofunctors $F\colon \CC\op\to\cat{Cat}$, we show that the isotropy $2$-group $\ZZ(F)$\/ is
    the $2$-group of pseudonatural autoequivalences of $\id[\CC]$. This generalizes the characterization of the covariant isotropy group of a presheaf topos
    in~\cite{hofstraetal:picard} (see also Example~\ref{ex:1d} for a statement and alternative proof of the result from loc.\ cit.)
    \item In Section~\ref{sec:monoidal} we compute the isotropy $2$-groups of monoidal categories, 
    showing that it behaves as one might guess based on the case of one-dimensional isotropy of strict monoidal 
    categories and strict monoidal functors. Specifically, we show that if one takes the $2$-category of monoidal categories, 
    strong monoidal functors and monoidal natural transformations, then the isotropy $2$-group $\ZZ(\CC)$\/ of a monoidal category 
    $\CC$ is equivalent to the \emph{Picard $2$-group of $\CC$}. 
    \item We then conclude by joining the previous to results together: for a pseudofunctor 
    $H\colon \CC\to\cat{MonCat}$, the isotropy $2$-group of $H$ is equivalent to the product of $\Aut(\id[\CC])$ and 
    the pseudolimit of the isotropy $2$-groups of $H(A)$ for each $A\in\CC$. 
    Consequently, we obtain a characterization of the inner autoequivalences of monoidal fibrations in the sense of~\cite{moeller:monoidalgrothendieck}.
  \end{itemize}

We point out that since we work with lax slices, pseudofunctors and pseudonatural transformations, 
our generalization of isotropy is genuinely bicategorical, rather than merely $\cat{Cat}$-enriched, 
and hence does not follow directly from an (hitherto undeveloped) enriched theory of isotropy. 
Despite this, we assume all of our $2$-categories to be strict. 
This is mostly for notational convenience and due to us not being aware of examples of interest 
where the underlying two-dimensional category fails to be strict---the basic definitions and results themselves generalize straightforwardly to bicategories.

\subsection*{Acknowledgements}
We wish to thank an anonymous reviewer whose detailed and constructive comments greatly improved the manuscript. Pieter Hofstra was supported by an NSERC discovery grant.

\section{One-dimensional warm-up}\label{sec:1d}

In this section we introduce a useful technique for computing isotropy. The $2$-categorical version will be crucial in later sections, 
but we first explain and illustrate it in the one-dimensional setting, since the technique is already interesting and useful there.
We first give an alternative characterization of $\ZZ$\/  in the situation where $\CC$ has binary coproducts; after that, we show how the computation of 
$\ZZ$\/ may be simplified by considering a dense subcategory of $\CC$.

We will generally denote the isotropy functor by $\ZZ$, except when we discuss several different categories and their isotropy functors simultaneously, at which point we will disambiguate by writing $\ZZ_\CC$, $\ZZ_\DD$ and so on. Observe that $\ZZ$ is not in general a functor from $\CC$ to the category of small groups, as the group of natural automorphisms of $P_X$ might be a proper class. For instance, take the category of sets and bijections between them, and adjoin freely a new initial object $0$. In the resulting category $\CC$, a natural automorphism of $0/\CC\to\CC$ can be uniquely specified by giving an automorphism of each cardinal, so that $\ZZ_\CC(0)$ is a large group.

\begin{definition} A category $\CC$ has \emph{small isotropy} if for each object $X$ of $\CC$, the class of natural automorphisms of the projection $X/\CC\to\CC$ is in fact a set, so that $\ZZ$ defines a functor $\CC\to\Grp$.
\end{definition}

Any small category has small isotropy. As we shall see, however (Theorem~\ref{thm:densitycontrolsisotropy}), it frequently happens that large categories still have small isotropy.

 Let us assume now that $\CC$ has binary coproducts. Then for every object $X$\/ of $\CC$ the projection $P_X\colon X/\CC\to\CC$ has a left adjoint 
 $L_X\colon \CC\to X/\CC$ that sends $A\in \CC$ to $i_X\colon X\to X+A$. 
 Taking mates induces an isomorphism $\Aut(P_X)\to \Aut(L_X)\op$ for each $X$, as taking mates reverses the order of composition. 
 When reinterpreting an element of isotropy as a natural automorphism $\alpha$\/ of $L_X$, we find that the component of $\alpha$\/ at an object $A$\/ is
an isomorphism $\alpha_A\colon X+A \to X+A$\/ that fixes $X$, in the sense that the diagram to the left in:
\[
\begin{tikzpicture}
\node (x) at (0,1) {$X$};
\node (a1) at (-1,-1) {$X+A$};
\node (a2) at (1,-1) {$X+A$};
\draw[->] (x) to node[left]  {$i_X$} (a1);
\draw[->] (x) to node[right]  {$i_X$} (a2);
\draw[->] (a1) to node[below]  {$\alpha_A$} (a2);
\end{tikzpicture}
\qquad
\begin{tikzpicture}
\node (a1) at (-1,1) {$X+A$};
\node (a2) at (1.5,1) {$X+A$};
\node (b1) at (-1,-1) {$X+B$};
\node (b2) at (1.5,-1) {$X+B$};
\draw[->] (a1) to node[above]  {$\alpha_A$} (a2);
\draw[->] (b1) to node[below]  {$\alpha_B$} (b2);
\draw[->] (a1) to node[left]  {$X+f$} (b1);
\draw[->] (a2) to node[right] {$X+f$} (b2);
\end{tikzpicture}
\]
commutes. Moreover, for any $f\colon A \to B$\/ we must have $\alpha_{B}(X+f)=(X+f)\alpha_A$, as in the diagram above right.

 Consequently, there is a unique way of extending the assignment $X\mapsto\Aut(L_X)\op$ into a functor $\LL$ that is naturally isomorphic to $\ZZ$ via these isomorphisms. 
 A straightforward calculation shows that $\LL$ is then defined on objects by  $\LL(X)=\Aut(L_X)\op$ and on morphisms by sending $f\colon X\to Y$ to the 
 function $\LL(f)\colon\Aut(L_X)\op\to\Aut(L_Y)\op$ which sends $t:=(t_A)_{A\in\CC}\in \Aut(L_X)\op$ to $\LL(f)(t)$ where $(\LL(f)t)_A\colon Y+A\to Y+A$ is defined by 
 commutativity of
\begin{equation}\label{eq:definingLonmorphisms}\vcenter{\hbox{\begin{tikzpicture}
  \node (a1) at (-1,1) {$A$};
  \node (a2) at (-1,-1) {$Y$};
   \node (b) at (1,0) {$Y+A$};
  \node (c) at (4,0) {$Y+A$};
  \node (d) at (1,1) {$X+A$};
  \node (e) at (4,1) {$X+A$};
  \draw[->] (a1) to node[below]  {$i_A$} (b);
  \draw[->] (a1) to node[above]  {$i_A$} (d);
  \draw[->] (d) to node[above]  {$t_A$} (e);
  \draw[->] (e) to node[right]  {$f+\id[A]$} (c);
  \draw[->] (a2) to node[above]  {$i_Y$} (b);
  \draw[->] (a2) to[in=-135,out=0] node[below]  {$i_Y$} (c);
 \draw[->] (b) to node[above] {$(\LL(f)t)_A$} (c);
\end{tikzpicture}}}
\end{equation}

Consider now a functor $F\colon\CC\to\DD$. We write $\hat{F}\colon \presh[D] \to \presh[C]$\/ for the functor that restricts along $F$.
Recall also that $F$\/ is called \emph{dense} if it satisfies any of the following equivalent conditions:

\begin{enumerate}[label=(\roman*)]
      \item for all cocontinuous functors  $G,H\colon \DD\to \EE$, restriction along $F$ induces a bijection $\nat(G,H)\to\nat(GF,HF)$;
      \item the restricted Yoneda embedding $\DD\xrightarrow{Y} \presh[D]\xrightarrow{\hat{F}}\presh[C]$ is fully faithful;
      \item for each $A,B\in\DD$, the function \[\DD(A,B)\to \presh[C](\DD(F-,A),\DD(F-,B))\] is a bijection;
      \item for every object $A$ of $\cat{D}$, the identity natural transformation $\DD(F-,A)\to \DD(F-,A)$ exhibits $A$ as the $\DD(F-,A)$-weighted colimit of $F$.
\end{enumerate}

That (i)-(iv) are equivalent is proven more generally for enriched categories \eg in~\cite[Chapter 5]{kelly:enriched}, and will serve as our template for bicategorical density in Section~\ref{sec:2ddensity}.

\begin{remark}Note that when $\CC$ and $\DD$ are large categories, the functor categories $\presh[C]$ and $\presh[D]$ are ``very large'' and might not even exist in some foundations. However, both the theorem and its proof could be rephrased so as to avoid such very large structures, at the cost of making them slightly more cumbersome to state. For instance, (iii) is essentially a rephrasing of (ii), and one could talk about the class of natural transformations $\DD(F(-),A)\to \DD(F(-),B))$ without invoking a very large ambient category in which it lives as a hom-class. Similar remarks apply to our treatment of density in the $2$-categorical case.
\end{remark}

We now show that the isotropy of a category is determined by its restriction along any dense functor. This result is similar to one in~\cite{funketal:crossedtoposes},
where it is shown that the contravariant isotropy functor of a presheaf topos sends colimits to limits and hence is determined by its values on the representables. 

\begin{theorem}\label{thm:densitycontrolsisotropy}
  Let $\cat{C}$ have binary coproducts and $K\colon \DD\to \CC$ be a dense functor. 
  Then restriction along $K$ defines an isomorphism $\LL\cong \LL_K$ where $\LL_K$ is defined 
  on objects by $\LL_K(X)=\Aut(L_X\circ K)\op$ and on morphisms as in~\eqref{eq:definingLonmorphisms} 
  with $A$ of the form $K(B)$ for $B\in \DD$. Thus $\ZZ\cong \LL_K$. 

  In particular, if there exists a dense $K\colon \DD\to\CC$ with small domain, then $\CC$ has small isotropy.  
  As a result, any locally presentable category has small isotropy. 
\end{theorem}

\begin{proof}
  As $L_X$ is a left adjoint and hence cocontinuous, density of $K$\/ implies that $\Aut(L_X)\op\cong \Aut(L_X\circ K)\op$; moreover, 
  this isomorphism is clearly natural in $X$. If we can choose a dense $K\colon \DD\to \CC$ with a small domain, then  
  $\Aut(L_X\circ K)$ and hence $\ZZ(X)$ is small for every $X$. 
  The last claim follows from the fact that any locally presentable category is cocomplete and has a small dense subcategory~\cite[Theorem 1.20, Example 1.24(i)]{adamekrosicky:locpres}.
\end{proof}

\begin{example}\label{ex:1d} Let us use Theorem~\ref{thm:densitycontrolsisotropy} to compute various isotropy groups of interest.
  \begin{itemize}
    \item
    
    For the category $\cat{Mon}$ of monoids, the full subcategory $V$\/ on the single object $F\{x,y\}$, 
    the free monoid on two generators, is a dense subcategory. 
    We can use this to re-cast the argument given in~\cite{bergman:inner} as follows. 
    
    For a monoid $M$, a natural automorphism of the composite
    \[ V \to \cat{Mon}\to M/\cat{Mon}\]
    consists of a single automorphism $\alpha$ of $i_M\colon M \to M+F\{x,y\}$ in $M/\cat{Mon}$ that is natural with respect to endomorphisms of $F\{x,y\}$. As such an automorphism must fix $M$, 
    it it is determined by a pair of elements $\alpha(x),\alpha(y)\in M+F\{x,y\}$. 
    These elements can be written as words $w_1(x),w_2(y)$ in elements of $M\sqcup\{x,y\}$. 
    Moreover, the naturality condition for the endomorphism of $F\{x,y\}$ determined by $x,y\mapsto x$ implies that $\alpha(x)$ is a word in $M\sqcup\{x\}$ 
    and similarly $\alpha(y)$ is a word in $M\sqcup\{y\}$. Naturality at the swap map $x\mapsto y$, $y\mapsto x$ 
    then implies that $w_1(x)=w_2(x)$, so $\alpha$ is determined by a single word $w(x)$ with $\alpha(x)=w(x)$ 
    and $\alpha(y)=w(y)$. 
    
    Next, naturality at $x,y\mapsto xy$ implies that $w(x)w(y)=\alpha(x)\alpha(y)=\alpha(xy)=w(xy)$. 
    In the word $w(x)w(y)$ each occurrence of $x$ is to the left of each occurrence of $y$, which also must hold in $w(xy)$. 
    Thus $w(x)$ contains $x$ at most once, and as $\alpha$ is an automorphism, it contains it exactly once. 
    Thus $\alpha(x)=w(x)=axb$ with $a,b\in M$. Now the equation $w(x)w(y)=w(xy)$ implies $ba=1$ in $M$. 
    On the other hand, naturality at $x,y\mapsto 1$ implies that $w(1)=1$ so that $ab=1$ and hence $b=a^{-1}$. 
    
    Conversely, it is easy to check that any invertible $a\in M$ determines such a natural automorphism, and hence that  
   $\ZZ_\cat{Mon}$ is naturally isomorphic to the functor 
    $M\mapsto\{\text{invertible elements of $M$}\}$, \ie the right adjoint of the inclusion $\Grp\to\cat{Mon}$.

    \item  Next, we consider a small category $\CC$, and compute the isotropy of $\presh[C]$. 
 Note that $y\colon \CC\to\presh[C]$ is dense. Given $F\in\presh[C]$, an element of $\alpha \in \ZZ(F)$ is then determined by a natural family of isomorphisms 
 $\alpha_A\colon F+yA\to F+yA$ fixing $F$. Since $\presh[C]$ is extensive, it follows that $\alpha_A$\/ must be of the form $\id[F]+\beta_A$, where $\beta_A\colon yA \to yA$\/ is
 an isomorphism. By Yoneda's lemma, we may identify $\beta_A$\/ with an automorphism of $A$; and since these automorphisms 
 are natural in $A$, they together form a natural automorphism of 
 the identity functor on $\CC$. (The group $\Aut(\id[\CC])$\/ is sometimes called the \emph{center} of $\cat{C}$.)
In conclusion,  $\ZZ$ on $\presh[C]$ is the constant functor with value $\Aut(\id[\CC])$. 

The above argument goes through for the category of sheaves on site $(\cat{C},J)$, 
provided that the topology $J$\/ is subcanonical. We therefore also recover~\cite[Proposition 3.12.]{parker:IsotropyofGrothendieckToposes}, which states
that the isotropy functor $\ZZ\colon \mathsf{Sh}(\cat{C},J) \to \Grp$\/ is constant with value $\Aut(\id[\CC])$, and which was first proved in loc.\ cit.\ by
syntactic means. However, both this result and \cite[Proposition 3.12.]{parker:IsotropyofGrothendieckToposes} are special cases of~\cite[Theorem 4.1]{parker:IsotropyofGrothendieckToposes} suggested by Richard Garner, which states that the covariant isotropy of any extensive category $\CC$ is constant at $\Aut(\id[\CC])$.
  \end{itemize}
\end{example}

While Theorem~\ref{thm:densitycontrolsisotropy} is phrased for dense functors, in most examples of interest one works with a full dense subcategory, so that the theorem above could be phrased in terms of dense families of generators. We will proceed to observe that (not necessarily dense) families of generators satisfy a weaker result, which is not quite quite sufficient to recover the isotropy group, but still helps in showing in particular instances that the isotropy trivializes. 

Let $\DD$ be a category with coproducts. Recall that a set $\{A_j\}_{j\in J}$ of objects of $\DD$ is called \emph{a family of generators} of $\DD$ if it satisfies any of the following equivalent conditions.

\begin{enumerate}[label=(\roman*)]
    \item For any parallel pair $f,g\colon A\to B$ of morphisms of $\DD$, if for all $j\in J$ and $h\colon A_i\to A$ we have $gh=fh$, then $g=h$;
    \item the restricted Yoneda embedding $\DD\xrightarrow{Y} \presh[D]\xrightarrow{\hat{K}}\presh[C]$ is fully faithful, where $K\colon\CC\to \DD$ is the inclusion of the full subcategory spanned by  $\{A_j\}_{j\in J}$;
    \item The canonical morphism \[f_A\colon \coprod_{j\in J, f\colon A_j\to A} A_j\longrightarrow A\] defined by $f_A\circ i_f=f$
    is an epimorphism.
\end{enumerate}

That these are equivalent is shown \eg in~\cite[Section 4.5]{borceux:vol1}. When the family is a singleton, we will just speak of a \emph{generator}. We will extend the list above by one further condition, and provide a proof as we have not seen this stated elsewhere. 

\begin{proposition}\label{prop:generators}
  Let $\DD$ be a category with coproducts, $\{A_j\}_{j\in J}$ a set of objects of $\DD$ and $K\colon \CC\to\DD$ the inclusion of the full subcategory spanned by $\{A_j\}_{j\in J}$. Then the following are equivalent
  \begin{enumerate}[label=(\alph*)]
   \item $\{A_j\}_{j\in J}$ is a family of generators of $\DD$.
   \item for all cocontinuous  $G,H\colon \DD\to \EE$, restriction along $K$ induces an injection $\nat(G,H)\to\nat(GK,HK)$.
\end{enumerate}
\end{proposition}

\begin{proof}
(b)$\Rightarrow$(a):
To prove (a), we will prove condition (ii) above. Now, any representable $\DD\op\to\cat{Set}$ is continuous, so that the corresponding opposite functor 
$\DD\to\cat{Set}\op$ is cocontinuous. As the Yoneda embedding is fully faithful,
 (a) implies that the composite $\DD\xrightarrow{Y} \presh[D]\xrightarrow{\hat{K}}\presh[C]$ is faithful.

(a)$\Rightarrow$(b):  Recall that $f$\/ is an epimorphism if and only if the diagram    
  \[\begin{tikzpicture}
    \matrix (m) [matrix of math nodes,row sep=2em,column sep=4em,minimum width=2em]
    {
     A &  B\\
      B &  B\\};
    \path[->]
    (m-1-1) edge node [left] {$f$} coordinate[midway] (a) (m-2-1)
           edge node [above] {$f$} (m-1-2)
    (m-1-2) edge node [right] {$\id$} (m-2-2)
    (m-2-1) edge node [below] {$\id$} coordinate[midway] (b) (m-2-2);
  \end{tikzpicture}\]
  is a pushout diagram. Consequently, if $G,H\colon \DD\to \EE$ are cocontinuous, they both preserve the coproduct 
  $C_A:=\coprod_{i\in I, f\colon A_i\to A} A_i$ of (iii) above, and the epicness of the canonical map $f_A\colon C_A\to A$. Consider natural transformations $\sigma,\tau\colon G\to H$ with $\sigma K=\tau K$. As the map $Gf_A\colon G (C_A)\to A$ is an epimorphism, to prove that $\sigma=\tau$ it suffices to show that $\sigma_{C_A}=\tau_{C_A}$ for each $A$. In turn, as $G$ preserves the coproduct structure of $C_A$, it suffices to show that $\sigma_{C_A}\circ G(i_f)=\tau_{C_A}\circ  G(i_f)$ for an arbitrary inclusion $i_f\colon A_j\to C_A$, but this follows from naturality of $\sigma$ and $\tau$ and the assumption $\sigma K=\tau K$. 
\end{proof}

\begin{theorem}\label{thm:generatorsandisotropy}
  Let $\DD$ have coproducts and $K\colon \CC\to \DD$ be the inclusion of a full subcategory spanned by a family of generators of $\DD$.
  Then there is a natural inclusion $\ZZ\hookrightarrow \LL_K$ where $\LL_K$ is defined as in Theorem~\ref{thm:densitycontrolsisotropy}.
\end{theorem}

\begin{proof} Given the isomorphism $\ZZ\cong\LL$, it suffices to show that restriction along $K$ defines an injection $\LL \to \LL_K$, 
but this follows from the cocontinuity of each $L_X$  and from condition (b) of Proposition~\ref{prop:generators}.
\end{proof}

In particular, if in the above $\LL_K$\/ is trivial, then so is $\ZZ$.

\begin{example} Let us see Theorem~\ref{thm:generatorsandisotropy} in action.
\begin{itemize} 
  \item Let $\cat{Grpd}$ be the $1$-category of (small) groupoids, and let $I$\/ be the ``walking isomorphism category'', \ie the indiscrete category on two objects. 
  Then  $I$ is a generator for $\cat{Grpd}$ as two functors between groupoids are equal if and only if they agree on morphisms. 
  For any groupoid $G$ the only automorphism of $G+I$ fixing $G$ has to send $I$ to itself. 
  Moreover, naturality with respect the constant maps $I\to I$ forces any such isomorphism to be the identity $I$, 
  so the only natural automorphism of $\{I\}\to \cat{Grpd}\xrightarrow{L_G} G/\cat{Grpd}$ is the identity. 
  Thus $\ZZ\colon\cat{Grpd}\to\Grp$ is constant at the trivial group. Again, this proof implicitly relies on $\cat{Grpd}$ being extensive, and as such is a special case of~\cite[Theorem 4.1]{parker:IsotropyofGrothendieckToposes}.
  \item The abelian group  $\mathbb{Z}$\/ is a generator for $\cat{Ab}$. 
  An automorphism  of $G+\mathbb{Z}=G\times\mathbb{Z}$ fixing $G$ is determined by the image 
  $(g,n)\in G\times\mathbb{Z}$ of $1\in\mathbb{Z}$.  As the inverse of this automorphism must also lie in 
  $G/\cat{Ab}$, we must have $g=0\in G$ and $n=\pm 1$. Consequently, if $K$ denotes the inclusion 
  $\{\mathbb{Z}\}\to\cat{Ab}$, then $\LL_K$ is the constant functor on $\mathbb{Z}_2$. 
  Furthermore, one can easily check that both $-\id$ and $\id$ define natural automorphisms of $G/\cat{Ab}\to\cat{Ab}$, so that $\ZZ_{\cat{Ab}}$ is also constant at $\mathbb{Z}_2$.
   \item We sketch how to use generating families to obtain the aforementioned result of~\cite{hofstraetal:picard} stating that for the $1$-category $\cat{StrMonCat}$\/ of strict monoidal categories and strict monoidal functors, the isotropy group of $\CC\in \cat{StrMonCat}$\/ is isomorphic to the strict Picard group of $\CC$. We will consider the generating family for $\cat{StrMonCat}$  given by
    \begin{itemize}
        \item the free monoidal category on two objects and
        \item the free monoidal category on $\bullet\rightarrow\bullet$. 
    \end{itemize}
    Let $\mathbf{F}$ be the full subcategory of $\cat{StrMonCat}$ spanned by these and $K\colon\mathbf{F}\to\cat{StrMonCat}$ be the inclusion. This is indeed a generating family  as two strict monoidal functors are equal iff they agree on morphisms. The reason to include the free monoidal subcategory on two generators in $\mathbf{F}$ as well is that the resulting composite \[\mathbf{F}\to\cat{StrMonCat}\xrightarrow{\ob}\cat{Mon}\] is dense. Consequently, given a monoidal category $\CC$, the diagram 
   \[
   \begin{tikzpicture}
     \node (a1) at (-1,1) {$\mathbf{F}$};
      \node (b) at (2,-1) {$\cat{Mon}$};
     \node (c) at (6,-1) {$\ob(\CC)/\cat{Mon}$};
     \node (d) at (2,1) {$\cat{StrMonCat}$};
     \node (e) at (6,1) {$\CC/\cat{StrMonCat}$};
     \draw[->] (a1) to node[below,sloped]  {$\ob\circ K$} (b);
     \draw[->] (a1) to node[above]  {$K$} (d);
     \draw[->] (d) to node[above]  {$L_\CC$} (e);
     \draw[->] (e) to node[right]  {$\ob$} (c);
    \draw[->] (b) to node[above] {$L_{\ob(\CC)}$} (c);
   \end{tikzpicture}
   \]
   commutes within isomorphism, inducing a homomorphism $\LL_K(\CC)\to\ZZ_\cat{Mon}(\ob(\CC))$. Furthermore, one can use the structure of the coproduct of $\CC$ with one of the free monoidal categories in question to show that this homomorphism $\LL_K(\CC)\to\ZZ_\cat{Mon}(\ob(\CC))$ is injective---ultimately, this follows from the fact that the relevant free monoidal categories are thin, so that whiskering an automorphism of $\mathbf{F}\xrightarrow{K}\cat{StrMonCat}\xrightarrow{L_\CC}\CC/\cat{StrMonCat}$ along $\CC/\cat{StrMonCat}\to \ob(\CC)/\cat{Mon}$ is injective. Precomposing this with the injection $\ZZ_{\cat{StrMonCat}}(\CC)\hookrightarrow \LL_K(\CC)$, we obtain an injection $\ZZ_{\cat{StrMonCat}}(\CC)\hookrightarrow \ZZ_\cat{Mon}(\ob(\CC))$ 
   from the isotropy group of $\CC$ to the strict Picard group of $\CC$. Since any strictly invertible object of $\CC$ induces an element of $\ZZ_\cat{StrMonCat}(\CC)$ via conjugation, this map is in fact an isomorphism.
\end{itemize}
\end{example}

To conclude this section on one-dimensional isotropy, we remark that much of what we have discussed can be generalized without much effort to 
the \emph{isotropy monoid} $\MM\colon \CC \to \Mon$, which is defined at an object $X$ of $\cat{C}$ by
\[ \MM(X)=\mathrm{Nat}(P_X,P_X), \]
the monoid of natural endomorphisms of the projection $P_X\colon X/\cat{C}\to \cat{C}$.
This generalization is already considered in Bergman~\cite{bergman:inner}, where it is shown that in several examples the monoid does not contain much interesting
information not already present in the isotropy group. However, it turns out that when generalizing to the two-dimensional setting 
it will often be convenient to establish results for $\MM$\/ and deduce the corresponding statements for $\ZZ$, using the fact that $\ZZ$\/ may now be re-expressed
as the composite
\[ \CC \to \Mon \to \Grp \]
where the functor $\Mon \to \Grp$\/ sends a monoid to its group of invertible elements.

\section{Two-dimensional background}\label{sec:2dbackground}

We refer the reader to~\cite{johnsonyau:2d} for the general definitions of a pseudofunctor, 
lax and pseudonatural transformations between pseudofunctors, and of a modification between lax transformations. 
However, we will spell out explicitly in the next section what lax and pseudonatural transformations and modifications 
between them amount to in the special case under consideration. 
For us, all $2$-categories are assumed to be strict, but this is mostly for convenience---the basic definitions and results themselves generalize 
straightforwardly to bicategories. For $2$-categories $\CC,\DD$, we will use $[\CC,\DD]$ to denote the $2$-category of pseudofunctors $\CC\to \DD$, 
pseudonatural transformations and modifications. We will use the fact that a pseudonatural transformation may be made into a 
pseudonatural equivalence (\ie an equivalence in $[\CC,\DD]$) if and only if each $1$-cell component of the pseudonatural transformation is an equivalence.

Recall that given a functor $F\colon\CC\to\DD$ between $1$-categories and an \emph{arbitrary} family of 
isomorphisms $\{\sigma_A\colon F(A)\to G(A)\}_{A\in\CC}$, there is a unique way of making $A\mapsto G(A)$ into 
a functor $G$ so that $\sigma$ gives an isomorphism $\sigma\colon F\to G$. We will repeatedly use the following two-dimensional analogues of this fact:
 
\begin{lemma}\label{lem:promote}
  \begin{itemize}
    \item Given a pseudofunctor $F\colon \CC\to\DD$, and arbitrary equivalence data \[\sigma_A\colon FA\leftrightarrows GA\colon \tau_A,\quad \Gamma_A\colon \id\xrightarrow{\cong} \sigma_A\tau_A,\quad \Sigma\colon \tau_A\sigma_A\xrightarrow{\cong} \id\] for each $A\in\CC$\/, there is an essentially unique\footnote{More precisely, any two such pseudofunctors are related by a unique invertible icon that is compatible with the chosen pseudonatural equivalences to $F$. Moreover, there is a canonical pseudofunctor $G$ that can be determined from the above data without making further arbitrary choices.} way of promoting $A\mapsto G(A)$ into a pseudofunctor $G$ so that  \[\sigma\colon F\leftrightarrows G\colon \tau,\quad \Gamma\colon \id\xrightarrow{\cong} \sigma\tau,\quad \Sigma\colon \tau\sigma\xrightarrow{\cong} \id\] defines a pseudonatural equivalence.

    \item Given a pseudonatural transformation $\sigma\colon F\to G$, and arbitrary invertible $2$-cells $\Gamma_A\colon \sigma_A\to\tau_A$ for each $A\in\CC$, there is a unique pseudonatural transformation $\tau$ whose $1$-cell components are given by $\{\tau_A\}_{A\in\CC}$ for which the family $\{\Gamma_A\}_{A\in\CC}$ defines a modification $\Gamma\colon\sigma\to\tau$. 
  \end{itemize}
\end{lemma}

\begin{definition} 
  Let $\CC$ be a $2$-category and $X$ an object of of $\CC$. The lax slice $\lslice{X}{C}$ of $\CC$ under $X$ is defined as follows. 
  The objects of  $\lslice{X}{C}$ are pairs $(A,f)$, where $A$ is an object and $f\colon X\to A$ is a $1$-cell of $\CC$. 
  Given two objects $(A,f)$, $(B,g)$ of  $\lslice{X}{C}$, a $1$-cell $f\to g$ is depicted in the triangle 
    \[\begin{tikzpicture}
    \node (x) at (0,1) {$X$};
    \node (a) at (-1,-1) {$A$};
    \node (b) at (1,-1) {$B$};
    \draw[->] (x) to coordinate[midway] (f) node[left] {$f$} (a);
    \draw[->] (x) to  coordinate[midway] (g)  node[right] {$g$} (b);
    \draw[->] (a) to node[below] {$j$} (b);
    \draw[shorttwocell] (f) to node[above] {$\sigma$} (g);
    \end{tikzpicture}\]
  and is given by a pair $(j,\sigma)$, where $j\colon A\to B$ is a $1$-cell and $\sigma\colon hf\to g$ is a $2$-cell of $\CC$. 
  A $2$-cell $(j,\sigma)\to (k,\tau)$ is given by a $2$-cell $\theta\colon j\to k$ in $\CC$ satisfying the equality of pasting diagrams
    \[\begin{tikzpicture}[baseline=-.3cm]
    \node (x) at (0,1) {$X$};
    \node (a) at (-1,-1) {$A$};
    \node (b) at (1,-1) {$B$};
    \draw[->] (x) to coordinate[midway] (f) node[left] {$f$} (a);
    \draw[->] (x) to coordinate[midway] (g) node[right] {$g$} (b);
    \draw[->] (a) to[out=35,in=145] coordinate[midway] (h1) node[above]{$k$} (b);
    \draw[->] (a) to[out=-35,in=215] coordinate[midway] (h2) node[below] {$j$} (b);
    \draw[shorttwocell] (f) to node[above] {$\tau$} (g);
    \draw[shorttwocell] (h2) to node[right] {$\theta$} (h1);
    \end{tikzpicture}
    =\begin{tikzpicture}[baseline=-.3cm]
    \node (x) at (0,1) {$X$};
    \node (a) at (-1,-1) {$A$};
    \node (b) at (1,-1) {$B$};
    \draw[->] (x) to coordinate[midway] (f) node[left] {$f$} (a);
    \draw[->] (x) to  coordinate[midway] (g) node[right] {$g$} (b);
    \draw[->] (a) to node[below] {$j$} (b);
    \draw[shorttwocell] (f) to node[above] {$\sigma$} (g);
    \end{tikzpicture}\]
\ie  $\tau \circ(\theta f)=\sigma$.

The pseudo-slice $\pslice{X}{C}$ is the locally full sub-$2$-category of $\lslice{X}{C}$ with the same objects,
 but whose $1$-cells are those $1$-cells $(h,\sigma)$ of $\lslice{X}{C}$\/ with $\sigma$\/ invertible in $\CC$.
\end{definition}

There is a strict $2$-functor $Q_X\colon\lslice{X}{C}\to \CC$ that sends a $2$-cell $\theta\colon (j,\sigma)\to (k,\tau)$  of  $\lslice{X}{C}$ to $\theta\colon j\to k$. 
We denote the (strict) inclusion $2$-functor $\pslice{X}{C}\hookrightarrow\lslice{X}{C}$ by $I_X$ and the composite $Q_XI_X$ by $P_X$. 

For us, the word \emph{$2$-group} refers to a monoidal groupoid in which every object $X$\/ admits a weak inverse, 
\ie an object $X^{-1}$ such that $X\otimes X^{-1}\cong I\cong X\otimes X^{-1}$. 
A morphism of $2$-groups is a strong monoidal functor, and a $2$-cell between them is given by a monoidal 
natural transformation (which is necessarily an isomorphism). We denote the $2$-category of (small) $2$-groups by $\cat{2Grp}$. 
If $\CC$\/ is a $2$-category and $A$ is an object of $\CC$, we denote by $\Aut(A)$ the $2$-group of autoequivalences of 
$A$\/ and isomorphisms between them, where the tensor product is composition. 
In particular, we will assume that the ambient $2$-category of the identity 
functor $\id[\CC]\colon\CC\to\CC$ is $[\CC,\CC]$, so that $\Aut(\id[\CC])$\/ denotes the $2$-group of 
pseudonatural autoequivalences of $\id[\CC]$ and of invertible modifications between them.

We will now discuss some issues we face when dealing with lax natural transformations. 
First of all, whiskering a lax natural transformation along a pseudofunctor on either side 
results in a well-defined lax natural transformation~\cite[Section 11.1]{johnsonyau:2d}. 
In fact, postwhiskering along $F\colon \CC\to \DD$ defines a pseudofunctor $[\EE,\CC]_{l}\to [\EE,\DD]_l$, where $[\CC,\DD]_l$ is the $2$-category of pseudofunctors $\CC\to \DD$, 
lax natural transformations and modifications, and this pseudofunctor is strict whenever $F$ is. Moreover, prewhiskering along $F\colon \CC\to \DD$ defines a strict $2$-functor $[\DD,\EE]_{l}\to [\CC,\EE]_l$. However, it is well known that that lax natural transformations do not form the $2$-cells of a bicategory, or even of a tricategory, 
 see \eg the discussion in~\cite{lack:icons}. In addition to other problems that disappear when restricting to strict $2$-functors between strict $2$-categories, 
   one has an issue with interchange: while whiskering lax natural transformations between $2$-functors is a well-defined operation, the meaning of 
\begin{equation}\label{eq:whisker}\vcenter{\hbox{
\begin{tikzpicture}
    \node (a) at (-1,-1) {$\CC$};
    \node (b) at (1,-1) {$\DD$};
    \node (c) at (3,-1) {$\EE$};
    \draw[->] (a) to[out=35,in=145] coordinate[midway] (h1) node[above]{$F$} (b);
    \draw[->] (a) to[out=-35,in=215] coordinate[midway] (h2) node[below] {$G$} (b);
    \draw[shorttwocell] (h1) to node[right] {$\sigma$} (h2);
    \draw[->] (b) to[out=35,in=145] coordinate[midway] (g1) node[above]{$H$} (c);
    \draw[->] (b) to[out=-35,in=215] coordinate[midway] (g2) node[below] {$K$} (c);
    \draw[shorttwocell] (g1) to node[right] {$\tau$} (g2);
  \end{tikzpicture}}}\end{equation}
is ambiguous in general when $\sigma$ and $\tau$ are lax natural, as the two interpretations of the picture---namely,  
$(K\sigma)\circ(\tau F)$ and  $(\tau G)\circ(H\sigma)$---do not agree. The best we can do in this situation is to use $\tau$\/ in order to obtain $2$-cells 
      \begin{equation}\label{eq:interchange}\vcenter{\hbox{
      \begin{tikzpicture}
    \matrix (m) [matrix of math nodes,row sep=2em,column sep=4em,minimum width=2em]
    {
     HF(A) & HG(A)\\
      KF(A) & KG(A) \\};
    \path[->]
    (m-1-1) edge node [left] {$\tau_{F(A)}$} coordinate[midway] (a) (m-2-1)
           edge node [above] {$H(\sigma_A)$} (m-1-2)
    (m-1-2) edge node [right] {$\tau_{G(A)}$} (m-2-2)
    (m-2-1) edge node [below] {$K(\sigma_A)$} coordinate[midway] (b) (m-2-2);
    \draw[shorttwocell](m-2-1) to node [above,sloped] {$\Psi[\sigma,\tau]_A$} node [below,sloped] {$\tau_{\sigma_A}$} (m-1-2);
  \end{tikzpicture}}}\end{equation}
resulting in a modification  $\Psi=\Psi[\sigma,\tau]\colon (K\sigma)\circ(\tau F)\to (\tau G)\circ(H\sigma)$ that is guaranteed to be invertible 
only when each $\tau_{\sigma_A}$ is invertible (\eg if $\tau$ is pseudonatural). The resulting 3-dimensional structure is roughly speaking a 
``tricategory with weak interchange'', and can be understood more formally in terms of the canonical self-enrichment stemming from the 
(lax) Gray tensor product~\cite[Section I.4]{gray:formalct}, or if one wants to allow for pseudofunctors between strict $2$-categories, 
in terms of categories ``weakly enriched in Gray'', see \eg the formalization in~\cite[Section 1.3]{verity:thesis}. 
However, for our purposes it is sufficient to work in a more pedestrian manner: when all lax natural transformations in sight 
are in fact pseudonatural, we are working in a tricategory, and when on occasion we need a lax natural transformation, 
we take care to only use invertibility of the (components of the) 
interchanger $\Psi[\sigma, \tau]_A$\/ in instances where $\tau_{\sigma_A}$\/ has been shown to be invertible. 
On occasion, we will also use this weak interchanger, so for future reference we will record a few trivial observations on it now: 

\begin{lemma}\label{lem:inter}
Consider lax natural transformations $\sigma, \tau$\/ as in~\eqref{eq:whisker}, and the interchanger $\Psi[\sigma,\tau]$\/ as defined by~\eqref{eq:interchange}.
  \begin{enumerate}
    \item When $\tau$\/ is the identity, the resulting modification $\Psi[\sigma,\tau]$\/ is the identity modification; more generally, when $\tau$\/ is invertible,
    then so is $\Psi[\sigma,\tau]$.
    \item When $\tau=\tau_1\circ \tau_2$\/ is a composite of two lax natural transformations, then
    $\Psi[\sigma,\tau]$\/ is the composite of the pasting diagram:
    \[
       \begin{tikzpicture}
    \matrix (m) [matrix of math nodes,row sep=2.5em,column sep=4em,minimum width=2em]
    {
     HF & HG\\
      LF & LG \\
      KF & KG \\
      };
    \path[->]
    (m-1-1) edge node [left] {$\tau_1.F$} coordinate[midway] (a) (m-2-1)
           edge node [above] {$H.\sigma$} (m-1-2)
    (m-1-2) edge node [right] {$\tau_1.G$} (m-2-2)
    (m-2-1) edge node [above] {$L.\sigma$} coordinate[midway] (b) (m-2-2);
      \path[->]
    (m-2-1) edge node [left] {$\tau_2.F$} coordinate[midway] (a) (m-3-1)
      (m-3-1)     edge node [below] {$K.\sigma$} (m-3-2)
      (m-2-2) edge node[right] {$\tau_1.G$} (m-3-2);
    \draw[shorttwocell](m-2-1) to node [above,sloped] {$\Psi[\sigma,\tau_1]$} (m-1-2);
     \draw[shorttwocell](m-3-1) to node [above,sloped] {$\Psi[\sigma,\tau_2]$} (m-2-2);
  \end{tikzpicture}
    \] 
    \item For a modification $\Gamma\colon \tau \to \tau'$, we have, for each object $A$, an equality of pasting diagrams
    \begin{align*}
      \begin{tikzpicture}[baseline=-.35cm]
    \matrix (m) [matrix of math nodes,row sep=2em,column sep=4em,minimum width=2em,ampersand replacement=\&]
    {
    HF(A) \& HF(A) \& HG(A)\\
    KF(A) \&  KF(A) \& KG(A) \\};
    \draw[=] (m-1-1) to (m-1-2);
    \draw[=] (m-2-1) to (m-2-2);
     \path[->]
    (m-1-1) edge node [left] {$\tau_{FA}$} coordinate[midway] (a) (m-2-1);
     \draw[shorttwocell](m-2-1) to node [above,sloped] {$\Gamma_{FA}$}  (m-1-2);
    \path[->]
    (m-1-2) edge node [left] {$\tau'_{FA}$} coordinate[midway] (a) (m-2-2)
           edge node [above] {$H(\sigma_A)$} (m-1-3)
    (m-1-3) edge node [right] {$\tau'_{GA}$} (m-2-3)
    (m-2-2) edge node [below] {$K(\sigma_A)$} coordinate[midway] (b) (m-2-3);
    \draw[shorttwocell](m-2-2) to node [above,sloped] {$\Psi[\sigma,\tau']_A$}  (m-1-3);
  \end{tikzpicture}  \\  = 
        \begin{tikzpicture}[baseline=-.3cm]
    \matrix (m) [matrix of math nodes,row sep=2em,column sep=4em,minimum width=2em,ampersand replacement=\&]
    {
HF(A) \& HG(A) \& HG(A) \\
    KF(A) \& KG(A) \& KG(A) \\};
    \draw[=] (m-1-2) to (m-1-3);
    \draw[=] (m-2-2) to (m-2-3);
     \path[->]
    (m-1-3) edge node [left] {$\tau'_{GA}$} coordinate[midway] (a) (m-2-3);
     \draw[shorttwocell](m-2-2) to node [above,sloped] {$\Gamma_{GA}$}  (m-1-3);
    \path[->]
    (m-1-1) edge node [left] {$\tau_{FA}$} coordinate[midway] (a) (m-2-1)
           edge node [above] {$H(\sigma_A)$} (m-1-2)
    (m-1-2) edge node [right] {$\tau_{GA}$} (m-2-2)
    (m-2-1) edge node [below] {$K(\sigma_A)$} coordinate[midway] (b) (m-2-2);
    \draw[shorttwocell](m-2-1) to node [above,sloped] {$\Psi[\sigma,\tau]_A$} node [below,sloped] {$\tau_{\sigma_A}$} (m-1-2);
  \end{tikzpicture}
    \end{align*}
  \end{enumerate}
\end{lemma}

Lastly, we will use string diagrams~\cite{marsden:stringdiagrams,selinger:graphicallanguages} to aid our reasoning, both when reasoning inside a $2$-category or a monoidal category, 
and when reasoning about pseudonatural transformations between $2$-functors. 
Strictly speaking, in the latter case one might want to use  three-dimensional surface diagrams,
 but for our purposes the two-dimensional string diagram calculus is sufficient: the main difference to the usual string diagrams is 
 that instead of equations between string diagrams we get isomorphisms between them, representing invertible modifications. 
 As we do not need to reason about equations between the resulting modifications, this two-dimensional reasoning 
 strikes a good balance between readability and rigor. 
 
 Our convention is to write string diagrams bottom to top, \eg 
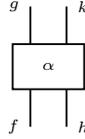
   \[
  \begin{pic}
    \setlength\minimummorphismwidth{6mm}
    \node[morphism] (f) at (0,0) {$\alpha$};
    \draw (f.south west) to +(0,-.5) node[left] {$f$};
    \draw (f.north west) to +(0,.5) node[left] {$g$};
    \draw (f.south east) to +(0,-.5) node[right] {$h$};
    \draw (f.north east) to +(0,.5) node[right] {$k$};
  \end{pic}\]
represents a $2$-cell $\alpha\colon h \circ f \Rightarrow k \circ g$. 
 
\section{Two-dimensional isotropy}\label{sec:2disotropy}

We now aim to define, for a $2$-category $\CC$, the isotropy $2$-group of $\CC$. It will be convenient to also introduce the notion of isotropy $2$-monoid; in fact, once we define
$\MM\colon \CC \to \cat{MonCat}$, the isotropy $2$-group $\ZZ\colon \CC \to \cat{2Grp}$\/ 
is then simply the composite of $\MM$\/ with the $2$-functor $\cat{MonCat} \to \cat{2Grp}$\/ that sends a monoidal category
to the subcategory of weakly invertible objects. Towards defining $\MM$\/ and $\ZZ$, we first introduce an approximation to these concepts as follows:

\begin{definition}
Given a $2$-category $\CC$\/ and an object $X$\/ of $\CC$, let $\Mps(X)$\/ be the monoidal category whose 
 objects are pseudonatural endomorphisms of $P_X\colon \pslice{X}{C}\to\CC$, whose 
 morphisms are given by modifications, and where the monoidal product is given by composition. We also let $\Zps(X)$\/ be the subcategory 
 of $\Mps$\/ on the pseudonatural autoequivalences of $P_X$\/ and the invertible modifications between them.
 \end{definition}
 
 Unfortunately, this definition does not quite work in general: 
 as the following example shows, the assignment $X \mapsto \Zps(X)$\/ (and similarly $X \mapsto \Mps(X)$\/ )
 need not be $2$-functorial.
 \begin{example}\label{ex:pseudoisotropyisnot2functorial}
  Let $\CC$ be the $2$-category generated by 
  \begin{itemize} 
    \item Two objects $X$ and $Y$
    \item Two $1$-cells $f,g\colon X\to Y$ and two $1$-cells $\alpha_f,\alpha_g\colon Y\to Y$
    \item A $2$-cell $\theta\colon f\to g$
  \end{itemize}
  subject to the constraints
  \begin{itemize}
    \item $\alpha_f,\alpha_g$ are isomorphisms
    \item $\alpha_f f=f$ and $\alpha_g g=g$
  \end{itemize}
  A pictorial representation of the generators of $\CC$ is given by
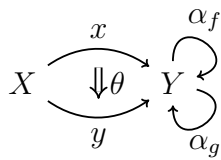
\[\begin{tikzpicture}
    \node (a) at (-1,-1) {$X$};
    \node (b) at (1,-1) {$Y$};
    \draw[->] (a) to[out=35,in=145] coordinate[midway] (h1) node[above]{$f$} (b);
    \draw[->] (a) to[out=-35,in=215] coordinate[midway] (h2) node[below] {$g$} (b);
    \draw[shorttwocell] (h1) to node[right] {$\theta$} (h2);
    \draw[->] (b) to[loop above,out=90,in=15,looseness=5] node[above] {$\alpha_f$} (b);
    \draw[->] (b) to[loop below,out=-15,in=-90,looseness=5] node[below] {$\alpha_g$} (b);
    \end{tikzpicture}\] 
  There is a $2$-natural automorphism $\alpha$ of $\pslice{X}{C}\to\CC$ whose components at 
  $\id[X]$, $f$ and $g$ are given by $\id[X]$, $\alpha_f$ and $\alpha_g$, respectively, with the remaining 
  components determined uniquely by these choices. However, there is no $2$-cell $\alpha_f\to \alpha_g$, 
  and consequently no modification from $\alpha_{-f}$ to $\alpha_{-g}$. This shows that in general one cannot get a $2$-functor sending $X$ to $\Zps(X)$.
\end{example}

We thus seek to modify the above definition in such a way that we obtain well-defined $2$-functors. The solution is to 
consider the lax slice, together with a slightly wider class of transformations, defined as follows.

\begin{definition}  Given a $2$-category $\CC$ and an object $X$ of $\CC$, we will call a lax natural endomorphism of
 $Q_X\colon\lslice{X}{C}\to\CC$ \emph{almost pseudonatural} if its restriction along $I_X\colon\pslice{X}{C}\to\lslice{X}{C}$ is pseudonatural.
\end{definition}

It is straightforward to show that the composite of almost pseudonatural transformations is again almost pseudonatural: if $\alpha, \beta$\/ are such that
$\alpha.I_X$\/ and $\beta.I_X$\/ are pseudonatural, then so is the composite $(\alpha \circ \beta).I_X=(\alpha.I_X) \circ (\beta.I_X)$.

\begin{example}\label{ex:constants}
   Assume $\CC$ has a pseudo-terminal object $1$. Then any $1$-cell $x\colon 1\to X$ in $\CC$ induces an endomorphism of 
   $Q_X\colon \lslice{X}{C}\to \CC$ whose $1$-cell component at $f\colon X\to A$ is given by 
   $A\xrightarrow{!}1\xrightarrow{x} X\xrightarrow{f} A$ and whose $2$-cell component at $(j,\sigma)\colon (A,f)\to (B,g)$ is given by
     \[\begin{tikzpicture}
        \matrix (m) [matrix of math nodes,row sep=2em,column sep=4em,minimum width=2em]
        {
         A & 1 & X & A\\
          B &  1 & X & B\\};
        \path[->]
        (m-1-1) edge node [left] {$j$} coordinate[midway] (a) (m-2-1)
               edge node [above] {$!$} (m-1-2)
        (m-1-2) edge node [right] {$\id$} (m-2-2)
                edge node [above] {$x$} (m-1-3)
        (m-2-2) edge node [below] {$x$} (m-2-3)
        (m-2-1) edge node [below] {$!$} coordinate[midway] (b) (m-2-2)
        (m-1-3) edge node[above] {$f$} (m-1-4)
                edge node [right] {$\id$} (m-2-3)
        (m-2-3) edge node[below] {$g$} (m-2-4)
        (m-1-4) edge node[right] {$j$} (m-2-4);
        \draw[shorttwocell](m-1-4) to node [above,sloped] {$\sigma$} (m-2-3); 
  \end{tikzpicture}\]
  where the leftmost square is to be filled with the canonical isomorphism. 
  Note that this composite $2$-cell is invertible whenever $j$ is, so that this endomorphism is almost pseudonatural. 
  In the case where $\CC$\/ is the $2$-category $\cat{MonCat_l}$ of monoidal categories, lax monoidal functors and monoidal natural transformations, 
  such endomorphisms of $P_\CC\colon \lslice{\CC}{MonCat_l}\to \cat{MonCat_l}$\/ thus correspond to monoids in $\CC$. Moreover, 
  modifications between two such endomorphisms correspond to monoid homomorphisms.
\end{example}

With this class of transformations we can now give the final definition of the isotropy $2$-monoid and isotropy $2$-group:

\begin{definition}
 Given a $2$-category $\cat{C}$\/ and an object $X$\/ of $\CC$, the  
  \begin{enumerate}
      \item  \emph{isotropy $2$-monoid} $\MM(X)$\/ at $X$\/ is the monoidal category 
      whose objects are almost pseudonatural endomorphisms of $Q_X\colon\lslice{X}{C}\to \CC$, 
      whose morphisms are modifications, and where the monoidal product is given by composition; 
      \item \emph{isotropy $2$-group} $\ZZ(X)$\/ at $X$\/ is the $2$-group whose objects are almost pseudonatural autoequivalences of  
      $Q_X\colon\lslice{X}{C}\to \CC$, whose morphisms are invertible modifications, and where the monoidal product is given by composition.
     \end{enumerate}
\end{definition}

In Section~\ref{sec:2disotropywithcoprods} we will prove that under suitable assumptions on $\CC$, we in fact have that $\MM \cong \Mps$\/ and $\ZZ \cong \Zps$.

In the remainder of this section we demonstrate that for general $\CC$, $\MM$\/ and $\ZZ$\/ are well-defined $2$-functors. 
 It suffices to show this for $X\mapsto \MM(X)$, 
as $\ZZ$\/ can be obtained from $\MM$\/ by postcomposing with the $2$-functor $\cat{MonCat}\to\cat{2Grp}$\/ 
that sends a monoidal category to its $2$-group of weakly invertible objects.  
We will deduce the functoriality of $X\mapsto \MM(X)$\/ formally, but we will interleave this with an explicit discussion of how 
$\MM$\/ acts on objects, morphisms and $2$-cells of $\CC$. In fact, the same argument will establish functoriality of the assignment sending $X$ to all lax natural endomorphisms of $Q_X$ and not just to the almost pseudonatural ones. However, as for our purposes $\MM(X)$ already is merely a convenient ambient structure, we do not consider this variant of the isotropy $2$-monoid further.

To begin, let us spell out the meaning of a lax natural endomorphism of 
$Q_X\colon\lslice{X}{C}\to \CC$: the data for such an automorphism $\alpha$\/ consists of
  \begin{itemize}
    \item for each object $(A,f)$ of $\lslice{X}{C}$, a $1$-cell $\alpha_{(A,f)}\colon A\to A$ in $\CC$
    \item for each $1$-cell $(j,\sigma)\colon (A,f)\to (B,g)$ of $\lslice{X}{C}$, a $2$-cell $\alpha_{(j,\sigma)}\colon j\alpha_{(A,f)}\to \alpha_{(B,g)}j$ in $\CC$, depicted by the pasting diagram
      \[\begin{tikzpicture}
    \matrix (m) [matrix of math nodes,row sep=2em,column sep=4em,minimum width=2em]
    {
     A & B\\
      A & B \\};
    \path[->]
    (m-1-1) edge node [left] {$\alpha_{(A,f)}$} coordinate[midway] (a) (m-2-1)
           edge node [above] {$j$} (m-1-2)
    (m-1-2) edge node [right] {$\alpha_{(B,g)}$} (m-2-2)
    (m-2-1) edge node [below] {$j$} coordinate[midway] (b) (m-2-2);
    \draw[shorttwocell](m-2-1) to node [above,sloped] {$\alpha_{(j,\sigma)}$} (m-1-2);
  \end{tikzpicture}\]
  or the string diagram 
  \[
  \begin{pic}
    \node[morphism] (f) at (0,0) {$\alpha_{(j,\sigma)}$};
    \draw (f.south west) to +(0,-.5) node[left] {$\alpha_{(A,f)}$};
    \draw (f.north west) to +(0,.5) node[left] {$j$};
    \draw (f.south east) to +(0,-.5) node[right] {$j$};
    \draw (f.north east) to +(0,.5) node[right] {$\alpha_{(B,g)}$};
  \end{pic}\]

  \end{itemize}
This data should satisfy the following conditions
\begin{itemize}
  \item ($2$-naturality): for any $2$-cell $\theta\colon (j,\sigma)\to (k,\tau)$ of $\lslice{X}{C}$ we have 
        \[
        \begin{pic}
          \node[morphism] (f) at (0,0) {$\alpha_{(j,\sigma)}$};
          \node[morphism,scale=0.5,font=\normalsize] (t) at (-.28,.75) {$\theta$};
          \draw (t.north) to +(0,.3) node[left] {$k$};
          \draw (f.south west) to +(0,-.5) node[left] {$\alpha_{(A,f)}$};
          \draw (f.north west) to  node[left] {$j$} (t.south);
          \draw (f.south east) to +(0,-.5) node[right] {$j$};
          \draw (f.north east) to node[right] {$\alpha_{(B,g)}$} ++(0,.2) to ++(0,.7) ;
        \end{pic}=  \begin{pic}
          \node[morphism] (f) at (0,0) {$\alpha_{(k,\tau)}$};
          \node[morphism,scale=0.5,font=\normalsize] (t) at (.29,-.75) {$\theta$};
            \draw (t.south) to +(0,-.3) node[right] {$j$};
          \draw (f.south west) to  ++(0,-.3) node[left] {$\alpha_{(A,f)}$} to ++(0,-.6);
          \draw (f.north west) to +(0,.5) node[left] {$k$};
          \draw (f.south east) to  node[right] {$k$} (t.north);
          \draw (f.north east) to +(0,.5) node[right] {$\alpha_{(B,g)}$};
        \end{pic}
        \]
  \item (unit constraint): for any identity $1$-cell $(\id[A],\id[f])\colon (A,f)\to (A,f)$ of $\lslice{X}{C}$ we have 
        \[
    \begin{pic}
    \node[morphism] (f) at (0,0) {$\alpha_{(\id[A],\id[f])}$};
    \draw (f.south) to +(0,-.5) node[left] {$\alpha_{(A,f)}$};
    \draw (f.north) to +(0,.5) node[right] {$\alpha_{(A,f)}$};
    \end{pic}=\begin{pic} \draw (0,-.75) node[right] {$\alpha_{(A,f)}$} to (0,.75);\end{pic}\]
  \item (respect for composites): for any $(j,\sigma)\colon (A,f)\to (B,g)$ and $(k,\tau)\colon (B,g)\to (C,h)$ we have
        \[\begin{pic}
            \node[morphism] (f) at (0,0) {$\alpha_{((k,\tau)\circ(j,\sigma))}$};
            \draw (f.south west) to +(0,-.5) node[left] {$\alpha_{(A,f)}$};
            \draw (f.north west) to +(0,.5) node[left] {$j$};
            \draw (f.north) to +(0,.5) node[left] {$k$};
            \draw (f.south) to +(0,-.5) node[right] {$j$};
            \draw (f.south east) to +(0,-.5) node[right] {$k$};
            \draw (f.north east) to +(0,.5) node[right] {$\alpha_{(C,h)}$};
          \end{pic}=
          \begin{pic}
            \node[morphism] (f) at (0,0) {$\alpha_{(j,\sigma)}$};
            \node[morphism] (g) at (.57,1) {$\alpha_{(k,\tau)}$};
            \draw (f.south west) to +(0,-.5) node[left] {$\alpha_{(A,f)}$};
            \draw (f.north west) to +(0,1.5) node[left] {$j$};
            \draw (f.south east) to +(0,-.5) node[right] {$j$};
            \draw (f.north east) to (g.south west); 
            \draw (g.north west) to +(0,.5) node[left] {$k$};
            \draw (g.south east) to +(0,-1.5) node[right] {$k$};
            \draw (g.north east) to +(0,.5) node[right] {$\alpha_{(C,h)}$};
          \end{pic}
          \]
\end{itemize}
A lax natural endomorphism $\alpha$ of $Q_X\colon\lslice{X}{C}\to \CC$ is almost pseudonatural if $\alpha_{[j,\sigma]}$ is invertible whenever $\sigma$ is invertible.

A modification $\Gamma\colon\alpha\to\beta$ between two such lax natural transformations $\alpha,\beta$\/ 
 consists of a $2$-cell $\Gamma_{(A,f)}\colon \alpha_{(A,f)}\to \beta_{(A,f)}$\/ in $\CC$ for each object $(A,f)$ of $\lslice{X}{C}$ such that 
  \[
  \begin{pic}
    \node[morphism] (f) at (0,0) {$\beta_{(j,\sigma)}$};
    \node[morphism,scale=0.5,font=\normalsize] (t) at (-.275,-.75) {$\Gamma_{(A,f)}$};
    \draw (f.south west) to  node[left] {$\beta_{(A,f)}$} (t.north);
    \draw (t.south) to +(0,-.3) node[left] {$\alpha_{(A,f)}$};
    \draw (f.north west) to +(0,.5) node[left] {$j$};
    \draw (f.south east) to ++(0,-.2) node[right] {$j$} to ++(0,-.7) ;
    \draw (f.north east) to +(0,.5) node[right] {$\beta_{(B,g)}$};
  \end{pic}=
    \begin{pic}
    \node[morphism] (f) at (0,0) {$\alpha_{(j,\sigma)}$};
     \node[morphism,scale=0.5,font=\normalsize] (t) at (.285,.75) {$\Gamma_{(B,g)}$};
    \draw (f.south west) to +(0,-.5) node[left] {$\alpha_{(A,f)}$};
    \draw (f.north west) to ++(0,.2) node[left] {$j$}  to ++(0,.7);
    \draw (f.south east) to +(0,-.5) node[right] {$j$};
    \draw (f.north east) to node[right] {$\alpha_{(B,g)}$} (t.south);
    \draw (t.north) to +(0,.3) node[right] {$\beta_{(B,g)}$};
  \end{pic}\]

This equips $\MM(X)$ with the structure of a category. This category is in fact a (strict) monoidal category, as a result of the 
general fact that for any two $2$-categories $\CC,\DD$ there is a strict $2$-category $[\CC,\DD]_{l}$ of pseudofunctors $\CC\to \DD$, 
lax natural transformations and modifications, so that for any functor $F\colon \CC\to \DD$ the category $[\CC,\DD]_{l}(F,F)$ is strictly monoidal. 
Moreover, it is easy to see that $\MM(X)$\/ is indeed a monoidal subcategory of the category of all lax automorphisms of $Q_X\colon\lslice{X}{C}\to \CC$. 

We next describe the action of $\MM$ on morphisms and $2$-cells. 
Given $x\colon X\to Y$ in $\CC$, we now define $\MM(x)\colon \MM(X)\to\MM(Y)$. 
Conceptually, this can be defined by observing that $x$ induces  $2$-functors $x^*\colon\lslice{Y}{C}\to\lslice{X}{C}$ and $x^{\bullet}\colon\pslice{Y}{C}\to\pslice{X}{C}$ fitting into a strictly commuting diagram
    \[\begin{tikzpicture}
    \node (x) at (0,1) {$\lslice{Y}{C}$};
    \node (a) at (-1,-1) {$\lslice{X}{C}$};
    \node (b) at (1,-1) {$\CC$};
    \node (c) at (-3,1) {$\pslice{Y}{C}$};
    \node (d) at (-4,-1) {$\pslice{X}{C}$};
    \draw[->] (x) to node[left] {$x^*$} (a);
    \draw[->] (x) to node[right] {$Q_Y$} (b);
    \draw[->] (a) to node[below] {$Q_X$} (b);
    \draw[->] (c) to node[left] {$x^\bullet$} (d);
    \draw[->] (c) to node[above] {$I_Y$} (x);
    \draw[->] (d) to node[below] {$I_X$} (a);
    \end{tikzpicture}\]
so that given $\alpha\in\MM(X)$ we can define $\MM(x)\alpha$ by whiskering along $x^*$: this is well-defined since the commutativity of the diagram implies that whiskering along $x^*$ preserves the almost-pseudonaturality.
As whiskering along $x^*$ also sends modifications to modifications and preserves composites of lax transformations on the nose, 
$\MM(x)$ is in fact a strict monoidal functor. 

Unwinding this definition amounts to defining the $1$-cells of $\MM(x)\alpha$ by setting $(\MM(x)\alpha)_{(A,f)}:=\alpha_{(A,fx)}$ for $(A,f)\in \lslice{Y}{C}$ 
and the $2$-cells by setting $(\MM(x)\alpha)_{j,\sigma}=\alpha_{(j,\sigma x)}$ for $(j,\sigma)\colon (A,f)\to (B,g)$ in $\lslice{Y}{C}$, \
\ie 
  \[
   \begin{pic}
    \node[morphism] (f) at (0,0) {$(\MM(x)\alpha)_{(j,\sigma)}$};
    \draw (f.south west) to +(0,-.5) node[left] {$(\MM(x)\alpha)_{(A,f)}$};
    \draw (f.north west) to +(0,.5) node[left] {$j$};
    \draw (f.south east) to +(0,-.5) node[right] {$j$};
    \draw (f.north east) to +(0,.5) node[right] {$(\MM(x)\alpha)_{(B,g)}$};
  \end{pic} \quad
  :=\quad\begin{pic}
    \node[morphism] (f) at (0,0) {$\alpha_{(j,\sigma x)}$};
    \draw (f.south west) to +(0,-.5) node[left] {$\alpha_{(A,fx)}$};
    \draw (f.north west) to +(0,.5) node[left] {$j$};
    \draw (f.south east) to +(0,-.5) node[right] {$j$};
    \draw (f.north east) to +(0,.5) node[right] {$\alpha_{(B,gx)}$};
  \end{pic}.\]
Given a modification $\Gamma\colon \alpha\to\beta$, the modification $\MM(x)\Gamma\colon \MM(x)\alpha\to\MM(x)\beta$ is 
then obtained by whiskering along $x^*$: explicitly this means that at an object $(A,f)$ of  $\lslice{Y}{C}$ we have $(\MM(x)\Gamma)_{(A,f)}=\Gamma_{(A,fx)}$. 

We next give the action of $\MM$ on a $2$-cell $\theta\colon x\to y$, resulting in a monoidal natural transformation 
$\MM(\theta)\colon \MM(x)\to\MM(y)$. First of all, $\theta$ induces a $2$-natural transformation 
$\theta^*\colon x^*\to y^*$ whose component $\theta^*_{(A,f)}$ at an object $(A,f)$ of $\lslice{Y}{C}$ is given by the $1$-cell $(\id[A],f. \theta)$
    \[\begin{tikzpicture}
    \node (x) at (0,3) {$X$};
    \node (y1) at (-1,1) {$Y$};
    \node (y2) at (1,1) {$Y$};
    \node (a) at (-2,-1) {$A$};
    \node (b) at (2,-1) {$A$};
    \draw[->] (x) to coordinate[midway] (h) node[left] {$x$} (y1);
    \draw[->] (x) to coordinate[midway] (j) node[right] {$y$} (y2);
    \draw[->] (y1) to coordinate[midway] (f) node[left] {$f$} (a);
    \draw[->] (y2) to  coordinate[midway] (g)  node[right] {$f$} (b);
    \draw[->] (a) to node[below] {$\id$} (b);
    \draw[->] (y1) to node[below] {$\id$} (y2);
    \draw[shorttwocell] (h) to node[above] {$\theta$} (j);
    \end{tikzpicture}\]
in $\lslice{X}{C}$. Note that this really requires us to work with the lax slice, as the components of $\theta^*$ live in the pseudoslice only if $\theta$ is invertible. Now, if we take $\alpha \in \MM(X)$, then we have a diagram
\[
\begin{tikzpicture}
    \node (a) at (-2,-1) {$\lslice{Y}{C}$};
    \node (b) at (1,-1) {$\lslice{X}{C}$};
    \node (c) at (3,-1) {$\CC$};
    \draw[->] (a) to[out=35,in=145] coordinate[midway] (h1) node[above]{$x^*$} (b);
    \draw[->] (a) to[out=-35,in=215] coordinate[midway] (h2) node[below] {$y^*$} (b);
    \draw[shorttwocell] (h1) to node[right] {$\theta^*$} (h2);
    \draw[->] (b) to[out=35,in=145] coordinate[midway] (g1) node[above]{$Q_X$} (c);
    \draw[->] (b) to[out=-35,in=215] coordinate[midway] (g2) node[below] {$Q_X$} (c);
    \draw[shorttwocell] (g1) to node[right] {$\alpha$} (g2);
  \end{tikzpicture}
\]
Consequently, we may consider the interchange modification $\Psi[\theta^*,\alpha]$\/ as defined in~\eqref{eq:interchange}; 
noting that $Q_X.\theta^*=1$, we get a modification 
  \[\MM(x)\alpha=\alpha x^*=\id\circ (\alpha x^*)= (Q_X \theta^*)\circ (\alpha x^*)\xrightarrow{\Psi[\theta^*,\alpha]} (\alpha y^*)\circ (Q_X\theta^*)=\MM(y)\alpha. \]
By Lemma~\ref{lem:inter} part 3, we have, for each modification $\Gamma\colon \alpha \to \beta$, a commutative square
\[ 
\begin{tikzpicture}
\node (a) at (-3,1) {$(Q_X \theta^*)\circ (\alpha x^*)$};
\node (b) at (2,1) {$ (\alpha y^*)\circ (Q_X\theta^*)$};
\node (c) at (-3,-1) {$(Q_X \theta^*)\circ (\beta x^*)$};
\node (d) at (2,-1) {$ (\beta y^*)\circ (Q_X\theta^*)$};
\draw[->] (a) to node[above] {$\Psi[\theta^*,\alpha]$} (b);
\draw[->] (a) to node[left] {$(Q_X \theta^*) \circ (\Gamma.x^*) $} (c);
\draw[->] (c) to node[below] {$\Psi[\theta^*,\beta]$} (d);
\draw[->] (b) to node[right] {$(\Gamma.y^*) \circ (Q_X \theta^*)$} (d);
\end{tikzpicture}
\]
from which it follows that $\MM(\theta)$\/ is a natural transformation. Moreover, the first two properties of interchangers stated in Lemma~\ref{lem:inter}
 imply that this natural transformation is a monoidal one: the tensor unit of $\MM(x)$\/ is the 
identity transformation for which this interchanger is the identity, and the interchanger at a composite 
agreeing with the composite of the interchangers implies that $\MM(\theta)_{\alpha\otimes \beta}=\MM(\theta)_\alpha\otimes \MM(\theta)_\beta$.

We now describe explicitly what $\MM(\theta)$ amounts to by describing its component modifications 
$\MM(\theta)_\alpha\colon\MM(x)\alpha\to \MM(y)\alpha$ for each $\alpha\in  \MM(X)$. Given $\alpha\in \MM(X)$, the data of such a modification consists of a 
$2$-cell $\alpha_{(A,fx)}=(\MM(x)\alpha)_{(A,f)}\to (\MM(y)\alpha)_{(A,f)}=\alpha_{(A,fy)}$ for any object $(A,f)$ of $\lslice{Y}{C}$. 
Observe now that $\theta$ induces a $1$-cell $(\id,f\theta)\colon (A,fx)\to (A,fy)$ in $\lslice{X}{C}$, so that $\alpha_{(\id,f\theta)}$ is a 
$2$-cell $\alpha_{(A,fx)}\to \alpha_{(A,fy)}$: it is these $2$-cells that form $\MM(\theta)_\alpha$.

We summarize all of the above in the following.

\begin{theorem}
For a $2$-category $\CC$, the assignment $X \mapsto \MM(X)$\/ underlies a well-defined $2$-functor $\MM\colon \CC \to \cat{MonCat}$, and
$X \mapsto \ZZ(X)$\/ underlies a well-defined $2$-functor $\ZZ\colon \CC \to \cat{2Grp}$. 
\end{theorem}

To conclude this section, we note that the isotropy $2$-group truly is a generalization of the one-dimensional isotropy group: 
if $\CC$\/ is a $1$-category viewed as a locally discrete $2$-category $d\CC$,  
then the resulting isotropy $2$-group coincides with viewing the isotropy group $Z\colon\CC \to \Grp$ as a locally discrete $2$-group. More precisely,
letting $d\colon \cat{Cat} \to \cat{2Cat}$\/ be the $2$-functor that is left adjoint to the $2$-functor $\cat{2Cat} \to \cat{Cat}$\/ that forgets the $2$-cells,
we obtain an extension of the composite
\[  \CC \to \Grp \to \cat{2Grp}\]
to a $2$-functor $\tilde{Z}\colon  d\CC \to \cat{2Grp}$; this extension is precisely $\ZZ\colon d\CC \to \cat{2Grp}$. 

\section{Isotropy in the presence of binary coproducts}\label{sec:2disotropywithcoprods}

To be able to compute $\ZZ$, we now move on generalize the tools of Section~\ref{sec:1d} to the $2$-categorical setting, focusing on coproducts in this section.
Throughout this section, we assume that $\CC$ has binary coproducts. To clarify the level of weakness intended,
 we refer to coproducts in the bicategorical sense, so that $\CC(A+B,-)$ represents $\CC(A,-)\times\CC(B,-)$ up to a pseudonatural equivalence, 
 which we assume to be fixed once and for all. We will write $[-,-]$ for the equivalence 
 \[ \hom_\CC(A,-)\times\hom_\CC(B,-)\to\hom_\CC(A+B,-),\]
  so that $[f,g]$ gives the chosen factorization of $f$ and $g$ up to isomorphism through the coproduct injections. 
  To avoid cluttering the diagrams, we will suppress the chosen invertible $2$-cells from the notation, 
  but we do not mean to assume that diagrams involving coproducts commute strictly. 

With our assumptions, for each $X$, the coproducts in $\CC$\/ induce a pseudofunctor $L_X\colon \CC\to\pslice{X}{C}$\/ sending    
\[\begin{tikzpicture}
    \node (a) at (-1,-1) {$A$};
    \node (b) at (1,-1) {$B$};
    \draw[->] (a) to[out=35,in=145] coordinate[midway] (h1) node[above]{$k$} (b);
    \draw[->] (a) to[out=-35,in=215] coordinate[midway] (h2) node[below] {$j$} (b);
    \draw[shorttwocell] (h2) to node[right] {$\theta$} (h1);
    \end{tikzpicture}\]
     in $\CC$ to 

    \[\begin{tikzpicture}
    \node (x) at (0,1) {$X$};
    \node at (0,.5) {$\cong$};
    \node (a) at (-2,-1) {$X+A$};
    \node (b) at (2,-1) {$X+B$};
    \draw[->] (x) to[out=200,in=80] coordinate[midway] (f) node[left] {$i_X$} (a);
    \draw[->] (x) to[out=-15,in=100] coordinate[midway] (g) node[right] {$i_X$} (b);
    \draw[->] (a) to[out=25,in=155] coordinate[midway] (h1) node[above]{$X+k$} (b);
    \draw[->] (a) to[out=-25,in=205] coordinate[midway] (h2) node[below] {$X+j$} (b);
    \draw[shorttwocell] (h2) to node[right] {$X+\theta$} (h1);
    \end{tikzpicture}\] 
 in  $\pslice{X}{C}$.

Now, $L_X$ is the left biadjoint of $\pslice{X}{C}\xrightarrow{P_X}\CC$.
 
In particular, taking pseudomates (which can be done once one chooses all the required data of the biadjunction, see~~\cite[Section 3.1]{lauda:frobenius}) induces a monoidal equivalence between the monoidal category of all pseudonatural endomorphisms of 
$L_X$\/ and the reverse of that of all pseudonatural endomorphisms of $P_X$. 
Here by the reverse $\DD\rev$\/ of a monoidal category $\DD$ we refer to the monoidal category obtained from $\DD$\/ by 
keeping the underlying category fixed but switching the order of the tensor product, so that $A\otimes\rev B:=B\otimes A$. 
Before showing that $\ZZ$\/ can be computed in terms of automorphisms of $L_X$, we first show that the difference between 
pseudo- and lax slices is immaterial in the presence of binary coproducts. 
Towards this aim, let us first consider the composite $I_XL_XQ_X\colon \lslice{X}{C} \to \lslice{X}{C}$, which sends an object $(A,f)$\/ to $i_X\colon X \to X+A$. 
We now observe that there is a canonical lax natural transformation $\tau\colon I_XL_XQ_X\to\id[\lslice{X}{C}]$\/ whose component at $(A,f)$\/ is the $1$-cell
 \[\begin{tikzpicture}
    \node (x) at (0,1) {$X$};
    \node (a) at (-1,-1) {$X+A$};
    \node (b) at (1,-1) {$A$};
    \draw[->] (x) to coordinate[midway] (f) node[left] {$i_X$} (a);
    \draw[->] (x) to  coordinate[midway] (g)  node[right] {$f$} (b);
    \draw[->] (a) to node[below] {$[f,\id]$} (b);
    \draw[shorttwocell] (f) to node[above] {$\cong$} (g);
    \end{tikzpicture}\] 
and with lax naturality square at $(j,\sigma)\colon (A,f)\to (B,g)$\/ given by 
  
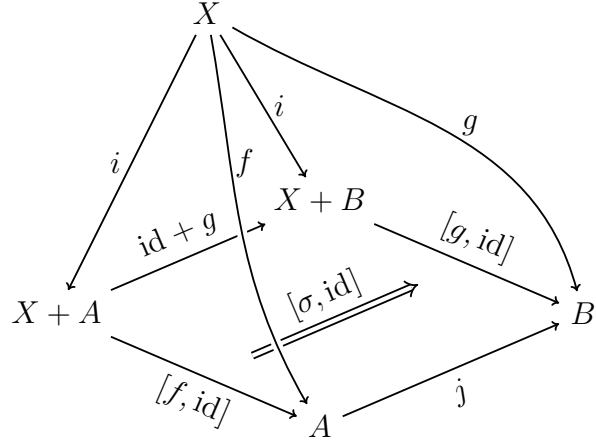
\[
  \begin{tikzpicture}
    \node (a) at (-2,0) {$X+A$};
     \node (b1) at (1.5,1.5) {$X+B$};
    \node (b2) at (1.5,-1.5) {$A$};
    \node (c) at (5,0) {$B$};
    \node (x) at (0,4) {$X$};
    \draw[->] (a) to node[above,sloped]  {$\id+g$} (b1);
    \draw[->] (a) to coordinate[midway](mu) node[below,sloped]  {$[f,\id]$} (b2);
    \draw[->] (b2) to node[below,sloped]  {$j$} (c);
   \draw[->] (b1) to coordinate[midway](nu)  node[above,sloped] {$[g,\id]$} (c);
   \draw[shorttwocell,shorten >=20pt, shorten <=20pt] (mu) to  node[above,sloped] {$[\sigma,\id]$} (nu);
   \draw[->] (x) to node[left]  {$i$} (a);
   \draw[->] (x) to node[right]  {$i$} (b1);
   \draw[->] (x) to[out=-30,in=105] node[right]  {$$} (c);
   \draw[draw=white,double=white,double distance=\pgflinewidth,ultra thick] (x) to[in=120,out=-80] (b2);
    \draw[->] (x) to[in=120,out=-80] node[above] {$$} (b2);
    \node (g) at (3.5,2.5) {$g$};
      \node (f) at (0.5,2) {$f$};
  \end{tikzpicture}
  \]

\begin{theorem}\label{thm:isowithcoproducts} If $\CC$ has binary coproducts, then whiskering along 
$\pslice{X}{C}\hookrightarrow\lslice{X}{C}$ induces an isomorphism of monoidal categories $\MM(X)\to\Mps(X)$ and of $2$-groups $\ZZ(X)\to\Zps(X)$
\end{theorem}

\begin{proof} 
  It suffices to prove the claim for $\MM$ and $\Mps(X)$---if two monoidal categories are isomorphic, 
  so are their associated $2$-groups on weakly invertible objects. We'll use $G_X$\/ to  denote the functor $\MM(X)\to\Mps(X)$ that acts by 
  whiskering along $I_X\colon \pslice{X}{C}\hookrightarrow\lslice{X}{C}$. Note first that $G_X$\/ is a strict monoidal functor.  
  We will first show that $G_X$ is an equivalence, by exhibiting a pseudoinverse $F_X\colon \Mps(X)\to \MM(X)$ and 
  natural isomorphisms $F_X G_X\cong \id$ and $G_X F_X\cong \id$. 
  The pseudoinverse  $F_X\colon \Mps(X)\to \MM(X)$\/ is defined by
  \[\begin{pic}
    \node[morphism] (f) at (0,0) {$F_X(\mu)$};
    \draw (f.south) to +(0,-.5) node[left] {$Q_X$};
    \draw (f.north) to +(0,.5) node[right] {$Q_X$};
    \end{pic}\quad:=\quad \begin{pic}
    \node[morphism] (f) at (.85,-1) {$\mu$};
    \node[morphism] (g) at (.85,0) {$=$};
    \setlength\minimummorphismwidth{14mm}
    \node[morphism] (h) at (0,1) {$\tau$};
    \draw (f.north) to (g.south);
    \draw ([xshift=-2pt]g.north west) to ([xshift=2.5pt]h.south east) node[below left] {$I_X$};
    \draw ([xshift=2pt]g.north east) to +(0,1.5) node[right] {$Q_X$};
    \draw (f.south) to[out=-90,in=0] ++(-0.425,-.5) to[out=180,in=-90] ++(-0.425,.5) to (h.south) node[below left] {$L_X$};
    \draw ([xshift=-2.5pt]h.south west) to +(0,-3) node[left] {$Q_X$};
  \end{pic}\]
  where $\tau$ is the lax natural transformation $I_XL_XQ_X\to\id[\lslice{X}{C}]$\/ described above.

  We now describe $F_X$ explicitly. Given $\mu\in \Mps(X)$\/ and an object $(A,f)$\/ of $\lslice{X}{C}$,  $F_X(\mu)_{(A,f)}$\/ is defined as the composite 
  \begin{equation}\label{eq:invmod} A\xrightarrow{i_A} X+A\xrightarrow{\mu_{i_X}}X+A\xrightarrow{[f,\id]}A.\end{equation}
  Given a $1$-cell $(j,\sigma)\colon (A,f)\to (B,g)$, the $2$-cell $F_X(\mu)_{(j,\sigma)}$\/ arises via the diagram 
    
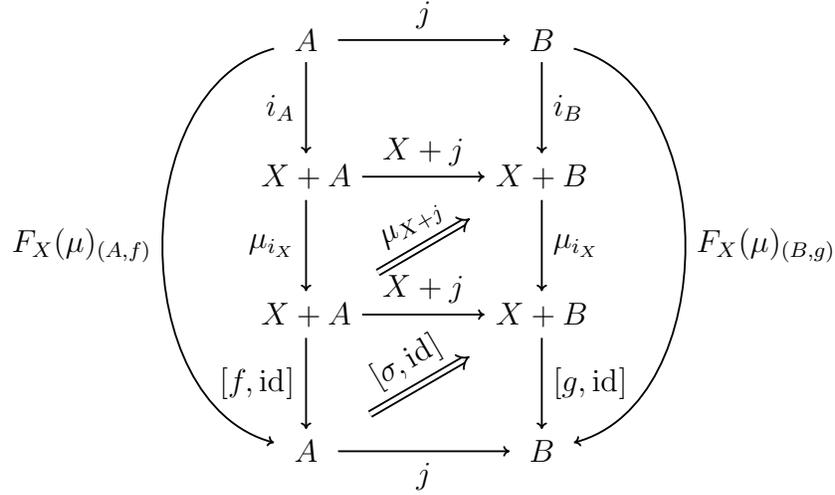
\[\begin{tikzpicture}
    \matrix (m) [matrix of math nodes,row sep=3em,column sep=4em,minimum width=2em]
    {
     A & B\\
     X+A & X+B \\
     X+A & X+B \\
     A & B \\};
    \path[->]
    (m-1-1) edge node [left] {$i_A$}  (m-2-1)
           edge node [above] {$j$} (m-1-2)
           edge[in=165,out=195] node [left] {$F_X(\mu)_{(A,f)}$} (m-4-1)
    (m-2-1) edge node[left] {$\mu_{i_X}$}  coordinate[midway] (e) (m-3-1)
            edge node[above]  {$X + j$}  coordinate[midway] (d) (m-2-2)
    (m-2-2) edge node[right] {$\mu_{i_X}$} (m-3-2)
    (m-1-2) edge node [right] {$i_B$} (m-2-2)
            edge[out=-15,in=15] node [right] {$F_X(\mu)_{(B,g)}$} (m-4-2)
    (m-3-1) edge node[left] {$[f,\id]$} (m-4-1)
            edge node[above] {$X+j$} coordinate[midway] (c)  (m-3-2)
    (m-3-2) edge node[right] {$[g,\id]$} (m-4-2)
    (m-4-1) edge node [below] {$j$} (m-4-2);
    \draw[shorttwocell,shorten >=14pt, shorten <=14pt] (m-3-1) to node[above,sloped] {$\mu_{X+j}$} (m-2-2);
     \draw[shorttwocell,shorten >=14pt, shorten <=14pt] (m-4-1) to node[above,sloped] {$[\sigma,\id]$} (m-3-2);
  \end{tikzpicture}\]
  On morphisms, $F_X$\/ acts by sending $\Gamma\colon \mu\to \nu$ to the modification 
  $F_X(\Gamma)\colon F_X(\mu)\to F_X(\nu)$\/ where $F_X(\Gamma)_{(A,f)}$ is defined by
  \[
  \begin{tikzpicture}
    \node (a) at (-2,0) {$A$};
    \node (c) at (5,0) {$A$};
    \node (d) at (0,0) {$X+A$};
    \node (e) at (3,0) {$X+A$};
    \draw[->] (a) to node[above]  {$i_A$} (d);
    \draw[->] (d) to[out=35,in=145] coordinate[midway](mu) node[above]  {$\mu_A$} (e);
     \draw[->] (d) to[out=-35,in=215] coordinate[midway](nu) node[below]  {$\nu_A$} (e);
    \draw[shorttwocell] (mu) to  node[right] {$\Gamma_A$} (nu);
    \draw[->] (e) to node[above]  {$[f,\id]$} (c);
  \end{tikzpicture}
  \]

Note that it follows from the description~\eqref{eq:invmod} by taking $\mu=1$, we get:
    \begin{equation}\label{eq:simplifytau}
    \begin{pic}
    \setlength\minimummorphismwidth{10mm}
    \node[morphism] (h) at (0,-1) {$=$};
    \setlength\minimummorphismwidth{14mm}
    \node[morphism] (e) at (-1,0) {$\tau$};
    \draw (h.north east) to ++(0,1) node[right] {$Q_X$};
    \draw (h.north west) to node[left] {$$} ([xshift=1pt]e.south east);
    \draw (h.south) to[out=-90,in=0] ++(-0.5,-.5) to[out=180,in=-90] ++(-0.5,.5) to (e.south);
    \draw ([xshift=-1pt]e.south west)  to +(0,-2) node[left] {$Q_X$};
  \end{pic}\cong\quad 
  \begin{pic} 
      \draw (0,0) to +(0,2) node[right] {$Q_X$};
      \end{pic}
  \end{equation}

 We now show that $F_X$ is the pseudoinverse of $G_X$. Note first that the string diagram representation of $G_X$\/ is   
  \[ \begin{pic} 
    \node[morphism] (f) at (0,0) {$G_X(\alpha)$};
    \draw (f.south) to +(0,-.5) node[left] {$P_X$};
    \draw (f.north) to +(0,.5) node[right] {$P_X$};
    \end{pic}\quad:=\begin{pic}
    \node[morphism] (f) at (0.4,0) {$\alpha$};
    \setlength\minimummorphismwidth{10mm}
    \node[morphism] (g) at (0,1) {$=$};
    \node[morphism] (h) at (0,-1) {$=$};
    \draw (g.south east) to node[left] {$Q_X$} (f.north);
    \draw (h.north east) to node[left] {$Q_X$} (f.south);
    \draw (g.south west) to node[left] {$I_X$} (h.north west);
    \draw (g.north) to node[right] {$P_X$}  ++(0,0.75);
    \draw (h.south) to node[right] {$P_X$}  ++(0,-.75);
  \end{pic}\] 

  Then the isomorphism $F_XG_X\to\id$ at $\alpha\in\MM(X)$ is the composite
  \[ \begin{pic} 
    \node[morphism] (f) at (0,0) {$F_X G_X(\alpha)$};
    \draw (f.south) to +(0,-.5) node[left] {$Q_X$};
    \draw (f.north) to +(0,.5) node[right] {$Q_X$};
    \end{pic}\quad:=
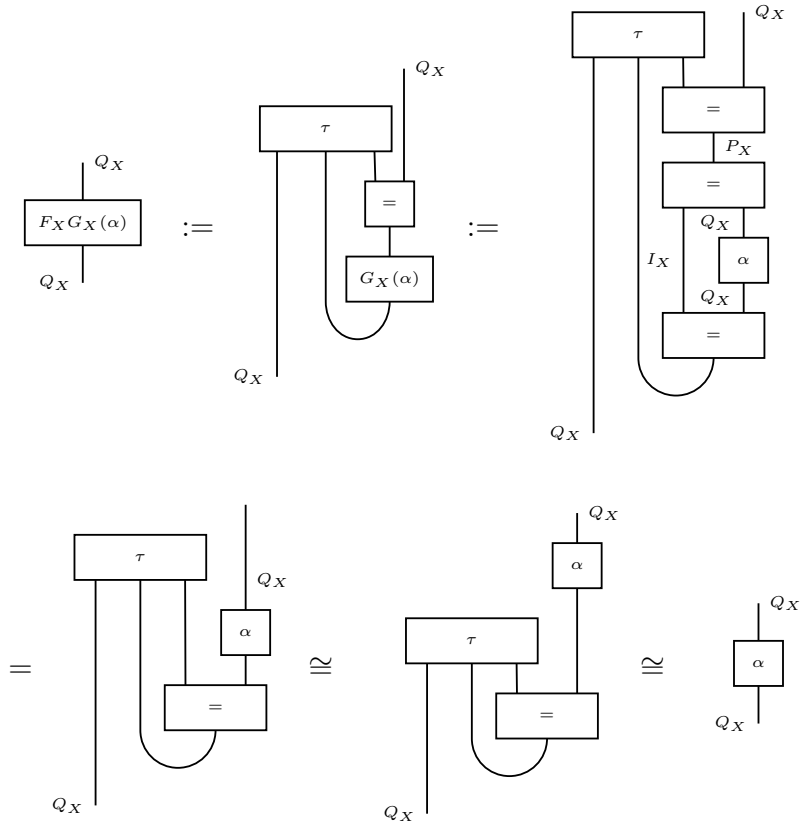
    \begin{pic}
    \node[morphism] (f) at (.85,-1) {$G_X(\alpha)$};
    \node[morphism] (g) at (.85,0) {$=$};
    \setlength\minimummorphismwidth{14mm}
    \node[morphism] (h) at (0,1) {$\tau$};
    \draw (f.north) to (g.south);
    \draw ([xshift=-2pt]g.north west) to ([xshift=2.5pt]h.south east);
    \draw ([xshift=2pt]g.north east) to +(0,1.5) node[right] {$Q_X$};
    \draw (f.south) to[out=-90,in=0] ++(-0.425,-.5) to[out=180,in=-90] ++(-0.425,.5) to (h.south);
    \draw ([xshift=-2.5pt]h.south west) to +(0,-3) node[left] {$Q_X$};
    \end{pic} :=\quad 
    \begin{pic}
    \node[morphism] (f) at (0.4,0) {$\alpha$};
    \setlength\minimummorphismwidth{10mm}
    \node[morphism] (g) at (0,1) {$=$};
    \node[morphism] (h) at (0,-1) {$=$};
     \node[morphism] (i) at (0,2) {$=$};
    \draw (g.south east) to node[left] {$Q_X$} (f.north);
    \draw (h.north east) to node[left] {$Q_X$} (f.south);
    \draw (g.south west) to node[left] {$I_X$} (h.north west);
    \draw (g.north) to node[right] {$P_X$}  (i.south);
    \setlength\minimummorphismwidth{14mm}
    \node[morphism] (e) at (-1,3) {$\tau$};
    \draw (h.south) to[out=-90,in=0] ++(-0.5,-.5) to[out=180,in=-90] ++(-0.5,.5) to (e.south);
    \draw (i.north west) to ([xshift=1pt]e.south east);
    \draw (i.north east) to +(0,1) node[right] {$Q_X$};
    \draw ([xshift=-1pt]e.south west)  to +(0,-5) node[left] {$Q_X$};
  \end{pic}\] 

  \[=\begin{pic}
    \node[morphism] (f) at (0.4,0) {$\alpha$};
    \setlength\minimummorphismwidth{10mm}
    \node[morphism] (h) at (0,-1) {$=$};
    \setlength\minimummorphismwidth{14mm}
    \node[morphism] (e) at (-1,1) {$\tau$};
    \draw (h.north east) to node[left] {$$} (f.south);
    \draw (h.north west) to node[left] {$$} ([xshift=1pt]e.south east);
    \draw (h.south) to[out=-90,in=0] ++(-0.5,-.5) to[out=180,in=-90] ++(-0.5,.5) to (e.south);
    \draw (f.south) to node[right] {$Q_X$} ++(0,2);
    \draw ([xshift=-1pt]e.south west)  to +(0,-3) node[left] {$Q_X$};
  \end{pic}
  \cong \quad 
  \begin{pic}
    \node[morphism] (f) at (0.4,1) {$\alpha$};
    \setlength\minimummorphismwidth{10mm}
    \node[morphism] (h) at (0,-1) {$=$};
    \setlength\minimummorphismwidth{14mm}
    \node[morphism] (e) at (-1,0) {$\tau$};
    \draw (h.north east) to node[left] {$$} (f.south);
    \draw (h.north west) to node[left] {$$} ([xshift=1pt]e.south east);
    \draw (h.south) to[out=-90,in=0] ++(-0.5,-.5) to[out=180,in=-90] ++(-0.5,.5) to (e.south);
    \draw (f.south) to  ++(0,1) node[right] {$Q_X$};
    \draw ([xshift=-1pt]e.south west)  to +(0,-2) node[left] {$Q_X$};
  \end{pic}\cong\quad 
  \begin{pic} 
      \node[morphism] (f) at (0,0) {$\alpha$};
      \draw (f.south) to +(0,-.5) node[left] {$Q_X$};
      \draw (f.north) to +(0,.5) node[right] {$Q_X$};
      \end{pic}
  \]
  where the penultimate isomorphism is an instance of invertible interchange (as $\alpha$ is almost pseudonatural), and the last one is~\eqref{eq:simplifytau}.
  
  Explicitly, the modification $F_XG_X(\alpha)\to\alpha$ at the index $(A,f)$ is given by
  \begin{equation}\label{diag:modificationattheroundtrip}\vcenter{\hbox{
  \begin{tikzpicture}
  \node (a) at (-2,0) {$A$};
  \node (b) at (0.5,-0.5) {$X+A$};
  \node (c0) at (3,1.25) {$X+A$};
  \node (c1) at (3,-1.25) {$X+A$};
   \node (d0) at (7,1.25) {$A$};
  \node (d1) at (7,-1.25) {$A$};
  \draw[->] (a) to node[below]  {$i_A$} (b);
  \draw[->] (a) to node[above]  {$i_A$} (c0);
  \draw[->] (a) to[in=165,out=45] node[above]  {$\id$} (d0);
  \draw[->] (a) to[in=-165,out=-45] node[below]  {$F_XG_X(\alpha)_{(A,f)}$} (d1);
  \draw[->] (b) to node[below]  {$G(\alpha)_{i_X}$} (c1);
  \draw[->] (c0) to node[below]  {$[f,\id]$} (d0);
  \draw[->] (c0) to node[left]  {$\alpha_{i_X}$} (c1); 
  \draw[->] (c1) to node[above]  {$[f,\id]$} (d1);
  \draw[->] (d0) to node[right]  {$\alpha_{(A,f)}$} (d1);
  \draw[shorttwocell,shorten >=10pt, shorten <=10pt] (c1) to  node[below,sloped] {\qquad$\alpha_{([f,\id],\cong)}$} (d0);
  \end{tikzpicture}}}
  \end{equation}

  The isomorphism $G_XF_X\to\id$ at $\mu\in\Mps(X)$ is the composite
  \[\begin{pic} 
    \node[morphism] (f) at (0,0) {$G_XF_X(\mu)$};
    \draw (f.south) to +(0,-.5) node[left] {$P_X$};
    \draw (f.north) to +(0,.5) node[right] {$P_X$};
    \end{pic}\quad:=\begin{pic}
  \node[morphism] (f) at (0.4,0) {$F_X(\mu)$};
  \setlength\minimummorphismwidth{10mm}
  \node[morphism] (g) at (0,1) {$=$};
  \node[morphism] (h) at (0,-1) {$=$};
  \draw (g.south east) to node[left] {$Q_X$} (f.north);
  \draw (h.north east) to node[left] {$Q_X$} (f.south);
  \draw (g.south west) to node[left] {$I_X$} (h.north west);
  \draw (g.north) to node[right] {$P_X$}  ++(0,0.75);
  \draw (h.south) to node[right] {$P_X$}  ++(0,-.75);
  \end{pic}:=\quad
  \begin{pic}
    \node[morphism] (f) at (.85,-1) {$\mu$};
    \node[morphism] (g) at (.85,0) {$=$};
     \node[morphism] (i) at (-.85,-2.5) {$=$};
    \setlength\minimummorphismwidth{14mm}
    \node[morphism] (h) at (0,1) {$\tau$};
     \setlength\minimummorphismwidth{22mm}
     \node[morphism] (j) at (0,2) {$=$};
    \draw (f.north) to (g.south);
    \draw ([xshift=-2pt]g.north west) to ([xshift=2.5pt]h.south east);
    \draw ([xshift=2pt]g.north east) to ([xshift=4.5pt]j.south east);
    \draw (f.south) to[out=-90,in=0] ++(-0.425,-.5) to[out=180,in=-90] ++(-0.425,.5) to (h.south);
    \draw ([xshift=-2.5pt]h.south west) to ([xshift=2pt]i.north east);
    \draw ([xshift=-2pt]i.north west) to ([xshift=-4.5pt]j.south west);
    \draw (i.south) to +(0,-.75);
    \draw (j.north) to +(0,.75);
    \end{pic}
  \] 

  \[=\quad 
  \begin{pic}
    \node[morphism] (f) at (.85,-1) {$\mu$};
    \node[morphism] (g) at (.85,0) {$=$};
     \node[morphism] (i) at (-.85,0) {$=$};
    \setlength\minimummorphismwidth{14mm}
    \node[morphism] (h) at (0,1) {$\tau$};
     \setlength\minimummorphismwidth{22mm}
     \node[morphism] (j) at (0,2) {$=$};
    \draw (f.north) to (g.south);
    \draw ([xshift=-2pt]g.north west) to ([xshift=2.5pt]h.south east);
    \draw ([xshift=2pt]g.north east) to ([xshift=4.5pt]j.south east);
    \draw (f.south) to[out=-90,in=0] ++(-0.425,-.5) to[out=180,in=-90] ++(-0.425,.5) to (h.south);
    \draw ([xshift=-2.5pt]h.south west) to ([xshift=2pt]i.north east);
    \draw ([xshift=-2pt]i.north west) to ([xshift=-4.5pt]j.south west);
    \draw (i.south) to +(0,-.75);
    \draw (j.north) to +(0,.75);
    \end{pic}
    = \begin{pic} 
    \node[morphism] (f) at (0,0) {$\mu$};
    \draw (f.south)  to[out=-90,in=0] ++(-0.5,-.5) to[out=180,in=-90] ++(-0.5,.5) to ++(0,1) to[out=90,in=0] ++(-0.5,.5) to[out=180,in=90] ++(-0.5,-.5) to ++(0,-2) node[left] {$P_X$};
    \draw (f.north) to +(0,1.5) node[right] {$P_X$};
    \end{pic}\cong 
    \begin{pic} 
    \node[morphism] (f) at (0,1) {$\mu$};
    \draw (f.south)  to ++(0,-1.5) to[out=-90,in=0] ++(-0.5,-.5) to[out=180,in=-90] ++(-0.5,.5) to ++(0,1) to[out=90,in=0] ++(-0.5,.5) to[out=180,in=90] ++(-0.5,-.5) to ++(0,-2) node[left] {$P_X$};
    \draw (f.north) to +(0,.5) node[right] {$P_X$};
    \end{pic}
    \cong
    \begin{pic} 
    \node[morphism] (f) at (0,0) {$\mu$};
    \draw (f.south) to +(0,-.5) node[left] {$P_X$};
    \draw (f.north) to +(0,.5) node[right] {$P_X$};
    \end{pic}\]
  where we used the equality 
  
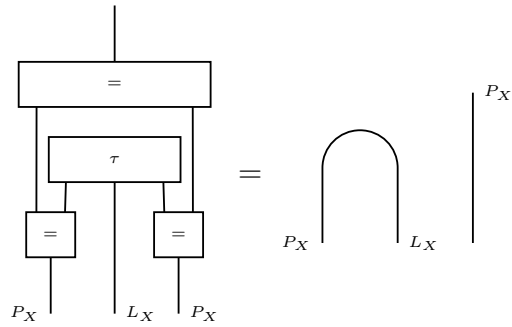
\[
  \begin{pic}
    \node[morphism] (g) at (.85,0) {$=$};
     \node[morphism] (i) at (-.85,0) {$=$};
    \setlength\minimummorphismwidth{14mm}
    \node[morphism] (h) at (0,1) {$\tau$};
     \setlength\minimummorphismwidth{22mm}
     \node[morphism] (j) at (0,2) {$=$};
    \draw ([xshift=-2pt]g.north west) to ([xshift=2.5pt]h.south east);
    \draw ([xshift=2pt]g.north east) to ([xshift=4.5pt]j.south east);
    \draw ([xshift=-2.5pt]h.south west) to ([xshift=2pt]i.north east);
    \draw ([xshift=-2pt]i.north west) to ([xshift=-4.5pt]j.south west);
    \draw (h.south) to ++(0,-1.75)  node[right] {$L_X$};
    \draw (g.south) to +(0,-.75) node[right] {$P_X$};
    \draw (i.south) to +(0,-.75)  node[left] {$P_X$};
    \draw (j.north) to +(0,.75);
    \end{pic}=    
  \begin{pic}
      \draw (-2,0) node[left] {$P_X$} to ++(0,1) to[in=180,out=90] ++(.5,.5) to[in=90,out=0] ++(.5,-.5) to ++(0,-1) node[right] {$L_X$};
      \draw (0,0) to +(0,2) node[right] {$P_X$};
  \end{pic}
  \]
  which is a consequence of~\eqref{eq:simplifytau}.
  
  Explicitly, the isomorphism $G_XF_X(\mu)\to\mu$ at $(A,f)$ is rather like $F_XG_X(\alpha)\to\alpha$ and is given by 
  
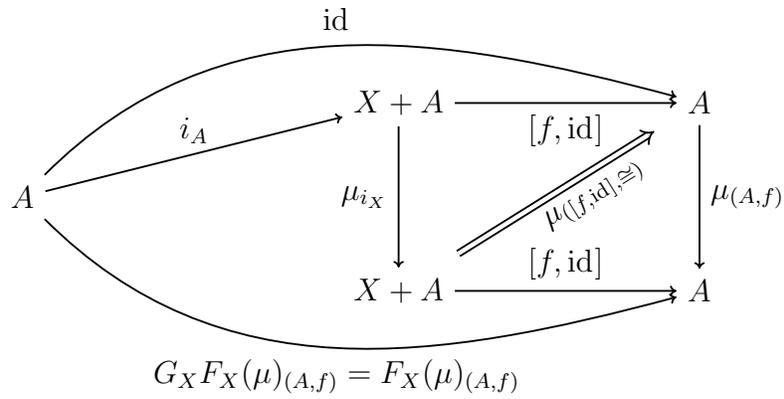
\[
  \begin{tikzpicture}
  \node (a) at (-2,0) {$A$};
  \node (c0) at (3,1.25) {$X+A$};
  \node (c1) at (3,-1.25) {$X+A$};
   \node (d0) at (7,1.25) {$A$};
  \node (d1) at (7,-1.25) {$A$};
  \draw[->] (a) to node[above]  {$i_A$} (c0);
  \draw[->] (a) to[in=165,out=45] node[above]  {$\id$} (d0);
  \draw[->] (a) to[in=-165,out=-45] node[below]  {$G_XF_X(\mu)_{(A,f)}=F_X(\mu)_{(A,f)}$} (d1);
  \draw[->] (c0) to node[below]  {$[f,\id]$} (d0);
  \draw[->] (c0) to node[left]  {$\mu_{i_X}$} (c1);
  \draw[->] (c1) to node[above]  {$[f,\id]$} (d1);
  \draw[->] (d0) to node[right]  {$\mu_{(A,f)}$} (d1);
  \draw[shorttwocell,shorten >=10pt, shorten <=10pt] (c1) to  node[below,sloped] {\qquad$\mu_{([f,\id],\cong)}$} (d0);
  \end{tikzpicture}
  \]

  We now use this to show that $G_X$ is in fact a strict monoidal isomorphism. 
  As we already know that $G_X$ is a strict monoidal functor and an equivalence, it suffices to show that $G_X$ is bijective on objects. 
  To see that $G_X$ is surjective on objects is straightforward. Abstractly, this follows from the fact that any essentially surjective 
  isofibration is in fact strictly surjective on objects. Explicitly, given $\mu\in\Mps(X)$, we can choose $\alpha$ and an invertible modification 
  $\Gamma\colon G_X(\alpha)\to \mu$. Now, there is a unique $\beta\in\MM(X)$ such that the $2$-cells of 
  $\Gamma$ induce an invertible modification $\alpha\to\beta$, and by construction $G_X(\beta)=\mu$.

  To see that $G_X$ is injective, note that the modification $F_XG_X(\alpha)\to\alpha$ being invertible implies that the 
  $2$-cells $\alpha_{(j,\sigma)}$ of $\alpha$  can be recovered by conjugating the $2$-cells $F_XG_X(\alpha)_{(j,\sigma)}$ with 
  the modification. Now, $F_XG_X(\alpha)_{(j,\sigma)}$ depends only on $G_X(\alpha)$, and by inspecting 
  diagram~\ref{diag:modificationattheroundtrip}, so do the component $2$-cells of the modification $F_XG_X(\alpha)\to\alpha$. 
  Hence $G(\alpha)=G(\beta)$ implies $\alpha=\beta$.
\end{proof}

As a result, when binary coproducts exist, obstructions like those of Example~\ref{ex:pseudoisotropyisnot2functorial} cannot arise. 
Indeed, the assignment sending $X$ to pseudonatural endomorphisms of $P_X\colon \pslice{X}{C}\to\CC$\/ becomes a $2$-functor, 
with the action on $1$-cells defined similarly as for $\MM$ (indeed, the action of $\Mps$ on $1$-cells is well-defined even without coproducts). 
The action on $2$-cells can be described as follows. Given a $2$-cell $\begin{pic}
    \node (a) at (-1,-1) {$X$};
    \node (b) at (1,-1) {$Y$};
    \draw[->] (a) to[out=35,in=145] coordinate[midway] (h1) node[above]{$x$} (b);
    \draw[->] (a) to[out=-35,in=215] coordinate[midway] (h2) node[below] {$y$} (b);
    \draw[shorttwocell] (h1) to node[right] {$\sigma$} (h2);
    \end{pic}$ in $\CC$
and $\alpha\in\Mps$,  we can define a modification $\Mps(\sigma)\alpha\colon \Mps(x)\alpha\to\Mps(y)\alpha$\/ by first 
extending $\alpha$ to $G(\alpha)$ and then restricting $\MM(\sigma)$ to obtain a modification  $\Mps(x)\alpha\to\Mps(y)\alpha$. 
Explicitly, the component of the resulting modification at $(A,f)$ is given by 
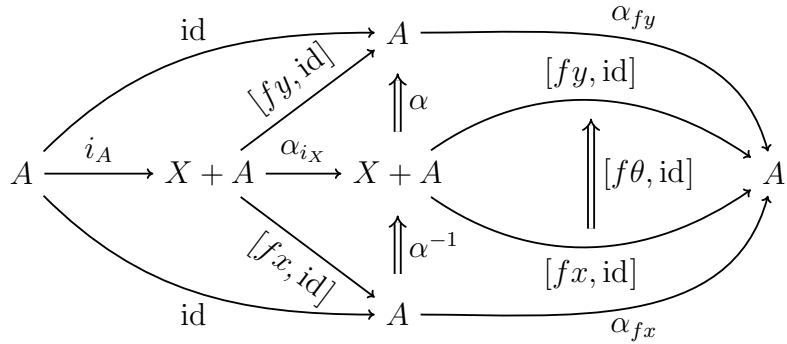
\[
\begin{tikzpicture}[scale=1.25]
  \node (a) at (-2,0) {$A$};
   \node (b1) at (2,1.5) {$A$};
  \node (b2) at (2,-1.5) {$A$};
  \node (c) at (6,0) {$A$};
  \node (d) at (0,0) {$X+A$};
  \node (e) at (2,0) {$X+A$};
  \draw[->] (a) to[in=180,out=45] node[above]  {$\id$} (b1); 
  \draw[->] (a) to node[above]  {$i_A$} (d);
  \draw[->] (a) to[in=180,out=-45] node[below]  {$\id$} (b2); 
  \draw[->] (d) to node[above]  {$\alpha_{i_X}$} (e);
  \draw[->] (d) to node[sloped,above]  {$[fy,\id]$} (b1);
  \draw[->] (d) to node[sloped,below]  {$[fx,\id]$} (b2);
  \draw[->] (e) to[out=35,in=145] coordinate[midway](mu) node[above]  {$[fy,\id]$} (c);
  \draw[->] (e) to[out=-35,in=215] coordinate[midway](nu) node[below]  {$[fx,\id]$} (c);
  \draw[shorttwocell] (nu) to  node[right] {$[f\theta,\id]$} (mu);
  \draw[->] (b1) to[in=105,out=0]  node[above] {$\alpha_{fy}$} (c);
  \draw[->] (b2) to[in=-105,out=0] node[below]  {$\alpha_{fx}$} (c);
  \draw[shorttwocell] (e) to  node[right] {$\alpha$} (b1);
  \draw[shorttwocell] (b2) to  node[right] {$\alpha^{-1}$} (e);
\end{tikzpicture}
\]

We can now also reformulate $\MM$\/ and $\ZZ$\/ in terms of pseudonatural endomorphisms/autoequivalences of $L_X\colon \CC\to\pslice{X}{C}$.
Let us define the following pseudofunctor $\NN\colon\CC \to \cat{MonCat}$:
 \begin{itemize}
  \item On objects, $\NN$ acts by sending $X$ to the reverse of the  monoidal category of all pseudonatural endomorphisms of $L_X\colon \CC\to\pslice{X}{C}$.
  \item For a morphism $x\colon X\to Y$\/ the resulting monoidal functor $\NN(x)\colon\NN(X)\to\NN(Y)$\/ is defined as follows. 
  Given $\mu\in\NN(X)$\/ we define $(\NN(x)\mu)_A$\/ as the $1$-cell $Y+A\to Y+A$\/  defined (up to isomorphism) by the diagram
    \[
    \begin{tikzpicture}
      \node (a1) at (-1,1) {$A$};
      \node (a2) at (-1,-1) {$Y$};
       \node (b) at (1,0) {$Y+A$};
      \node (c) at (4,0) {$Y+A$};
      \node (d) at (1,1) {$X+A$};
      \node (e) at (4,1) {$X+A$};
      \node at (1,-.5) {$\cong$};
        \node at (1,.5) {$\cong$};
      \draw[->] (a1) to node[below]  {$i_A$} (b);
      \draw[->] (a1) to node[above]  {$i_A$} (d);
      \draw[->] (d) to node[above]  {$\mu_A$} (e);
      \draw[->] (e) to node[right]  {$x+A$} (c);
      \draw[->] (a2) to node[above]  {$i_Y$} (b);
      \draw[->] (a2) to[in=-135,out=0] node[below]  {$i_Y$} (c);
     \draw[->] (b) to node[above] {$(\LL(x)\mu)_A$} (c);
    \end{tikzpicture}
    \]
     Given a modification $\Gamma\colon\mu\to\nu$\/ in $\NN(X)$, the resulting modification $\NN(x)\Gamma$\/ has the $A$-th component defined 
    \[
    \begin{tikzpicture}
      \node (a1) at (-2,0) {$A$};
      \node (a2) at (1.5,-1.5) {$Y+A$};
       \node (b) at (1.5,1.5) {$Y+A$};
      \node (c) at (6,0) {$Y+A$};
      \node (d) at (0,0) {$X+A$};
      \node (e) at (3,0) {$X+A$};
      \draw[->] (a1) to node[above]  {$i_A$} (b);
      \draw[->] (a1) to node[below]  {$i_A$} (d);
      \draw[->] (a1) to[out=-75,in=180] node[below]  {$i_A$} (a2);
      \draw[->] (a2) to[out=0,in=270] node[below]  {$(\LL(f)\nu)_A$} (c);
      \draw[->] (d) to[out=35,in=145] coordinate[midway](mu) node[above]  {$\mu_A$} (e);
       \draw[->] (d) to[out=-35,in=215] coordinate[midway](nu) node[below]  {$\nu_A$} (e);
      \draw[shorttwocell] (mu) to  node[right] {$\Gamma_A$} (nu);
      \draw[->] (e) to node[below]  {$f+A$} (c);
     \draw[->] (b) to node[above,sloped] {$(\LL(f)\mu)_A$} (c);
    \end{tikzpicture}
    \]
    on $A$\/
    and on $i_Y$\/ by the identity $2$-cell. 
  \item Given a $2$-cell $\theta$, $\begin{pic}
    \node (a) at (-1,-1) {$X$};
    \node (b) at (1,-1) {$Y$};
    \draw[->] (a) to[out=35,in=145] coordinate[midway] (h1) node[above]{$x$} (b);
    \draw[->] (a) to[out=-35,in=215] coordinate[midway] (h2) node[below] {$y$} (b);
    \draw[shorttwocell] (h1) to node[right] {$\theta$} (h2);
    \end{pic}$ in $\CC$
  the resulting monoidal natural transformation $\NN(\theta)\colon \NN(x)\to \NN(y)$\/ consists of a modification  
  $\NN(\theta)_\mu\colon \NN(x)\mu\to\NN(y)\mu$\/ for each $\mu\in\NN(X)$. On an object $A$, the required $2$-cell is defined on $A$\/ by 
    \[
    \begin{tikzpicture}
      \node (a1) at (-2,0) {$A$};
      \node (a2) at (2,-1.5) {$Y+A$};
       \node (b) at (2,1.5) {$Y+A$};
      \node (c) at (6,0) {$Y+A$};
      \node (d) at (0,0) {$X+A$};
      \node (e) at (2,0) {$X+A$};
      \draw[->] (a1) to[in=180,out=45] node[above]  {$i_A$} (b);
      \draw[->] (a1) to node[above]  {$i_A$} (d);
      \draw[->] (a1) to[in=180,out=-45] node[below]  {$i_A$} (a2);
      \draw[->] (a2) to[in=-90,out=0] node[below]  {$(\LL(y)\mu)_A$} (c);
      \draw[->] (d) to node[above]  {$\mu_A$} (e);
      \draw[->] (e) to[out=35,in=145] coordinate[midway](mu) node[above]  {$x+A$} (c);
      \draw[->] (e) to[out=-35,in=215] coordinate[midway](nu) node[below]  {$y+A$} (c);
      \draw[shorttwocell] (mu) to  node[right] {$\theta+A$} (nu);
     \draw[->] (b) to[in=90,out=0]  node[above] {$(\LL(x)\mu)_A$} (c);
    \end{tikzpicture}
    \]
  \end{itemize}
With this definition, we then obtain:

\begin{corollary}\label{cor:LispseudoequivtoZ} 
  Let $\CC$ be a $2$-category with binary coproducts. Then $\MM$\/ is pseudonaturally equivalent to the pseudofunctor $\NN$\/ defined above.
   Moreover, there is a pseudonatural equivalence between $\ZZ$\/ and the functor $\LL$\/ that sends $X$\/ to the reverse of the 
   $2$-group of all pseudonatural autoequivalences of $L_X\colon \CC\to\pslice{X}{C}$\/ and invertible modifications between them, 
   with the action of $\LL$\/ on 1\/ and $2$-cells defined as for $\NN$.
\end{corollary}

\begin{proof}
  It suffices to prove the claim for $\MM$\/ and $\NN$, as the claim for $\ZZ$\/ and $\LL$\/ follows from it. 
  Taking pseudomates (again, see \eg\cite[Section 3.1]{lauda:frobenius}) induces a monoidal equivalence between the monoidal category of all pseudonatural endomorphisms of 
  $P_X$\/ and the reverse of the monoidal category of all pseudonatural endomorphisms of $L_X$. 
  Composing this with the isomorphisms from Theorem~\ref{thm:isowithcoproducts} then gives an equivalence 
  $\NN(X)\cong \ZZ(X)$\/ for each $X$, and using these pointwise equivalences the assignment 
  $X\mapsto \NN(X)$\/ can be promoted (via Lemma~\ref{lem:promote}) into a pseudofunctor $\NN'$\/ pseudonaturally equivalent to $\ZZ$. 
  We now show that $\NN'$\/ is again equivalent to our desired $\NN$\/ from the statement of the corollary. 
  In fact, as $\NN'$\/ already has the desired action on objects, it suffices to show that each monoidal functor $\NN'(x)$\/  is 
  (monoidally) isomorphic to $\NN(x)$, for these monoidal isomorphisms can then be used to promote $\NN$\/ into a 
  pseudofunctor for which these monoidal isomorphisms then give an invertible icon $\NN'\to\NN$. 

  Now,  $\NN'(x)\mu$\/ is defined by
  \[\begin{pic} 
    \node[morphism] (f) at (0,0) {$\mu$};
    \node[morphism] (e) at (-1.15,-1.25) {$=$};
    \setlength\minimummorphismwidth{26mm}
    \node[morphism] (g) at (-0.15,1.5) {$=$};
    \draw (f.south)  to[out=-90,in=-180] ++(0.5,-.5) to[out=0,in=-90] ++(0.5,.5) to  ++(0,1.1) node[right] {$P_X$} to ([xshift=3pt]g.south east); 
    \draw (f.north) to[out=90,in=0] ++(-0.5,.5) to[out=180,in=90] ++(-0.5,-.5) to ([xshift=1pt]e.north east);
    \draw ([xshift=-1pt]e.north west) to +(0,0.2) node[left] {$x^*$} to ([xshift=-3pt]g.south west);
    \draw (g.north) to[out=90,in=180] ++(1,1) to[out=0,in=90] ++(1,-1) to ++(0,-4.5);
    \draw (e.south) to[out=-90,in=0] ++(-0.5,-.5) to[out=180,in=-90] ++(-0.5,.5) to ++(0,4.5);
    \end{pic}\]
  However, to know an element  $\mu\in \NN(Y)$\/ within isomorphism, it is enough to know the restrictions of the 1\/ and 
  $2$-cell components along inclusions $A\to X+A$\/ for each $A$. Hence $\NN'(X)\mu$\/ is determined by
  \begin{equation}\tag{$*$}\label{pic:mateofLxmu}\begin{pic} 
    \node[morphism] (f) at (0,0) {$\mu$};
    \node[morphism] (e) at (-1.15,-1.25) {$=$};
    \setlength\minimummorphismwidth{26mm}
    \node[morphism] (g) at (-0.15,1.5) {$=$};
    \draw (f.south)  to[out=-90,in=-180] ++(0.5,-.5) to[out=0,in=-90] ++(0.5,.5) to  ++(0,1.1) node[right] {$P_X$} to ([xshift=3pt]g.south east); 
    \draw (f.north) to[out=90,in=0] ++(-0.5,.5) to[out=180,in=90] ++(-0.5,-.5) to ([xshift=1pt]e.north east);
    \draw ([xshift=-1pt]e.north west) to +(0,0.2) node[left] {$x^*$} to ([xshift=-3pt]g.south west);
    \draw (g.north) to[out=90,in=180] ++(1,1) to[out=0,in=90] ++(1,-1) to ++(0,-4.5) to[out=-90,in=180] ++(0.5,-.5) to[out=0,in=-90] ++(0.5,.5) to ++(0,6);
    \draw (e.south) to[out=-90,in=0] ++(-0.5,-.5) to[out=180,in=-90] ++(-0.5,.5) to ++(0,4.5);
    \end{pic}\end{equation}
  We now compute as follows
  \[\begin{pic} 
    \node[morphism] (f) at (0,0) {$\mu$};
    \node[morphism] (e) at (-1.15,-1.25) {$=$};
    \setlength\minimummorphismwidth{26mm}
    \node[morphism] (g) at (-0.15,1.5) {$=$};
    \draw (f.south)  to[out=-90,in=-180] ++(0.5,-.5) to[out=0,in=-90] ++(0.5,.5) to  ++(0,1.1) node[right] {$P_X$} to ([xshift=3pt]g.south east); 
    \draw (f.north) to[out=90,in=0] ++(-0.5,.5) to[out=180,in=90] ++(-0.5,-.5) to ([xshift=1pt]e.north east);
    \draw ([xshift=-1pt]e.north west) to +(0,0.2) node[left] {$x^*$} to ([xshift=-3pt]g.south west);
    \draw (g.north) to[out=90,in=180] ++(1,1) to[out=0,in=90] ++(1,-1) to ++(0,-4.5) to[out=-90,in=180] ++(0.5,-.5) to[out=0,in=-90] ++(0.5,.5) to ++(0,6);
    \draw (e.south) to[out=-90,in=0] ++(-0.5,-.5) to[out=180,in=-90] ++(-0.5,.5) to ++(0,4.825);
    \end{pic}
    \quad \cong \quad
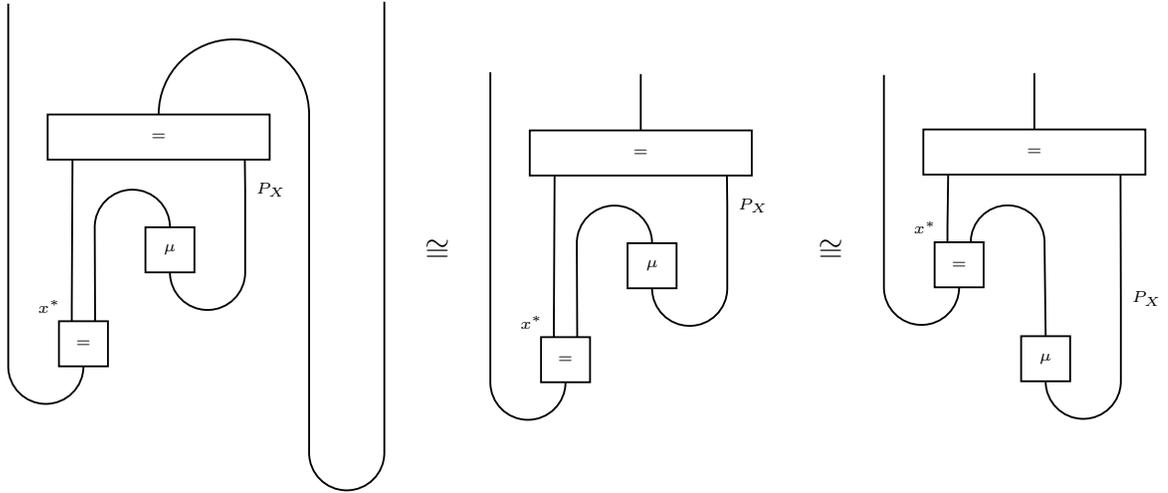
  \begin{pic} 
    \node[morphism] (f) at (0,0) {$\mu$};
    \node[morphism] (e) at (-1.15,-1.25) {$=$};
    \setlength\minimummorphismwidth{26mm}
    \node[morphism] (g) at (-0.15,1.5) {$=$};
    \draw (f.south)  to[out=-90,in=-180] ++(0.5,-.5) to[out=0,in=-90] ++(0.5,.5) to  ++(0,1.1) node[right] {$P_X$} to ([xshift=3pt]g.south east); 
    \draw (f.north) to[out=90,in=0] ++(-0.5,.5) to[out=180,in=90] ++(-0.5,-.5) to ([xshift=1pt]e.north east);
    \draw ([xshift=-1pt]e.north west) to +(0,0.2) node[left] {$x^*$} to ([xshift=-3pt]g.south west);
    \draw (g.north) to ++(0,.75);
    \draw (e.south) to[out=-90,in=0] ++(-0.5,-.5) to[out=180,in=-90] ++(-0.5,.5) to ++(0,4.125);
    \end{pic}
        \quad \cong \quad
  \begin{pic} 
    \node[morphism] (f) at (0,-1.25) {$\mu$};
    \node[morphism] (e) at (-1.15,0) {$=$};
    \setlength\minimummorphismwidth{26mm}
    \node[morphism] (g) at (-0.15,1.5) {$=$};
    \draw (f.south)  to[out=-90,in=-180] ++(0.5,-.5) to[out=0,in=-90] ++(0.5,.5) to  ++(0,1.1) node[right] {$P_X$} to ([xshift=3pt]g.south east); 
    \draw ([xshift=1pt]e.north east) to[in=180,out=90] ++(0.49,0.49) to[in=90,out=0] ++(0.49,-0.49) to  (f.north);
    \draw ([xshift=-1pt]e.north west) to +(0,0.2) node[left] {$x^*$} to ([xshift=-3pt]g.south west);
    \draw (g.north) to ++(0,.75);
    \draw (e.south) to[out=-90,in=0] ++(-0.5,-.5) to[out=180,in=-90] ++(-0.5,.5) to ++(0,2.825);
    \end{pic}
    \]
  To simplify this further, we observe first that $x\colon X\to Y$\/ induces a pseudonatural transformation 
  $\sigma(x)\colon L_X\to x^*L_Y$\/ whose $1$-cell component at an object $A\in\CC$\/ is given by 
      \[\begin{tikzpicture}
    \node (x) at (0,1) {$X$};
    \node (a) at (-1.5,-1) {$X+A$};
    \node (c) at (.75,0) {$Y$};
    \node (b) at (1.5,-1) {$Y+A$};
    \draw[->] (x) to coordinate[midway] (f) node[left] {$i_X$} (a);
    \draw[->] (x) to node[right] {$\ x$} coordinate[midway] (g)   (c);
    \draw[->] (c) to node[right] {$\ i_Y$} (b);
    \draw[->] (a) to node[below] {$x+A$} (b);
    \node at (0,-.254) {$\cong$};
    \end{tikzpicture}\]
  There is an invertible modification
  
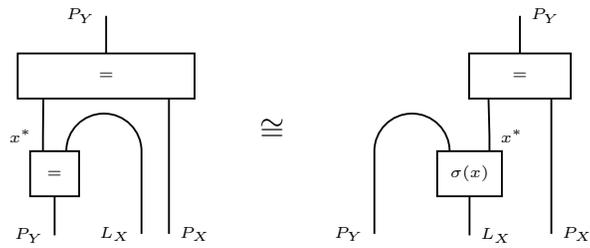
\[
  \begin{pic} 
    \node[morphism] (f) at (0,0) {$=$};
    \setlength\minimummorphismwidth{20mm}
    \node[morphism] (e) at (0.685,1.3) {$=$};
    \draw (f.south) to +(0,-.52) node[left] {$P_Y$};
    \draw ([xshift=-1pt]f.north west)  to +(0,0.2) node[left] {$x^*$} to ([xshift=-1pt]e.south west);
    \draw ([xshift=1pt]f.north east) to[out=90,in=180] ++(.5,.5) to[out=0,in=90] ++(.5,-.5) to ++(0,-1.1) node[left] {$L_X$};
    \draw (e.north) to +(0,.5) node[left] {$P_Y$};
    \draw ([xshift=1pt]e.south east) to +(0,-1.8) node[right] {$P_X$};
    \end{pic} \quad\cong\quad
    \begin{pic} 
    \node[morphism] (f) at (0,0) {$\sigma(x)$}; 
    \setlength\minimummorphismwidth{10mm} 
    \node[morphism] (e) at (0.67,1.3) {$=$};
    \draw (f.south) to +(0,-.52) node[right] {$L_X$};
    \draw ([xshift=-1.75pt]f.north west) to[out=90,in=0] ++(-.5,.5) to[out=180,in=90] ++(-.5,-.5) to ++(0,-1.1) node[left] {$P_Y$}; 
    \draw ([xshift=1.75pt]f.north east)  to +(0,0.2) node[right] {$x^*$} to ([xshift=-.5pt]e.south west);
    \draw ([xshift=.5pt]e.south east) to +(0,-1.8) node[right] {$P_X$};
    \draw (e.north) to +(0,.5) node[right] {$P_Y$};
    \end{pic}
  \]
  corresponding to the isomorphisms
  \[\begin{tikzpicture}
    \node (a) at (0,0) {$X+A$};
    \node (b) at (4,0) {$A$};
    \draw[draw=none] (a) to[out=35,in=145] coordinate[midway] (h1) (b);
    \node (c) at (h1) {$Y+A$};
    \node at (2,0) {$\cong$};
    \draw[->] (a) to[out=45,in=180] node[above,sloped] {$x+A$} (c);
    \draw[->] (c) to [out=0,in=135] node[above,sloped] {$[f,\id]$} (b);
    \draw[->] (a) to[out=-35,in=215] coordinate[midway] (h2) node[below] {$[fx,\id]$} (b);
    \end{tikzpicture}\]
  for $(A,f)\in\pslice{Y}{C}$, 
  where the top path corresponds to a $1$-cell of the RHS and the bottom path to a $1$-cell of the LHS.
  We can then continue the calculation and conclude that \ref{pic:mateofLxmu} is isomorphic to 
  \[  \begin{pic} 
    \node[morphism] (f) at (0,-1.25) {$\mu$};
    \node[morphism] (e) at (-1.15,0) {$=$};
    \setlength\minimummorphismwidth{26mm}
    \node[morphism] (g) at (-0.15,1.5) {$=$};
    \draw (f.south)  to[out=-90,in=-180] ++(0.5,-.5) to[out=0,in=-90] ++(0.5,.5) to  ++(0,1.1) node[right] {$P_X$} to ([xshift=3pt]g.south east); 
    \draw ([xshift=1pt]e.north east) to[in=180,out=90] ++(0.49,0.49) to[in=90,out=0] ++(0.49,-0.49) to  (f.north);
    \draw ([xshift=-1pt]e.north west) to +(0,0.2) node[left] {$x^*$} to ([xshift=-3pt]g.south west);
    \draw (g.north) to ++(0,.75);
    \draw (e.south) to[out=-90,in=0] ++(-0.5,-.5) to[out=180,in=-90] ++(-0.5,.5) to ++(0,2.825);
    \end{pic} 
    \quad\cong\quad
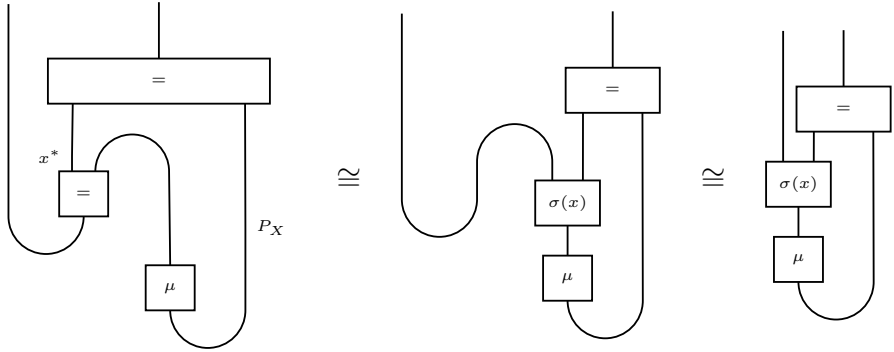
    \begin{pic} 
    \node[morphism] (f) at (0,0) {$\mu$};
    \node[morphism] (e) at (0,1) {$\sigma(x)$};
    \setlength\minimummorphismwidth{9mm}
    \node[morphism] (g) at (0.6,2.5) {$=$};
    \draw (e.north east) to ([xshift=-1pt]g.south west);
    \draw (e.north west) to ++(0,.25) to[out=90,in=0] ++(-.5,.5) to[in=90,out=180] ++(-.5,-.5) to ++(0,-.5)  to[out=-90,in=0] ++(-0.5,-.5) to[out=180,in=-90] ++(-0.5,.5) to ++(0,2.47);
    \draw (f.south) to[out=-90,in=180] ++(.5,-.5) to [out=0,in=-90] ++(.5,.5) to ([xshift=1pt]g.south east);
    \draw (f.north) to (e.south);
    \draw (g.north) to +(0,.75);
    \end{pic}
    \quad\cong\quad 
    \begin{pic} 
    \node[morphism] (f) at (0,0) {$\mu$};
    \node[morphism] (e) at (0,1) {$\sigma(x)$};
    \setlength\minimummorphismwidth{9mm}
    \node[morphism] (g) at (0.6,2) {$=$};
    \draw (e.north east) to ([xshift=-1pt]g.south west);
    \draw (e.north west) to ++(0,1.74);
    \draw (f.south) to[out=-90,in=180] ++(.5,-.5) to [out=0,in=-90] ++(.5,.5) to ([xshift=1pt]g.south east);
    \draw (f.north) to (e.south);
    \draw (g.north) to +(0,.75);
    \end{pic}
    \]
  Unwinding the definition of the right hand side, we see that $\NN'(x)\mu$\/ is isomorphic to $\NN(x)\mu$. 
  Moreover, this isomorphism is clearly natural in $\mu$\/ and it is monoidally natural as it respects composition of endomorphisms. 
  It is straightforward to check that the resulting action on $2$-cells is as desired.
\end{proof}

Inspired by a comment by an anonymous reviewer, we hereby conjecture that $\CC$ having lax limits of arrows is an alternative sufficient condition for the conclusion of Theorem~\ref{thm:isowithcoproducts} to hold. To explain why this is plausible, recall first that that the lax limit of an arrow $j\colon A\to B$ consists of a diagram 
\[\begin{tikzpicture}
\node (x) at (0,1) {$L(j)$};
\node (a) at (-1,-1) {$A$};
\node (b) at (1,-1) {$B$};
\draw[->] (x) to coordinate[midway] (f) node[left] {$\pi_1$} (a);
\draw[->] (x) to  coordinate[midway] (g)  node[right] {$\pi_2$} (b);
\draw[->] (a) to node[below] {$j$} (b);
\draw[shorttwocell] (f) to node[above] {$$} (g);
\end{tikzpicture}\]
that is universal as such, so that any other such triangle with $j$ on the bottom is induced by precomposing by a unique map into $L(j)$, together with a two-dimensional universal property that we do not need to recall here. In particular, if $(j,\sigma)$ defines a $1$-cell $(A,f)\to (B,g)$ in  $\lslice{X}{C}$, there is a unique map $\bar{\sigma}\colon X\to L(j)$ such that
\[ \begin{tikzpicture}[baseline=-.3cm]
    \node (x) at (0,1) {$X$};
    \node (a) at (-1,-1) {$A$};
    \node (b) at (1,-1) {$B$};
    \draw[->] (x) to coordinate[midway] (f) node[left] {$f$} (a);
    \draw[->] (x) to  coordinate[midway] (g)  node[right] {$g$} (b);
    \draw[->] (a) to node[below] {$j$} (b);
    \draw[shorttwocell] (f) to node[above] {$\sigma$} (g);
    \end{tikzpicture}\qquad=\qquad\begin{tikzpicture}[baseline=-.3cm] 
\node (y) at (0,2) {$X$};
\node (x) at (0,1) {$L(j)$};
\node (a) at (-1,-1) {$A$};
\node (b) at (1,-1) {$B$};
\draw[->] (y) to node[right] {$\bar{\sigma}$} (x);
\draw[->] (x) to coordinate[midway] (f) node[left] {$\pi_1$} (a);
\draw[->] (x) to  coordinate[midway] (g)  node[right] {$\pi_2$} (b);
\draw[->] (a) to node[below] {$j$} (b);
\draw[shorttwocell] (f) to node[above] {$$} (g);
\end{tikzpicture}\]
Similarly, the identity on $j$ induces an obvious diagonal map $\Delta\colon A\to L(j)$. Now, if  $\mu\in \Mps(X)$\/, we can use $\mu$ and the above maps in $\pslice{X}{C}$ to form the pasting diagram
\[
\begin{tikzpicture}[scale=1.25]
  \node (a) at (-2,0) {$A$};
   \node (b1) at (2,2) {$B$};
  \node (b2) at (2,-2) {$A$};
  \node (c1) at (6,-.5) {$A$};
  \node (c2) at (6,.5) {$B$};
  \node (d) at (0,0) {$L(j)$};
  \node (e) at (2,0) {$L(j)$};
  \draw[->] (a) to[in=180,out=45] node[above]  {$j$} (b1); 
  \draw[->] (a) to node[above]  {$\Delta$} (d);
  \draw[->] (a) to[in=180,out=-45] node[below]  {$\id$} (b2); 
  \draw[->] (d) to node[above]  {$\mu_{\bar{\sigma}}$} (e);
  \draw[->] (d) to node[sloped,above]  {$\pi_2$} (b1);
  \draw[->] (d) to node[sloped,below]  {$\pi_1$} (b2);
   \draw[->] (c1) to node[right]  {$j$} (c2);
  \draw[->] (e) to[out=35,in=145] coordinate[midway](mu) node[above]  {$\pi_2$} (c2);
  \draw[->] (e) to[out=-35,in=215] coordinate[midway](nu) node[below]  {$\pi_1$} (c1);
  \draw[shorttwocell] (nu) to  node[right] {$$} (mu);
  \draw[->] (b1) to[in=105,out=0]  node[above] {$\mu_{g}$} (c2);
  \draw[->] (b2) to[in=-105,out=0] node[below]  {$\mu_{f}$} (c1);
  \draw[shorttwocell] (e) to  node[right] {$\mu$} (b1);
  \draw[shorttwocell] (b2) to  node[right] {$\mu^{-1}$} (e);
\end{tikzpicture}
\]
defining a $2$-cell $j\mu_f\to \mu_g j$. In this way, any element $\mu$ of  $\Mps(X)$ can be used to obtain the data required of an element $\tilde{\mu}$ of $\MM(X)$. In fact, due to strictness of the limit notion, it is enough if $\mu$ is defined on the strict slice. Now, a $2$-cell $\theta j\to k$ induces a $1$-cell $L(k)\to L(j)$ (note the contravariance), and given composable $A\xrightarrow{j}B\xrightarrow{k} C$, we get comparison maps $L(j)\to L(k)$ and $L(kj)\to L(k)$: we expect that using such comparison maps would let one obtain a direct but diagrammatically involved proof that this extended $\tilde{\mu}$ is indeed an element of $\MM(X)$, and that the map $\mu\mapsto \tilde{\mu}$ provides a strict inverse for whiskering along $\pslice{X}{C}\hookrightarrow\lslice{X}{C}$. However, we have not discovered a conceptual proof for this conjecture, and as the required diagrams for the proposed direct argument get rather large and this result is not needed in the sequel, we have not pursued this conjecture further. This discussion motivates a stronger result that was conjectured by the reviewer: namely, that $\CC$ having powers (cotensors) by the walking arrow category $\cat{2}$ would be sufficient. This is indeed a stronger result as powers by $\cat{2}$ coincide with lax limits of identity morphisms, hence weakening the assumption. In this case, to obtain the required $2$-cell for $(j,\sigma)\colon (A,f)\to (B,g)$, one would first factorize it as $(A,f)\xrightarrow{(j,\id)} (B,jf)\xrightarrow{(\id,\sigma)} (B,g)$. Now, as lax natural transformations must respect composition and the arrow $(A,f)\xrightarrow{(j,\id)} (B,jf)$ already lives in $\pslice{X}{C}$, one can then use $L(\id[B])\cong B^{\cat{2}}$ as above to obtain a $2$-cell $\mu_{jf}\to \mu_g$, so that powers by $\cat{2}$ are already sufficient to produce the data required to extend an element of $\Mps(X)$ to $\MM(X)$. However, in this case we see neither a conceptual argument nor an obvious calculational path forward, leaving this conjecture wide open. 

\section{Isotropy and two-dimensional density}\label{sec:2ddensity}

We now seek to generalize the results about isotropy and dense functors from Section~\ref{sec:1d} to the two-dimensional setting.
As we are unaware of a bicategorical treatment of density (as opposed to the strict $\cat{Cat}$-enriched notion) 
we will adapt (i)-(iv) of \cite[Theorem 5.1]{kelly:enriched_reprint} to the bicategorical setting. 
The word ``cocontinuous'' below is to be understood in the bicategorical sense \ie as preserving bicolimits 
but not necessarily strict $2$-categorical limits: for instance, any left biadjoint is cocontinuous in this sense. 

\begin{proposition}\label{prop:2d-density} 
  Let $\CC$ and $\DD$ be bicategories. The following conditions are equivalent for a pseudofunctor $F\colon \CC\to \DD$. 
    \begin{enumerate}[label=(\roman*)]
      \item for all cocontinuous pseudofunctors  $G,H\colon \DD\to \EE$, restriction along $F$ induces an equivalence of categories $\hom(G,H)\to\hom(GF,HF)$;
      \item the restricted Yoneda embedding $\DD\xrightarrow{Y} [\DD\op,\cat{Cat}]\xrightarrow{\hat{F}}[\CC\op,\cat{Cat}]$\/ is a local equivalence;
      \item for all objects $A,B$ of $\DD$, the functor \[\DD(A,B)\to [\CC\op,\cat{Cat}](\DD(F-,A),\DD(F-,B))\] is an equivalence;
      \item for every object $A$ of $\DD$, the identity $2$-natural transformation $\DD(F-,A)\to \DD(F-,A)$ exhibits $A$ as the $\DD(F-,A)$-weighted colimit of $F$.
    \end{enumerate}
\end{proposition}

\begin{proof} 
  (i)$\Rightarrow$(ii):

  Any representable $\DD\op\to\cat{Cat}$ is continuous, so that the corresponding opposite pseudofunctor $\DD\to\cat{Cat}\op$ is cocontinuous. As the Yoneda embedding is local equivalence, (i) implies that the composite 
  \[ \DD\xrightarrow{Y} [\DD\op,\cat{Cat}]\xrightarrow{\hat{F}}[\CC\op,\cat{Cat}]\]
    is a local equivalence.

  (iii) is just a rephrasing of (ii), so that (ii) and (iii) are equivalent. The same holds for (iv), as unwinding definitions we see that (iv) amounts to saying that sending $f\colon A\to B$ to $\DD(F-,A)\xrightarrow{\DD(\id,f)}\DD(F-,B)$ gives an equivalence \[\DD(A,B)\simeq [\CC\op,\cat{Cat}](\DD(F-,A),\DD(F-,B))\] 
  that is pseudonatural in $B$: this implies (iii) and is clearly implied by (ii). 
  
  (iv)$\Rightarrow$(i): 
   Let $G,H\colon \DD\to \EE$ be cocontinuous pseudofunctors. Then $G$ preserves the colimits of (iv), so that for any $A\in\DD$ and $B\in\EE$, we have an equivalence
    \[\EE(GA,B)\simeq [\CC\op,\cat{Cat}](\DD(F-,A),\EE (GF-,B))\]
   that sends $f\colon GA\to B$ to the composite 
   \[\DD(F-,A)\xrightarrow{\hat{G}} \EE(GF-,GA)\xrightarrow{\EE(\id,f)}\EE(GF-,B).\]

   For a fixed $A$, setting $B=HA$ then gives an equivalence
  \[\EE(GA,HA)\simeq [\CC\op,\cat{Cat}](\DD(F-,A),\EE(GF-,HA))\]
  As every $A\in\DD$ is such a colimit, this will then hold uniformly in $A$, giving us an equivalence
  \[ [\DD,\EE](G,H)\simeq [\CC\op\times \DD,\cat{Cat}](\DD(F,\id),\EE(GF,H))\] 
  which sends $\tau\colon G\to H$ to the composite 
  \[\DD(F,\id)\xrightarrow{\hat{G}}\EE(GF,G)\xrightarrow{\EE(\id,\tau)}\EE(GF,H).\]

  Moreover, by~\cite[Proposition 8.6]{garnershulman:enriched}  there is always an equivalence
  \[[\CC,\EE](GF,HF)\simeq [\CC\op\times \DD,\cat{Cat}](\DD(F,\id),\EE(GF,H))\]
  given by sending $\tau\colon GF\to HF$ to 
\[\DD(F,\id)\xrightarrow{\hat{H}}\EE(HF,H)\xrightarrow{\EE(\tau,\id)}\EE(GF,H).\]
  As equivalences between $1$-categories satisfy the two-out-of-three property, to prove (i) it hence suffices to show that the triangle
    \[
      \begin{tikzpicture}
    \node (x) at (0,4) {$[\DD,\EE](G,H)$};
    \node (a) at (0,0) {$[\CC,\EE](GF,HF)$};
    \node (b) at (9.5,0) {$ [\CC\op\times \DD,\cat{Cat}](\DD(F,\id),\EE(GF,H))$};
    \draw[->] (x) to node[left] {$[F,\id]$} (a);
    \draw[->] (x) to node[above,sloped] {$\simeq$ by (iv) and cocontinuity of G} (b);
    \draw[->] (a) to node[below] {$\simeq$ by~\cite[Proposition 8.6]{garnershulman:enriched}} (b);
    \end{tikzpicture} 
    \]
  commutes up to natural isomorphism. Using the above explicit descriptions of the two equivalences, we see that to produce such an isomorphism is to produce for each pseudonatural transformation $\tau\colon G\to H$ an invertible modification 
    \[
    \begin{tikzpicture}
    \matrix (m) [matrix of math nodes,row sep=2em,column sep=4em,minimum width=2em]
    {
     \DD(F,\id) & \EE(GF,G)\\
      \EE(HF,H) & \EE(GF,H) \\};
    \path[->]
    (m-1-1) edge node [left] {$\hat{H}$} coordinate[midway] (a) (m-2-1)
           edge node [above] {$\hat{G}$} (m-1-2)
    (m-1-2) edge node [right] {$\EE(\id,\tau)$} (m-2-2)
    (m-2-1) edge node [below] {$\EE(\tau F,\id)$} coordinate[midway] (b) (m-2-2);
    \draw[shorttwocell](m-2-1) to node [above,sloped] {$\cong$} (m-1-2);
  \end{tikzpicture}
  \]
  that is natural in $\tau$. In turn, the data of such a modification consists of a natural isomorphism
    \[
    \begin{tikzpicture}
    \matrix (m) [matrix of math nodes,row sep=2em,column sep=4em,minimum width=2em]
    {
     \DD(FA,B) & \EE(GFA,GB)\\
      \EE(HFA,B) & \EE(GFA,HB) \\};
    \path[->]
    (m-1-1) edge node [left] {$H$} coordinate[midway] (a) (m-2-1)
           edge node [above] {$G$} (m-1-2)
    (m-1-2) edge node [right] {$\tau_B\circ (-)$} (m-2-2)
    (m-2-1) edge node [below] {$()\circ \tau_{FA}$} coordinate[midway] (b) (m-2-2);
    \draw[shorttwocell](m-2-1) to node [above,sloped] {$\cong$} (m-1-2);
  \end{tikzpicture}
  \]
  for each $A\in\CC\op,B\in\DD$. To unwind the final definition, such a natural isomorphism consists of an invertible 
  $2$-cell
    \[\begin{tikzpicture}
    \matrix (m) [matrix of math nodes,row sep=2em,column sep=4em,minimum width=2em]
    {
     GFA & GB \\
      HFA &  HB \\};
    \path[->]
    (m-1-1) edge node [left] {$\tau_{FA}$} coordinate[midway] (a) (m-2-1)
           edge node [above] {$G(f)$} (m-1-2)
    (m-1-2) edge node [right] {$\tau_B$} (m-2-2)
    (m-2-1) edge node [below] {$H(f)$} coordinate[midway] (b) (m-2-2);
    \draw[shorttwocell](m-2-1) to node [above,sloped] {$\cong$} (m-1-2);
  \end{tikzpicture}\] 
  for every $f\colon FA\to B$. For that, we may use the $2$-cell components of $\tau$. It is now straightforward but tedious to verify that pseudonaturality of $\tau$ guarantees that these data are suitably coherent, concluding the proof. 
\end{proof}

\begin{definition} A pseudofunctor is called $\emph{dense}$ if it satisfies the equivalent conditions of Proposition~\ref{prop:2d-density}. 
\end{definition}

The above notion might be more properly called bicategorical or $2$-categorical density, as 
opposed to strict $2$-categorical or $\cat{Cat}$-enriched density as discussed in~\cite[Theorem 5.1]{kelly:enriched_reprint}. 
Sometimes the notions coincide: for instance, the full subcategory on $\cat{1}$ is dense in $\cat{Cat}$\/
 in both senses, while the full subcategory on $\cat{0}$ is not dense in either sense. 
 We will discuss an example after Lemma~\ref{lem:stuffisdenseinmoncat} showing that $\cat{Cat}$-enriched density does not imply bicategorical density.
We also point out that, just as in the one-dimensional case, there is a weaker version satisfying a similar result for isotropy. However, we refrain from formulating that since
we will not need it for the present purposes.

We will now generalize Theorem~\ref{thm:densitycontrolsisotropy} and show that dense 
pseudofunctors let one compute isotropy like in the one-dimensional case.

\begin{corollary}\label{cor:2d-densitycontrolsisotropy}
  Let the $\CC$\/ be a $2$-category with binary coproducts and $K\colon \DD\to \CC$ be a dense pseudofunctor. 
  Then $\ZZ_\CC$ is pseudonaturally equivalent to a pseudofunctor $\LL_K$ that sends $X\in\CC$ to $\Aut(L_X\circ K)$\/ and that 
  on $1$- and $2$-cells acts as in Corollary~\ref{cor:LispseudoequivtoZ}. 
If, moreover, $K$\/ is small, then $\CC$\/ has essentially small 
 isotropy, in the sense that $\ZZ_\CC$\/ is pseudonaturally equivalent to a pseudofunctor taking values in the $2$-category of small $2$-groups. 
\end{corollary}
 
\begin{proof}
  This follows from Corollary~\ref{cor:LispseudoequivtoZ} and Proposition~\ref{prop:2d-density}.
\end{proof}

Using this, we can now compute various isotropy $2$-groups of interest.

\begin{theorem}\label{thm:groupoidshaveno2disotropy} The isotropy $2$-groups of groupoids vanish. 
More precisely, if $\cat{Grpd}$ is the $2$-category of groupoids, then $\ZZ\colon \cat{Grpd}\to \cat{2Grp}$ is 
pseudonaturally equivalent to the constant functor at the terminal $2$-group.
\end{theorem}

\begin{proof} 
  Let $\cat{Grpd}$ be the $2$-category of (small) groupoids, and let $\cat{1}$ be the terminal groupoid. 
  Then the full subcategory on $\cat{1}$ is $2$-categorically dense in $\cat{Grpd}$. 
  For any groupoid $\cat{G}$, the only autoequivalence of $\cat{G}+\cat{1}$ that fixes $\cat{G}$ (up to isomorphism) 
  has to fix $\cat{1}$ exactly. Hence $\LL_\cat{1}(G)$ and consequently $\ZZ(G)$ is equivalent to the trivial $2$-group.
\end{proof}

Note that this result really gives a pseudonatural equivalence and not an isomorphism. 
Indeed, given a groupoid $G$, one can obtain an element $\alpha\in\ZZ(G)$ that is isomorphic 
but not equal to $\id$ by choosing for each $f\colon G\to H$ an autoequivalence $\alpha_f\colon H\to H$\/
 and an isomorphism $\Gamma_f\colon \alpha_f\cong \id[H]$, as there is then a unique way of promoting 
 $\alpha$ into an element of $\ZZ(G)$ so that $\Gamma$ defines an invertible modification $\alpha\to\id$. 
 Thus $\ZZ(G)$ is a large $2$-group that is equivalent to the trivial $2$-group. Consequently, isotropy is trivial for 
 groupoids whether one organizes them into a $1$-category or into a $2$-category.
  However, if one looks at groupoids and \emph{cofunctors} between them, one gets a category with nontrivial isotropy, as shown in~\cite{garner:innerofgroupoids}. 

For the following application, we consider the $2$-category $[\CC\op,\cat{Cat}]$\/ of $\CC$-indexed categories. 

\begin{theorem}\label{thm:prestacks} For a $2$-category $\CC$, the isotropy $2$-functor of $[\CC\op,\cat{Cat}]$\/ is 
pseudonaturally equivalent to the constant functor with value $\Aut(\id[\CC])$. 
\end{theorem}

\begin{proof} The proof is similar to the one-dimensional case: the Yoneda embedding 
$y\colon\CC\to [\CC\op,\cat{Cat}]$\/ is dense so we can apply Theorem~\ref{cor:2d-densitycontrolsisotropy}. 
For a pseudofunctor $F\colon\CC\to\cat{Cat}$, an element of $\LL_y(F)$\/ is then given by autoequivalences 
$F+yA\to F+yA$ that restrict to the identity on $F$ and are pseudonatural in $A$. 
Consequently, they must map each $yA$\/ to itself, and hence by bicategorical Yoneda Lemma we have $\LL_y(F)\simeq \Aut(\id[\CC])$. 
\end{proof}

In particular, the isotropy of $\cat{Cat}$ is trivial. Moreover, for $\CC$\/ is a $1$-category and denoting by $\cat{Fib}(\CC)$\/ the $2$-category of cloven fibrations
over $\CC$, the isotropy $2$-functor $\ZZ\colon \cat{Fib}(\CC) \to \cat{2Grp}$\/ is constant, and equivalent to the group $\Aut(\id[\CC])$\/
viewed as a locally discrete $2$-group. 

Given that the covariant isotropy of any extensive category $\CC$ is constant at $\Aut(\id[\CC])$~\cite[Theorem 4.1]{parker:IsotropyofGrothendieckToposes}, we conjecture that the previous two theorems are special cases of a similar result for extensive $2$-categories.

\section{Monoidal categories}\label{sec:monoidal}

In this section we work out the motivating example of this paper, namely the Picard $2$-group of a monoidal category. 
We begin by setting up some notation and recalling some basic constructions on and facts about monoidal categories. 
We let $\cat{MonCat}$ denote the $2$-category of strict monoidal categories, 
strong monoidal functors and monoidal natural transformations---the strictness on objects is purely a matter 
of technical convenience as allowing for general monoidal categories results in a biequivalent $2$-category. 
The locally full subcategory on strict monoidal functors is denoted by  $\cat{MonCat}_s$. 

The forgetful functor $U\colon \cat{MonCat}_s\to\cat{Cat}$\/ has a left 2-adjoint $F\colon\cat{Cat}\to\cat{MonCat}_s$. 
For a given $\CC\in\cat{Cat}$, the monoidal category $F\CC$ can be described explicitly as follows: 
the objects and morphisms of $F\CC$\/ are finite sequences of objects and morphisms of $\CC$\/ respectively, 
with the domain and codomain of a sequence of morphisms given by the sequences of domains and codomains 
(so that $F\CC$\/ only has morphisms between sequences of the same length). 
Composition and identities in $F\cat{C}$\/ are induced from those of $\CC$\/ and the monoidal product is 
given by concatenation of sequences. In particular, the free monoidal category $F\cat{n}$\/
 for the discrete category $\cat{n}$\/ on $n$\/ objects results in a discrete monoidal category, 
 whose monoid on objects is given by the free monoid on $n$\/ generators. We will write $A_1,\dots A_n$\/ for the generators 
 of $F\cat{n}$, and on occasion will also denote the generator of $F\cat{1}$\/ merely by $A$, and the generators of $F\cat{2}$\/ by $A$\/ and $B$.

The universal property of $F\cat{n}$ states that a strict monoidal functor $F\cat{n}\to\CC$\/ is 
uniquely specified by a list of  $n$\/ objects of $\CC$---the images of the generators---and a 
monoidal natural transformation is uniquely determined by its components at the generators. 
We will routinely use this fact by denoting a morphism $F\cat{n}\to\CC$\/ sending the generator 
$A_i$ to $X_i\in\CC$\/ by $F\cat{n}\xrightarrow{A_i\mapsto X_i}\CC$\/ or simply by $F\cat{n}\xrightarrow{X_1,\dots, X_n}\CC$. 
 
In fact, $F\cat{n}$\/ (and $F\CC$\/ more generally) enjoys also a universal property for \emph{strong} 
monoidal functors out of $F\cat{n}$. More precisely, the inclusion 
\[\cat{MonCat}_s(F\CC,\DD)\to\cat{MonCat}(F\CC,\DD)\] is an equivalence of categories\footnote{This follows from~\cite[Corollary 5.6]{bkp:2dmonads} which states that free algebras for a $2$-monad are flexible.}.

Now, basic results in two-dimensional monad theory~\cite{bkp:2dmonads} imply that 
$\cat{MonCat}_s$\/ is cocomplete in the strict $\cat{Cat}$-enriched sense whereas $\cat{MonCat}$\/ has all bicolimits.  In particular,  $\cat{MonCat}$\/ has binary coproducts.  As we will only need coproducts of the form $\CC+F\cat{n}$\/ in our later computations, 
 we will only describe an explicit construction of the coproduct in this somewhat simpler case.

\paragraph{Construction of $\CC+F\cat{n}$}\label{para:coprodwithfreemonoidalcat}
 The underlying monoid on objects of $\CC+F\cat{n}$\/ is given by the coproduct of the corresponding monoids on $\CC$\/ and $F\cat{n}$. 
 As (the set of objects of) $F\cat{n}$\/ is a free monoid, every object of $\CC+F\cat{n}$\/ can hence be written as an alternating product 
 $X_0\otimes (\bigotimes_{i=1}^m A_{j_i}\otimes X_i)$, where $m\geq 0$, each $X_i$ is an object of $\CC$\/ (possibly the tensor unit) 
 and each $A_{j_i}$\/ is some generator of $F\cat{n}$. Moreover, as $F\cat{n}$\/ is discrete and the underlying monoid is free, 
 the number $m$\/ of generators of $F\cat{n}$\/ occurring in an object of $\CC+F\cat{n}$\/ is well-defined, 
 and consequently morphisms of $\CC+F\cat{n}$\/ are particularly simple to describe: there can be a morphism from  
 $X_0\otimes (\bigotimes_{i=1}^m A_{j_i}\otimes X_i)$ to $Y_0\otimes (\bigotimes_{i=1}^k A_{l_i}\otimes Y_i)$\/ only if
  $m=k$\/ and $A_{j_i}=A_{l_i}$\/ for each $i=1,\dots ,m$, in which case any such morphism is necessarily of the form 
  $f_0\otimes(\bigotimes_{i=1}^m \id\otimes f_i)$, with $f_i\colon X_i\to Y_i$ in $\CC$. 
  This is perhaps better explained via a string diagram: a generic morphism in $\CC+F\cat{n}$\/ is of the form 
  \[\begin{pic}
      \node[morphism] (f) at (0,1) {$f_0$};
      \draw (f.south) to +(0,-.7) node[left] {$X_0$};
      \draw (f.north) to +(0,.7) node[left] {$Y_0$};
      \draw (1,0) node[left] {$A_{j_1}$} to ++(0,2);
      \node[morphism] (g) at (2,1) {$f_1$};
      \draw (g.south) to +(0,-.7) node[left] {$X_1$};
      \draw (g.north) to +(0,.7) node[left] {$Y_1$};
     \node  at (3.1,1) {$\large\dots$};
      \draw (4,0) node[left] {$A_{j_n}$} to ++(0,2);
      \node[morphism] (h) at (5,1) {$f_n$};
      \draw (h.south) to +(0,-.7) node[left] {$X_n$};
      \draw (h.north) to +(0,.7) node[left] {$Y_n$};
  \end{pic}\]
for $n\geq 0$. Moreover, two parallel such morphisms are equal in $\CC+F\cat{n}$\/ if and only if their components in $\CC$\/ are equal, 
and this is true even if $\CC$\/ contains objects $X$\/ for which $X\otimes(-)$ in $\CC$\/ fails to be faithful (like the zero object of $\cat{Ab}$). In particular, when $n=0$, this description gives an isomorphism $\CC\cong \CC+F\cat{0}$. There are obvious strict monoidal functors from $\CC$ and $F\cat{n}$ into  $\CC+F\cat{n}$ making it into a (weak) coproduct in $\cat{MonCat}$. Moreover, $\CC+F\cat{n}$ as constructed above is also a \emph{strict coproduct} of $\CC$ and $F\cat{n}$ in $\cat{MonCat}_s$ so that one can uniquely specify strict monoidal functors out of it by specifying strict monoidal functors out of $\CC$ and out of $F\cat{n}$. In our case, we will often be working with strict monoidal functors $\CC+F\cat{n}\to \CC+F\cat{m}$ that restrict to the identity on $\CC$. As such maps are determined by their action on the generators of $F\cat{n}$, we will often abuse notation and write things like $\CC+F\cat{n}\xrightarrow{A_i\mapsto X_i} \CC+F\cat{m}$ instead of  $\CC+F\cat{n}\xrightarrow{[\id, A_i\mapsto X_i]} \CC+F\cat{m}$.

We now identify a suitable dense sub-$2$-category of $\cat{MonCat}$. 

\begin{lemma}\label{lem:stuffisdenseinmoncat}
  The inclusion $K$\/ of the full sub $2$-category $\mathbf{F}=\mathbf{F}_{0,1,2,3}$\/ of $\cat{MonCat}$\/ on the free monoidal categories on 
  $0,1,2$ and $3$\/ objects is dense in $\cat{MonCat}$. 
\end{lemma}

\begin{proof}
  By Proposition~\ref{prop:2d-density}, it suffices to verify for each $\CC,\DD\in\cat{MonCat}$, that the functor 
  \[ \mathcal{F}\colon \cat{MonCat}(\CC,\DD)\to [\mathbf{F},\cat{Cat}](\cat{MonCat}(K(-),\CC),\cat{MonCat}(K(-),\DD))\]
  is an equivalence. 
  
    First of all, a monoidal functor $G\colon\CC\to \DD$\/ induces a pseudonatural (in fact, a $2$-natural) transformation 
  \[ G\circ -\colon \cat{MonCat}(K(-),\CC)\to\cat{MonCat}(K(-),\DD),\] 
  and a monoidal natural transformation $\sigma\colon G\to H$\/ induces a modification $\sigma\circ -\colon G\circ -\to H\circ - $. 
  To see that the functor $\mathcal{F}$\/ is faithful, assume $\sigma\neq\tau$, so that $\sigma_X\neq \tau_X$\/ for some object $X$ of $\CC$. 
  Considering the unique strict monoidal functor $F\cat{1}\to\CC$\/ that maps the generator to $X$\/ we then see that $\sigma\circ -\neq\tau\circ -$. 
  To see that $\mathcal{F}$\/ is full, consider a modification $\Gamma\colon G\circ -\to H\circ - $. 
  In particular, the component of $\Gamma$ at $F\cat{1}$\/ gives a natural transformation 
  
  \[ \begin{tikzpicture}
    \node (a) at (-2,-1) {$\cat{MonCat}(F\cat{1},\CC)$};
    \node (b) at (3,-1) {$\cat{MonCat}(F\cat{1},\DD)$};
    \draw[->] (a) to[out=35,in=145] coordinate[midway] (h1) node[above]{$G\circ -$} (b);
    \draw[->] (a) to[out=-35,in=215] coordinate[midway] (h2) node[below] {$H\circ -$} (b);
    \draw[shorttwocell] (h1) to node[right] {$\Gamma_1$} (h2);
    \end{tikzpicture}\]
     that corresponds to a natural transformation $G\to H$. By considering maps $F\cat{2}\to \CC$, we can then verify that 
     $\sigma$ is a monoidal natural transformation so that $\Gamma=\sigma\circ -$.

  It remains to show that $\mathcal{F}$\/ is essentially surjective, so consider a pseudonatural transformation 
  \[ V\colon \cat{MonCat}(K(-),\CC)\to\cat{MonCat}(K(-),\DD).\] Using the universal property of $F\cat{1}$\/ 
  we can use $V$\/ to reconstruct a functor $G\colon\CC\to\DD$. More abstractly, we have equivalences 
  $\cat{MonCat}(F\cat{1},\cat{C})\simeq \cat{MonCat}_s(F\cat{1},\cat{C})\cong \cat{Cat}(1,U\cat{C})\cong U\cat{C}$, 
  so that using these equivalences a functor $\cat{MonCat}(F\cat{1},\cat{C})\to\cat{MonCat}(F\cat{1},\cat{D})$ gives rise 
  to a functor $G\colon U\cat{C}\to U\cat{D}$. Using the universal properties of $F\cat{0}$ and  $F\cat{2}$ we can endow this functor with
   the data of a strong monoidal functor. Finally, using the universal property of $F\cat{3}$ we can check the coherence 
   conditions and verify that the resulting functor is indeed strong monoidal.  It follows that $V\cong F\circ -$. 
\end{proof}

One would expect that the above result is deducible abstractly from the fact that the $2$-monad giving rise to $\cat{MonCat}$ has a specific finite presentation, much like for any algebraic theory that can be given by at most $n$-ary operations, the full subcategory on free models with $\leq n$ generators is dense. However, proving general results connecting presentations of $2$-monads to dense subcategories of algebras would take us too far afield. 

In passing, we note that we are now in a position to observe that $\cat{Cat}$-enriched density and bicategorical density are not equivalent. 
As any strong monoidal functor $F\cat{n}\to F\cat{m}$ is necessarily strictly monoidal, we may view  $\mathbf{F}$\/ both as a sub $2$-category of 
$\cat{MonCat}$\/ and of $\cat{MonCat}_s$. As a sub $2$-category of $\cat{MonCat}_s$, one can show that $\mathbf{F}$ is $\cat{Cat}$-dense. 
However, it is not dense in the bicategorical sense: due to flexibility of $F\cat{n}$, we have a 
pseudonatural equivalence $\cat{MonCat}(F\cat{m},-)\simeq \cat{MonCat}_s(F\cat{m},-)$,
 and consequently we can see that pseudonatural transformations 
 \[ \cat{MonCat}_s(K(-),\CC)\to \cat{MonCat}_s(K(-),\DD)\] 
 (when the domain of the functors is $\mathbf{F}$) correspond to strong monoidal functors by the above Lemma, 
 whereas $2$-natural transformations correspond to strict monoidal functors by $\cat{Cat}$-density. 
 As the categories of strong and strict monoidal functors $\CC\to\DD$ are not in general equivalent, 
 we see that  $\mathbf{F}$ is $\cat{Cat}$-dense in $\cat{MonCat}_s$\/ but not bicategorically dense.

We now turn to the main result of this section.
Let 
\[ \Picb\colon\cat{MonCat}\to\cat{2Grp}\]  
be the $2$-functor that sends a monoidal category $\CC$\/ to the $2$-group on the weakly invertible elements of $\CC$\/ 
and isomorphisms between them, and a strong monoidal functor to its restriction on such elements. 
As the component of a monoidal natural transformation between strong monoidal functors at any dualizable 
object is invertible\footnote{See~\cite[Proposition 5.2.3]{rivano:tannakiennes} for an early source, or~\cite[Proposition 7]{daypastro:frobenius} for a generalization of this to Frobenius monoidal functors.} this functor is also well-defined on $2$-cells.

\begin{theorem}\label{thm:picard} 
The isotropy $2$-group of monoidal categories is given by the Picard $2$-group. 
More precisely, there is a pseudonatural equivalence $\ZZ(-)\simeq \Picb\colon\cat{MonCat} \to \cat{2Grp}$.
\end{theorem}

Before giving the proof, we establish a technical lemma concerning the pseudonaturality squares of objects of $\LL_K(\CC)$.
Let $K\colon\mathbf{F} \to \cat{MonCat}$\/ be the dense inclusion from Lemma~\ref{lem:stuffisdenseinmoncat}, and let $\CC$\/ be a monoidal category. 
Consider an object $\alpha$\/ of $\LL_K(\CC)$. Explicitly, $\alpha$\/ comprises four monoidal functors $\alpha_i\colon \CC+F\cat{i} \to \CC+F\cat{i}$ and
for each $H\colon F\cat{n} \to F\cat{m}$\/ an invertible $2$-cell
   \[\begin{tikzpicture}
    \matrix (m) [matrix of math nodes,row sep=2em,column sep=4em,minimum width=2em]
    {
     \CC+F\cat{n} &  \CC+F\cat{m}\\
       \CC+F\cat{n} &  \CC+F\cat{m} \\};
    \path[->]
    (m-1-1) edge node [left] {$\alpha_n$} coordinate[midway] (a) (m-2-1)
           edge node [above] {$\id+H$} (m-1-2)
    (m-1-2) edge node [right] {$\alpha_m$} (m-2-2)
    (m-2-1) edge node [below] {$\id+H$} coordinate[midway] (b) (m-2-2);
    \draw[shorttwocell](m-2-1) to node [above,sloped] {$\alpha_H$} (m-1-2);
   \end{tikzpicture}\]
satisfying certain axioms.

\begin{lemma}\label{lem:psnat} Let $X$ be a weakly invertible element of $\CC\in\cat{MonCat}$, and let $\alpha$\/ be an object of $\LL_K(\CC)$  such that each $\alpha_i$ is a strict monoidal functor that strictly restricts to the identity under $\CC$, and satisfies  $\alpha_i(A_j)=X\otimes A_j\otimes X^{-1}$ for each generator $A_j$ of $F\cat{i}$.

 Then $\alpha$\/ is completely determined by the $2$-cells
$\alpha_H$\/ where $H$\/ ranges over the functors 
\[ \{F\cat{1}\to F\cat{0}\} \cup \{F\cat{1}\xrightarrow{A\mapsto A_i}F\cat{n}|n>1\}.\]
\end{lemma}
   In other words, if $\alpha'$\/ is another element of $\LL_K(\CC)$ that coincides with $\alpha$\/ on its $1$-cell components 
   (so that $\alpha'_i=\alpha_i$) and on the pseudonaturality squares for the monoidal functors 
   $\{F\cat{1}\to F\cat{0}\}\cup\{F\cat{1}\xrightarrow{A\mapsto A_i}F\cat{n}|n>1\}$,
    then $\alpha'=\alpha$. 

\begin{proof}
Throughout, we use the fact that to specify the $2$-cell $\alpha_H$\/ for some $H\colon F\cat{n} \to F\cat{m}$, it suffices to specify, for each generator $A_i$\/ of $F\cat{n}$,
the component $\alpha_{H,A_i}\colon (\id[\CC]+H)(\alpha_n(A_i)) \to \alpha_m(HA_i)$. This is because $\alpha_H$ 
 lives in $\pslice{\CC}{\cat{MonCat}}$ and hence is trivial under $\CC$. In particular, there is only one choice for the $2$-cell $\alpha_H$ when the domain of $H$ is $\CC+F\cat{0}\cong\CC$.

 As $F\cat{1}\to F\cat{0}\to F\cat{1}=F\cat{1}\xrightarrow{A\mapsto I} F\cat{1}$ and $\alpha$\/ is pseudonatural, we have the equation
     \[\begin{tikzpicture}[baseline=-.3cm]
    \matrix (m) [matrix of math nodes,row sep=2em,column sep=3em,minimum width=2em]
    {
     \CC+F\cat{1} & \CC+F\cat{0} &  \CC+F\cat{1}\\
      \CC+F\cat{1} & \CC+F\cat{0} &  \CC+F\cat{1}\\};
    \path[->]
    (m-1-1) edge node [left] {$\alpha_1$} coordinate[midway] (a) (m-2-1)
           edge node [above] {$$} (m-1-2)
    (m-1-2) edge node [right] {$\alpha_0$} (m-2-2)
            edge node [above] {$$} (m-1-3)
    (m-1-3) edge node [right] {$\alpha_1$} (m-2-3)
    (m-2-1) edge node [below] {$A\mapsto A_2$}  (m-2-2)
    (m-2-2) edge node [below] {$A,I$} (m-2-3);
    \draw[shorttwocell](m-2-1) to node [above,sloped] {$\alpha_I$} (m-1-2);
    \draw[shorttwocell](m-2-2) to node [above,sloped] {$\alpha$} (m-1-3);
  \end{tikzpicture}
    =
  \begin{tikzpicture}[baseline=-.3cm]
    \matrix (m) [matrix of math nodes,row sep=2em,column sep=3em,minimum width=2em]
    {
     \CC+F\cat{1} & \CC+F\cat{1}\\
      \CC+F\cat{1} & \CC+F\cat{1} \\};
    \path[->]
    (m-1-1) edge node [left] {$\alpha_1$} coordinate[midway] (a) (m-2-1)
           edge node [above] {$I$} (m-1-2)
    (m-1-2) edge node [right] {$\alpha_1$} (m-2-2)
    (m-2-1) edge node [below] {$I$} (m-2-2);
    \draw[shorttwocell](m-2-1) to node [above,sloped] {$\alpha_I$} (m-1-2);
  \end{tikzpicture}.\]
As the 2-cell components of $\alpha$ for maps out of $F\cat{0}$ are already fixed, this equation implies that the component of $\alpha$ at $F\cat{1}\to F\cat{0}$ determines the component of $\alpha$ at $F\cat{1}\xrightarrow{A\mapsto I} F\cat{1}$.

 Consider now the map $F\cat{2}\xrightarrow{A_1\mapsto A,A_2\mapsto I}F\cat{1}$. As $\alpha$\/ is pseudonatural, we have the equations
    \[\begin{tikzpicture}[baseline=-.3cm]
    \matrix (m) [matrix of math nodes,row sep=2em,column sep=3em,minimum width=2em]
    {
     \CC+F\cat{1} & \CC+F\cat{2} &  \CC+F\cat{1}\\
      \CC+F\cat{1} & \CC+F\cat{2} &  \CC+F\cat{1}\\};
    \path[->]
    (m-1-1) edge node [left] {$\alpha_1$} coordinate[midway] (a) (m-2-1)
           edge node [above] {$A\mapsto A_1$} (m-1-2)
    (m-1-2) edge node [right] {$\alpha_2$} (m-2-2)
            edge node [above] {$A,I$} (m-1-3)
    (m-1-3) edge node [right] {$\alpha_1$} (m-2-3)
    (m-2-1) edge node [below] {$A\mapsto A_1$}  (m-2-2)
    (m-2-2) edge node [below] {$A,I$} (m-2-3);
    \draw[shorttwocell](m-2-1) to node [above,sloped] {$\alpha_{A_1}$} (m-1-2);
    \draw[shorttwocell](m-2-2) to node [above,sloped] {$\alpha_{A,I}$} (m-1-3);
  \end{tikzpicture}
    =
  \begin{tikzpicture}[baseline=-.3cm]
    \matrix (m) [matrix of math nodes,row sep=2em,column sep=3em,minimum width=2em]
    {
     \CC+F\cat{1} & \CC+F\cat{1}\\
      \CC+F\cat{1} & \CC+F\cat{1} \\};
    \path[->]
    (m-1-1) edge node [left] {$\alpha_1$} coordinate[midway] (a) (m-2-1)
           edge node [above] {$\id$} (m-1-2)
    (m-1-2) edge node [right] {$\alpha_1$} (m-2-2)
    (m-2-1) edge node [below] {$\id$} (m-2-2);
    \node at (0,0) {$=$};
  \end{tikzpicture}\]
  and
    \[\begin{tikzpicture}[baseline=-.3cm]
    \matrix (m) [matrix of math nodes,row sep=2em,column sep=3em,minimum width=2em]
    {
     \CC+F\cat{1} & \CC+F\cat{2} &  \CC+F\cat{1}\\
      \CC+F\cat{1} & \CC+F\cat{2} &  \CC+F\cat{1}\\};
    \path[->]
    (m-1-1) edge node [left] {$\alpha_1$} coordinate[midway] (a) (m-2-1)
           edge node [above] {$A\mapsto A_2$} (m-1-2)
    (m-1-2) edge node [right] {$\alpha_2$} (m-2-2)
            edge node [above] {$A,I$} (m-1-3)
    (m-1-3) edge node [right] {$\alpha_1$} (m-2-3)
    (m-2-1) edge node [below] {$A\mapsto A_2$}  (m-2-2)
    (m-2-2) edge node [below] {$A,I$} (m-2-3);
    \draw[shorttwocell](m-2-1) to node [above,sloped] {$\alpha_{A_1}$} (m-1-2);
    \draw[shorttwocell](m-2-2) to node [above,sloped] {$\alpha_{A,I}$} (m-1-3);
  \end{tikzpicture}
    =
  \begin{tikzpicture}[baseline=-.3cm]
    \matrix (m) [matrix of math nodes,row sep=2em,column sep=3em,minimum width=2em]
    {
     \CC+F\cat{1} & \CC+F\cat{1}\\
      \CC+F\cat{1} & \CC+F\cat{1} \\};
    \path[->]
    (m-1-1) edge node [left] {$\alpha_1$} coordinate[midway] (a) (m-2-1)
           edge node [above] {$I$} (m-1-2)
    (m-1-2) edge node [right] {$\alpha_1$} (m-2-2)
    (m-2-1) edge node [below] {$I$} (m-2-2);
    \draw[shorttwocell](m-2-1) to node [above,sloped] {$\alpha_I$} (m-1-2);
  \end{tikzpicture}\]
  so that the component of $\alpha$\/ at $F\cat{2}\xrightarrow{A_1\mapsto A,A_2\mapsto I}F\cat{1}$\/ is 
  determined at the generators of $F\cat{2}$ by the above data, and hence determined altogether. 
  The same holds for $F\cat{2}\xrightarrow{A_1\mapsto I,A_2\mapsto A}F\cat{1}$. 

  We now show that the same holds for  $\alpha_{A\otimes B}$, \ie the component of $\alpha$\/ at 
  $F\cat{1}\xrightarrow{A\mapsto A_1\otimes A_2}F\cat{2}$. This $2$-cell is determined by its component at $A$\/ 
  which is a $1$-cell $X\otimes A_1\otimes A_2\otimes X^{-1}\to X\otimes A_1\otimes X^{-1}\otimes X\otimes A_2\otimes X^{-1}$\/ 
  in $\CC+F\cat{2}$. This $1$-cell is by \nameref{para:coprodwithfreemonoidalcat} of the form 
  \[\begin{pic}
      \node[morphism] (f) at (0,1) {$f$};
      \draw (f.south) to +(0,-.7) node[left] {$X$};
      \draw (f.north) to +(0,.7) node[left] {$X$};
      \draw (1,0) node[left] {$A_1$} to ++(0,2);
      \node[morphism,width=6mm] (g) at (2.5,1) {$g$};
      \draw ([xshift=-1pt]g.north west) to  +(0,.7) node[left] {$X^{-1}$};
      \draw ([xshift=1pt]g.north east) to  +(0,.7) node[right] {$X$};
      \draw (4,0) node[left] {$A_2$} to ++(0,2);
      \node[morphism] (h) at (5,1) {$h$};
      \draw (h.south) to +(0,-.7) node[left] {$X^{-1}$};
      \draw (h.north) to +(0,.7) node[left] {$X^{-1}$};
  \end{pic}\]
  for some $f,g,h$.

  Now, pseudonaturality of $\alpha$\/ implies the equations
      \[\begin{tikzpicture}[baseline=-.3cm]
    \matrix (m) [matrix of math nodes,row sep=2em,column sep=3em,minimum width=2em]
    {
     \CC+F\cat{1} & \CC+F\cat{2} &  \CC+F\cat{1}\\
      \CC+F\cat{1} & \CC+F\cat{2} &  \CC+F\cat{1}\\};
    \path[->]
    (m-1-1) edge node [left] {$\alpha_1$} coordinate[midway] (a) (m-2-1)
           edge node [above] {$A\otimes B$} (m-1-2)
    (m-1-2) edge node [right] {$\alpha_2$} (m-2-2)
            edge node [above] {$A,I$} (m-1-3)
    (m-1-3) edge node [right] {$\alpha_1$} (m-2-3)
    (m-2-1) edge node [below] {$A\otimes I$}  (m-2-2)
    (m-2-2) edge node [below] {$A,I$} (m-2-3);
    \draw[shorttwocell](m-2-1) to node [above,sloped] {$\alpha_{A\otimes B}$} (m-1-2);
    \draw[shorttwocell](m-2-2) to node [above,sloped] {$\alpha_{A,I}$} (m-1-3);
  \end{tikzpicture}
    =
  \begin{tikzpicture}[baseline=-.3cm]
    \matrix (m) [matrix of math nodes,row sep=2em,column sep=3em,minimum width=2em]
    {
     \CC+F\cat{1} & \CC+F\cat{1}\\
      \CC+F\cat{1} & \CC+F\cat{1} \\};
    \path[->]
    (m-1-1) edge node [left] {$\alpha_1$} coordinate[midway] (a) (m-2-1)
           edge node [above] {$\id$} (m-1-2)
    (m-1-2) edge node [right] {$\alpha_1$} (m-2-2)
    (m-2-1) edge node [below] {$\id$} (m-2-2);
    \node at (0,0) {$=$};
  \end{tikzpicture}\]
  and 
        \[\begin{tikzpicture}[baseline=-.3cm]
    \matrix (m) [matrix of math nodes,row sep=2em,column sep=3em,minimum width=2em]
    {
     \CC+F\cat{1} & \CC+F\cat{2} &  \CC+F\cat{1}\\
      \CC+F\cat{1} & \CC+F\cat{2} &  \CC+F\cat{1}\\};
    \path[->]
    (m-1-1) edge node [left] {$\alpha_1$} coordinate[midway] (a) (m-2-1)
           edge node [above] {$A\otimes B$} (m-1-2)
    (m-1-2) edge node [right] {$\alpha_2$} (m-2-2)
            edge node [above] {$I,A$} (m-1-3)
    (m-1-3) edge node [right] {$\alpha_1$} (m-2-3)
    (m-2-1) edge node [below] {$A\otimes I$}  (m-2-2)
    (m-2-2) edge node [below] {$I,A$} (m-2-3);
    \draw[shorttwocell](m-2-1) to node [above,sloped] {$\alpha_{A\otimes B}$} (m-1-2);
    \draw[shorttwocell](m-2-2) to node [above,sloped] {$\alpha_{I,A}$} (m-1-3);
  \end{tikzpicture}
    =
  \begin{tikzpicture}[baseline=-.3cm]
    \matrix (m) [matrix of math nodes,row sep=2em,column sep=3em,minimum width=2em]
    {
     \CC+F\cat{1} & \CC+F\cat{1}\\
      \CC+F\cat{1} & \CC+F\cat{1} \\};
    \path[->]
    (m-1-1) edge node [left] {$\alpha_1$} coordinate[midway] (a) (m-2-1)
           edge node [above] {$\id$} (m-1-2)
    (m-1-2) edge node [right] {$\alpha_1$} (m-2-2)
    (m-2-1) edge node [below] {$\id$} (m-2-2);
    \node at (0,0) {$=$};
  \end{tikzpicture}\]
  As these natural transformations are monoidal, evaluating these equations at the generating object $A$ of $F\cat{1}$ gives the equations
    \[\begin{pic}
      \draw[dashed] (-.5,1.5) to ++(0,1) to ++(3,0) to ++(0,-1) to ++(-3,0);
      \node at (-1,2) {$\alpha_{A,I}$};
      \node[morphism] (f) at (0,1) {$f$};
      \draw (f.south) to +(0,-.7) node[left] {$X$};
      \node[morphism] (i) at (0,2) {$$};
      \draw (f.north) to (i.north);
      \draw (i.north) to +(0,.5) node[left] {$X$};
      \draw (1,0) node[left] {$A$} to ++(0,2.8);
      \node[morphism,width=10mm] (g) at (2.42,1) {$g$};
       \node[morphism] (j) at (2,2) {$$};
      \node[morphism,width=10mm] (k) at (3.3,2) {$\alpha_{A,I}$};
      \draw ([xshift=-1pt]g.north west) to  (j.north);
      \draw (j.north) to +(0,.5) node[right] {$X^{-1}$};
      \draw ([xshift=1pt]g.north east) to  ([xshift=-1pt]k.north west);
      \node[morphism] (h) at (3.75,1) {$h$};
      \draw (h.south) to +(0,-.7) node[right] {$X^{-1}$};
      \draw (h.north) to ([xshift=1pt]k.north east);
  \end{pic}
  =
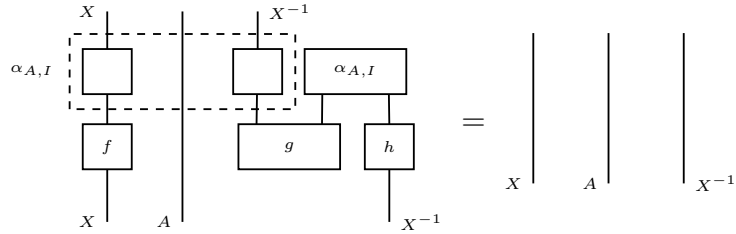
  \begin{pic}
    \draw (0,0) node[left] {$X$} to ++(0,2);
    \draw (1,0) node[left] {$A$} to ++(0,2);
    \draw (2,0) node[right] {$X^{-1}$} to ++(0,2);
  \end{pic}\]
  and
  \[
  \begin{pic}
      \draw[dashed] (1.25,1.5) to ++(0,1) to ++(3,0) to ++(0,-1) to ++(-3,0);
      \node at (4.75,2) {$\alpha_{I,A}$};
      \node[morphism] (f) at (0,1) {$f$};
      \draw (f.south) to +(0,-.7) node[left] {$X$};
      \node[morphism,width=10mm] (i) at (0.44,2) {$\alpha_{I,A}$};
      \draw (f.north) to ([xshift=-1pt]i.south west);
      \node[morphism,width=10mm] (g) at (1.31,1) {$g$};
      \node[morphism] (j) at (1.75,2) {$$};
       \draw (j.north) to +(0,.5) node[left] {$X$};
      \draw ([xshift=-1pt]g.north west) to ([xshift=1pt]i.south east);
      \draw ([xshift=1pt]g.north east) to (j.south);
      \draw (2.75,0) node[left] {$A$} to ++(0,2.8);
      \node[morphism] (h) at (3.75,1) {$h$};
      \draw (h.south) to +(0,-.7) node[right] {$X^{-1}$};
      \node[morphism] (k) at (3.75,2) {$$};
      \draw (h.north) to (k.north);
      \draw (k.north) to +(0,.5) node[right] {$X^{-1}$};
  \end{pic}
   = 
  \begin{pic}
    \draw (0,0) node[left] {$X$} to ++(0,2);
    \draw (1,0) node[left] {$A$} to ++(0,2);
    \draw (2,0) node[right] {$X^{-1}$} to ++(0,2);
  \end{pic}
  \]
  As $\alpha_{A,I}$ and $\alpha_{I,A}$ are invertible and have already been determined by the data at hand, these equations determine first $f$ and $h$ and consequently $g$, fixing $\alpha_{A\otimes B}$. 

  We now show that this is enough to determine the component of $\alpha$ at $F\cat{1}\xrightarrow{A\to W} F\cat{n}$ for any word $W$ in the generators of $F\cat{n}$, with the case $n=0$ holding by assumption. If $W=I$ and $n\geq 1$, then this already follows from pseudonaturality of $\alpha$ at $F\cat{1}\xrightarrow{A\to I} F\cat{n}=F\cat{1}\to F\cat{0}\to F\cat{n}$. If $W$ is a word of length 1, \ie a generator, then for $n>1$ it is fixed by assumption and for $n=1$ by the fact that $\alpha_{\id}=\id$ by pseudonaturality. Now, if $W=W_1\otimes W_2$ with components of $\alpha$ already determined at $A\mapsto W_1$ and at $A\mapsto W_2$, then the equality $F\cat{1}\xrightarrow{A\to A_i} F\cat{2}\xrightarrow{A_i\mapsto W_i} F\cat{n}=F\cat{1}\xrightarrow{A\to W_i}F\cat{n}$ and pseudonaturality of $\alpha$ determine the component of $\alpha$ at $F\cat{2}\xrightarrow{A_i\mapsto W_i} F\cat{n}$. As a result, pseudonaturality of $\alpha$ at $F\cat{1}\xrightarrow{A_1\otimes A_2} F\cat{2}\xrightarrow{A_i\mapsto W_i} F\cat{n}$ determines the component of $\alpha$ at $F\cat{1}\xrightarrow{A\to W} F\cat{n}$. Hence by induction the component of $\alpha$ is determined at $F\cat{1}\xrightarrow{A\to W} F\cat{n}$ for any word $W$ in the generators of $F\cat{n}$, and hence at any map $F\cat{1}\to F\cat{n}$. Consequently, for any map $F\cat{m}\to F\cat{n}$, the component of $\alpha$ at that map is determined for each generator of $F\cat{m}$, and hence altogether. Thus the pseudonaturality squares of $\alpha$ are uniquely determined by the $2$-cells $\{\alpha_{I}\}\cup\{\alpha_{n,A_i}|n>1\}$ corresponding to the monoidal functors $\{F\cat{1}\xrightarrow{A\mapsto I}F\cat{0}\}\cup\{F\cat{1}\xrightarrow{A\mapsto A_i}F\cat{n}|n>1\}$ as desired. 
\end{proof}

With this result, we can now prove Theorem~\ref{thm:picard}.

\begin{proof}[Proof of Theorem~\ref{thm:picard}]
  Invoking Lemma~\ref{lem:stuffisdenseinmoncat} and Corollary~\ref{cor:2d-densitycontrolsisotropy}, 
  it is sufficient to show that $\Picb \simeq \LL_K$, where $K\colon \mathbf{F} \to \cat{MonCat}$.
 Since there is a 2-equivalence between $2$-groups and coherent $2$-groups~\cite{baezlauda:2groups} \ie 
 $2$-groups where every object comes equipped with a chosen adjoint inverse, 
 we might as well assume that this is the case both for each $\Pic{C}$ and $\LL_K(\CC)$. 
 We now define a monoidal functor $\mathcal{G}=\mathcal{G}_\CC\colon\Pic{C}\to\LL_K(\CC)$ for each $\CC$. 
Consider $X\in\Pic{C}$\/ with inverse $X^{-1}$\/ and isomorphisms 
\[ i\colon I\to X^{-1}\otimes X \; \text{ and } \; e\colon X\otimes X^{-1}\to I \]
   satisfying the zigzag equations. We construct an element $\mathcal{G}(X)\in\LL_K(X)$\/ as follows.
    For $F\cat{n}$ with $n=0,1,2,3$, we define a strict monoidal functor $\mathcal{G}(X)_n$\/  by 
\[ \mathcal{G}(X)_n     =  [\id[\CC] , A_i\mapsto X\otimes A_i\otimes X^{-1}]\colon  \CC+F\cat{n}\to \CC+F\cat{n}.\]
     Given a monoidal functor  $H\colon F\cat{n}\to F\cat{m}$ (necessarily strict), the pseudonaturality square
   \[\begin{tikzpicture}
    \matrix (m) [matrix of math nodes,row sep=2em,column sep=4em,minimum width=2em]
    {
     \CC+F\cat{n} &  \CC+F\cat{m}\\
       \CC+F\cat{n} &  \CC+F\cat{m} \\};
    \path[->]
    (m-1-1) edge node [left] {$\mathcal{G}(X)_n$} coordinate[midway] (a) (m-2-1)
           edge node [above] {$\id+H$} (m-1-2)
    (m-1-2) edge node [right] {$\mathcal{G}(X)_m$} (m-2-2)
    (m-2-1) edge node [below] {$\id+H$} coordinate[midway] (b) (m-2-2);
    \draw[shorttwocell](m-2-1) to node [above,sloped] {\small $\mathcal{G}(X)_{H}$} (m-1-2);
   \end{tikzpicture}\]
   is given by the identity on $\CC$ and on the generator $A_i$  we define $\alpha_H$ by writing 
   $H(A_i)=\bigotimes_{j=1}^k W_j$ with each $W_j$ a generator of F\cat{m} and splitting into cases as follows:
      \begin{itemize}
        \item If $k=0$, \ie $H(A_i)=I$, then $(\id[\CC]+H)\circ\mathcal{G}_n(A_i)=X\otimes X^{-1}$\/ and 
        $\mathcal{G}_m\circ (\id[\CC]+H) (A_i)=I$, so we define $\mathcal{G}(X)_H(A_i)$ to be $e\colon  X \otimes X^{-1} \to I$.
        \item If $k=1$, so that $H(A_i)=A_j$ for some $j$, then 
        $(\id[\CC]+H)\circ\mathcal{G}_n(A_i)=X\otimes A_j\otimes X^{-1}=\mathcal{G}_m\circ (\id[\CC]+H) (A_i)$\/ and we define $\mathcal{G}(X)_H(A_i)$ to be the identity.
        \item If $k>1$, then $(\id[\CC]+H)\circ\mathcal{G}_n(A_i)=X\otimes (\bigotimes_{j=1}^n W_j)\otimes X^{-1}$\/ 
        and $\mathcal{G}_m\circ (\id[\CC]+H) (A_i)=\bigotimes_{j=1}^n (X\otimes W_j\otimes X^{-1})$, so we define $\mathcal{G}(X)_H(A_i)$\/
         by repeated applications of $i\colon I \to X^{-1} \otimes X$:
                \[\begin{pic}
                  \draw (0,0) node[left] {$X$} to ++(0,2);
                  \draw (1,0) node[left] {$W_1$} to ++(0,2);
                   \draw (2,2) node[left] {$X^{-1}$} to ++(0,-.5) to[in=180,out=-90] ++(.5,-.5) to[in=-90,out=0] ++(.5,.5) to ++(0,.5) node[right] {$X$};
                   \draw (4,0) node[left] {$W_2$} to ++(0,2);
                   \draw (5,2) node[left] {$X^{-1}$} to ++(0,-.5) to[in=180,out=-90] ++(.5,-.5) to[in=-90,out=0] ++(.5,.5) to ++(0,.5) node[right] {$X$};
                   \node  at (7,1.5) {$\large\dots$};
                  \draw (8,2) node[left] {$X^{-1}$} to ++(0,-.5) to[in=180,out=-90] ++(.5,-.5) to[in=-90,out=0] ++(.5,.5) to ++(0,.5) node[right] {$X$};
                   \draw (10,0) node[right] {$W_n$} to ++(0,2);
                   \draw (11,0) node[right] {$X^{-1}$} to ++(0,2);
                \end{pic}\]  
      \end{itemize}
   These are invertible morphisms because $i$\/ and $e$\/ are, and the zigzag equations imply that 
   this defines a pseudonatural transformation $\mathcal{G}(X)$. 
   Moreover, $\mathcal{G}(X)$ is an autoequivalence with inverse given by $\mathcal{G}(X^{-1})$.
   Any (by definition) invertible morphism $f\colon X\to Y$\/ in $\Pic{C}$\/ defines a modification 
   $\mathcal{G}(f)\colon \mathcal{G}(X)\to\mathcal{G}(Y)$, whose $n$-th $2$-cell component $\mathcal{G}(X)_n\to\mathcal{G}(Y)_n$\/
    is defined by identity on $\CC$ and by the diagram 
     \begin{equation}\label{eq:G1cell}
     \begin{pic} 
      \node[morphism] (f) at (0,0) {$f$};
      \draw (f.south) to +(0,-.85) node[left] {$X$};
      \draw (f.north) to +(0,.85) node[left] {$Y$};
      \draw (1,-1.15) node[left] {$A_j$} to ++(0,2.3);
      \node[morphism] (g) at (3,0) {$f^{-1}$};
      \draw (g.south)  to[in=0,out=-90] ++(-.5,-.5) to[in=-90,out=180] ++(-.5,.5) to ++(0,1.45)  node[right] {$Y^{-1}$};
      \draw (g.north)  to[in=180,out=90] ++(.5,.5) to[out=0,in=90] ++(.5,-.5) to  ++(0,-1.45) node[right] {$X^{-1}$};
      \end{pic}\end{equation}
    on the generator $A_i$\/ of $F\cat{n}$. 
    Hence $\mathcal{G}$\/ defines a functor $\Pic{C}\to\LL_K(\CC)$.  It is straightforward to check that that 
    $\mathcal{G}(X\otimes Y)\cong \mathcal{G}(Y)\circ\mathcal{G}(X)$, so that $\mathcal{G}$\/ is strong monoidal. 
    It is also straightforward, albeit more tedious, to verify that the resulting functors $\Pic{C}\to\LL_K(\CC)$\/ are pseudonatural in $\CC$.

   Hence it suffices to show that each functor $\mathcal{G}\colon\Pic{C}\to\LL_K(\CC)$ is fully faithful and essentially surjective.
   By virtue of how $\mathcal{G}$\/ acts on morphisms~\eqref{eq:G1cell} and by \nameref{para:coprodwithfreemonoidalcat}, we see that the functor $\mathcal{G}$\/ is faithful. 
   To see that $\mathcal{G}$\/ is full, consider a modification $\Gamma\colon \mathcal{G}(X)\to\mathcal{G}(Y)$. The component $\Gamma_1\colon \mathcal{G}(X)_1\to \mathcal{G}(Y)_1$\/ at the generating object $A$ of 
   $F\cat{1}$\/ is a map $X\otimes A\otimes X^{-1}\to Y\otimes A\otimes Y^{-1}$ and thus has to be induced by maps $f\colon X\to Y$ and $g\colon X^{-1}\to Y^{-1}$ by \nameref{para:coprodwithfreemonoidalcat}. Now, as $\Gamma$ is a modification we have the equation 
         \[\begin{tikzpicture}[baseline=-.3cm]
    \matrix (m) [matrix of math nodes,row sep=3em,column sep=4.5em,minimum width=2em]
    {
     \CC+F\cat{1} & \CC+F\cat{1}\\
    \CC+F\cat{1} &  \CC+F\cat{1} \\};
    \path[->]
    (m-1-1) edge node [left] {$\mathcal{G}(X)$} (m-2-1)
           edge node [above] {$I$} (m-1-2)
    (m-1-2) edge[in=135,out=-135] node [left] {$\mathcal{G}(X)$} coordinate[midway] (a) (m-2-2)
            edge[in=45,out=-45] node [right] {$\mathcal{G}(Y)$} coordinate[midway] (b) (m-2-2)
    (m-2-1) edge node [below] {$I$} (m-2-2);
    \begin{scope}[transform canvas={yshift=.4em}]
      \draw[shorttwocell,shorten >=20pt, shorten <=20pt](m-2-1) to node [above,sloped] {\small$\mathcal{G}(X)_I$} (m-1-2);
    \end{scope}
    \draw[shorttwocell] (a) to node[above] {$\Gamma_1$} (b);
  \end{tikzpicture}
  \enspace=\enspace
  \begin{tikzpicture}[baseline=-.3cm]
    \matrix (m) [matrix of math nodes,row sep=3em,column sep=4.5em,minimum width=2em]
    {
     \CC+F\cat{1} & \CC+F\cat{1}\\
    \CC+F\cat{1} &  \CC+F\cat{1} \\};
    \path[->]
    (m-1-1) edge[in=135,out=-135] node [left] {$\mathcal{G}(X)$} coordinate[midway] (a) (m-2-1)
            edge[in=45,out=-45] node [right] {$\mathcal{G}(Y)$} coordinate[midway] (b) (m-2-1)
           edge node [above] {$I$} (m-1-2)
    (m-1-2) edge node [right] {$\mathcal{G}(Y)$} (m-2-2)
    (m-2-1) edge node [below] {$I$} (m-2-2);
    \begin{scope}[transform canvas={yshift=-.4em}]
      \draw[shorttwocell,shorten >=20pt, shorten <=20pt](m-2-1) to node [below,sloped] {\small$\mathcal{G}(Y)_I$} (m-1-2);
    \end{scope}
    \draw[shorttwocell] (a) to node[above] {$\Gamma_1$} (b);
  \end{tikzpicture}\]
  which for the generating object of $F\cat{1}$\/ amounts to 
    \[\begin{pic} 
    \node[morphism] (f) at (0,0) {$f$};
     \node[morphism] (g) at (1,0) {$g$};
    \draw (f.south) to +(0,-.5) node[left] {$X$};
    \draw (g.south) to +(0,-.5) node[right] {$X^{-1}$};
    \draw (f.north) to[in=180,out=90] ++(.5,.5) to[in=90,out=0] (g.north);
    \end{pic}=  
    \begin{pic}
      \draw (-2,0) node[left] {$X$} to ++(0,1) to[in=180,out=90] ++(.5,.5) to[in=90,out=0] ++(.5,-.5) to ++(0,-1) node[right] {$X^{-1}$};
    \end{pic}\]  
   so that $g$\/ must satisfy
      \[\begin{pic} 
      \node[morphism] (f) at (0,0) {$g$};
      \draw (f.south) to +(0,-.85) node[left] {$X^{-1}$};
      \draw (f.north) to +(0,.85) node[left] {$Y^{-1}$};
      \end{pic}
      \enspace=\enspace
      \begin{pic} 
      \node[morphism] (g) at (3,0) {$f^{-1}$};
      \draw (g.south)  to[in=0,out=-90] ++(-.5,-.5) to[in=-90,out=180] ++(-.5,.5) to ++(0,1.45)  node[right] {$Y^{-1}$};
      \draw (g.north)  to[in=180,out=90] ++(.5,.5) to[out=0,in=90] ++(.5,-.5) to  ++(0,-1.45) node[right] {$X^{-1}$};
      \end{pic}\]
  and hence $\Gamma_1=\mathcal{G}(f)_1$. Moreover, $\Gamma_0=\id$ as this $2$-cell must be trivial over $\CC\cong \CC+F\cat{0}$, and
   $\Gamma_n$\/ for $n>1$\/ is uniquely determined by its action on the generators of $F\cat{n}$, 
  whence it is uniquely determined by whiskering it with the maps $\CC+F\cat{1}\xrightarrow{A_i} \CC+F\cat{n}$\/ sending the generator of $F\cat{1}$\/ to these generators, 
  so that $\Gamma_n$\/ is in fact determined by $\Gamma_1$. As $\Gamma_1=\mathcal{G}(f)_1$\/ by construction, we have $\Gamma=\mathcal{G}(f)$\/ as desired.

  We next show that $\mathcal{G}$\/ is essentially surjective, so let $\alpha\in \LL_K(\CC)$.  
 Now, $\alpha_1$ is determined within monoidal isomorphism by $\alpha_1(A)$, which by \nameref{para:coprodwithfreemonoidalcat} is given by a word in $w(A)$\/
  in the objects of $\CC$\/ and $A$. Because $\alpha$\/ is pseudonatural, for each object $Y$\/ of $F\cat{2}$, the square
  \[
  \begin{tikzpicture}
     \matrix (m) [matrix of math nodes,row sep=3em,column sep=4.5em,minimum width=2em]
    {
     \CC+F\cat{1} & \CC+F\cat{2}\\
    \CC+F\cat{1} &  \CC+F\cat{2} \\};
    \path[->] 
   (m-1-1) edge node[left] {$\alpha_1$} (m-2-1)
   (m-2-1) edge node[below] {$\CC+Y$} (m-2-2)
   (m-1-1) edge node[above] {$\CC+Y$} (m-1-2)
   (m-1-2) edge node[right] {$\alpha_2$} (m-2-2);
  \end{tikzpicture}
  \]
 must commute up to natural isomorphism. This means that 
   $\alpha_2$\/ acts, up to isomorphism, by $w(-)$\/ \ie by replacing instances of $A$\/ in the word $w(A)$\/ with the corresponding object of $F\cat{n}$. 
   The same reasoning applies to $\alpha_3$. Hence we assume (replacing $\alpha_n$\/ by a monoidally isomorphic functor if necessary) 
   that $\alpha_n$\/ is exactly given by $w(-)$\/ on $F\cat{2}$. As $\alpha_1$\/ is an equivalence, this word has to contain $A$\/ at least once. 
   As $\alpha_2$\/ is strong monoidal, we must have $w(A)\otimes w(B)\cong w(A\otimes B)$. 
   Now, each instance of $A$\/ is to the left of each instance of $B$\/ in the word $w(A)\otimes w(B)$, and hence this has to hold true of  
   $w(A\otimes B)$. Hence $w(A)$ has to contain $A$\/ exactly once, so that $w(A)=X\otimes A\otimes Y$\/ for some $X,Y\in\CC$. 
   The strong monoidal structure of $\alpha_2$\/ then gives us a morphism 
   \[X\otimes A\otimes Y\otimes X\otimes B\otimes Y= w(A)\otimes w(B)\to w(A\otimes B)=X\otimes A\otimes B\otimes Y\]
   which by \nameref{para:coprodwithfreemonoidalcat} induces an isomorphism  $Y\otimes X\to I$\/ in $\CC$.  Moreover, the unit constraint of $\alpha_2$ gives us an isomorphism $I\to w(I)=X\otimes I\otimes Y=X\otimes Y$ which again must live in $\CC$. As a result, the object $X$ is weakly invertible in $\CC$ and  $Y\cong X^{-1}$. We may hence assume that $\alpha_n$\/ is equal to the strict monoidal functor $\mathcal{G}(X)_n$. 

  However, we are not yet finished as the pseudonaturality squares of $\alpha$ may differ from those of $\mathcal{G}(X)$. We will produce an invertible modification $\Gamma\colon\alpha\to \mathcal{G}(X)$ using the following strategy. First of all, any invertible $\Gamma_0,\Gamma_1,\Gamma_2$ and $\Gamma_3$ will induce a modification $\Gamma\colon \alpha\to\beta$ for a unique pseudonatural $\beta$. We will proceed by finding some $\Gamma_0,\Gamma_1,\Gamma_2,\Gamma_3$ for which we can conclude that $\beta=\mathcal{G}(X)$ using Lemma~\ref{lem:psnat}.  This means that it suffices to find $\Gamma_0,\Gamma_1,\Gamma_2,\Gamma_3$ so that the resulting $\beta$ agrees with $\mathcal{G}(X)$ on the pseudonaturality $2$-cells required by Lemma~\ref{lem:psnat}, and conversely, we can use the desired equations in order to find $\Gamma$.

  For $\Gamma_0$, as this  $2$-cell must live in $\pslice{\CC}{\cat{MonCat}}$, we can only set it to be the identity so we do so. Our first nontrivial step is to find an isomorphism $\Gamma_1\colon \alpha_1=\mathcal{G}(X)_1\to\mathcal{G}(X)_1$ satisfying the pasting diagram 
    \[\begin{tikzpicture}[baseline=-.3cm]
    \matrix (m) [matrix of math nodes,row sep=3em,column sep=4.5em,minimum width=2em]
    {
     \CC+F\cat{1} & \CC+F\cat{0}\\
    \CC+F\cat{1} &  \CC+F\cat{0} \\};
    \path[->]
    (m-1-1) edge node [left] {$\mathcal{G}(X)$} (m-2-1)
           edge node [above] {$I$} (m-1-2)
    (m-1-2) edge[in=135,out=-135] node [left] {$\mathcal{G}(X)$} coordinate[midway] (a) (m-2-2)
            edge[in=45,out=-45] node [right] {$\mathcal{G}(X)$} coordinate[midway] (b) (m-2-2)
    (m-2-1) edge node [below] {$I$} (m-2-2);
    \begin{scope}[transform canvas={yshift=.4em,xshift=-.4em}]
      \draw[shorttwocell,shorten >=20pt, shorten <=20pt](m-2-1) to node [above,sloped] {$\alpha_I$} (m-1-2);
    \end{scope}
    \draw[shorttwocell] (a) to node[above] {$\Gamma_0$} (b);
  \end{tikzpicture}
  \enspace=\enspace
  \begin{tikzpicture}[baseline=-.3cm]
    \matrix (m) [matrix of math nodes,row sep=3em,column sep=4.5em,minimum width=2em]
    {
     \CC+F\cat{1} & \CC+F\cat{0}\\
    \CC+F\cat{1} &  \CC+F\cat{0} \\};
    \path[->]
    (m-1-1) edge[in=135,out=-135] node [left] {$\mathcal{G}(X)$} coordinate[midway] (a) (m-2-1)
            edge[in=45,out=-45] node [right] {$\mathcal{G}(X)$} coordinate[midway] (b) (m-2-1)
           edge node [above] {$I$} (m-1-2)
    (m-1-2) edge node [right] {$\mathcal{G}(X)$} (m-2-2)
    (m-2-1) edge node [below] {$I$} (m-2-2);
    \begin{scope}[transform canvas={yshift=-.4em,xshift=.4em}]
      \draw[shorttwocell,shorten >=20pt, shorten <=20pt](m-2-1) to node [below,sloped] {\small $\mathcal{G}(X)_I$} (m-1-2);
    \end{scope}
    \draw[shorttwocell] (a) to node[above] {$\Gamma_1$} (b);
  \end{tikzpicture}\]
  This is to ensure that the resulting pseudonatural transformation $\beta$ agrees with $\mathcal{G}(X)$ for $F\cat{1}\to F\cat{0}$.

  Now, any automorphism of $\mathcal{G}(X)_1$ in $\pslice{\CC}{\cat{MonCat}}$ is induced by an automorphism of  $X\otimes A\otimes X^{-1}$ which by \nameref{para:coprodwithfreemonoidalcat} is necessarily of the form $f\otimes \id[A]\otimes g$. Moreover, we have already fixed $\Gamma_0$ to be the identity. Thus finding such a $\Gamma_1$ amounts to finding invertible $f,g$ such that 
    \[\begin{pic} 
    \node[morphism] (f) at (0,0) {$f$};
     \node[morphism] (g) at (1,0) {$g$};
    \draw (f.south) to +(0,-.5) node[left] {$X$};
    \draw (g.south) to +(0,-.5) node[right] {$X^{-1}$};
    \draw (f.north) to[in=180,out=90] ++(.5,.5) to[in=90,out=0] (g.north);
    \end{pic}=  
    \begin{pic}
      \node[morphism,width=8mm] (f) at (0,0) {$\alpha_I$};
      \draw([xshift=-1pt]f.south west) to +(0,-.5) node[left] {$X$};
      \draw([xshift=1pt]f.south east) to +(0,-.5) node[right] {$X^{-1}$};
    \end{pic}\]
  As 
    \[\begin{pic}
      \node[morphism,width=8mm] (f) at (0,0) {$\alpha_I$};
      \draw(f.south west) to[in=0,out=-90] ++(-.5,-.5) to[in=-90,out=180] ++(-.5,.5) to ++(0,1) to[in=0,out=90] ++(-.5,.5) to[in=90,out=180] ++(-.5,-.5) to ++(0,-2) node[left] {$X$};
      \draw(f.south east) to +(0,-1) node[right] {$X^{-1}$};
    \end{pic}=
    \begin{pic}
      \node[morphism,width=8mm] (f) at (0,0) {$\alpha_I$};
      \draw([xshift=-1pt]f.south west) to +(0,-.5) node[left] {$X$};
      \draw([xshift=1pt]f.south east) to +(0,-.5) node[right] {$X^{-1}$};
    \end{pic}\]
  this can be solved by setting $f=\id[X]$ and defining $g$ via
  \[\begin{pic}
      \node[morphism,width=8mm] (f) at (0,0) {$\alpha_I$};
      \draw(f.south west) to[in=0,out=-90] ++(-.5,-.5) to[in=-90,out=180] ++(-.5,.5) to ++(0,1) node[left] {$X^{-1}$} ;
      \draw(f.south east) to +(0,-.75) node[right] {$X^{-1}$};
    \end{pic}\]
    Now that we've defined $\Gamma_1$, consider the desired pasting equalities
        \[\begin{tikzpicture}[baseline=-.3cm]
    \matrix (m) [matrix of math nodes,row sep=3em,column sep=4.5em,minimum width=2em]
    {
     \CC+F\cat{1} & \CC+F\cat{n}\\
    \CC+F\cat{1} &  \CC+F\cat{n} \\};
    \path[->]
    (m-1-1) edge node [left] {$\mathcal{G}(X)$} (m-2-1)
           edge node [above] {$A_i$} (m-1-2)
    (m-1-2) edge[in=135,out=-135] node [left] {$\mathcal{G}(X)$} coordinate[midway] (a) (m-2-2)
            edge[in=45,out=-45] node [right] {$\mathcal{G}(X)$} coordinate[midway] (b) (m-2-2)
    (m-2-1) edge node [below] {$A_i$} (m-2-2);
    \begin{scope}[transform canvas={yshift=.4em,xshift=-.4em}]
      \draw[shorttwocell,shorten >=20pt, shorten <=20pt](m-2-1) to node [above,sloped] {\small$\alpha_{A_i}$} (m-1-2);
    \end{scope}
    \draw[shorttwocell] (a) to node[above] {$\Gamma_n$} (b);
  \end{tikzpicture}
  \enspace=\enspace
  \begin{tikzpicture}[baseline=-.3cm]
    \matrix (m) [matrix of math nodes,row sep=3em,column sep=4.5em,minimum width=2em]
    {
     \CC+F\cat{1} & \CC+F\cat{n}\\
    \CC+F\cat{1} &  \CC+F\cat{n} \\};
    \path[->]
    (m-1-1) edge[in=135,out=-135] node [left] {$\mathcal{G}(X)$} coordinate[midway] (a) (m-2-1)
            edge[in=45,out=-45] node [right] {$\mathcal{G}(X)$} coordinate[midway] (b) (m-2-1)
           edge node [above] {$A_i$} (m-1-2)
    (m-1-2) edge node [right] {$\mathcal{G}(X)$} (m-2-2)
    (m-2-1) edge node [below] {$A_i$} (m-2-2);
    \begin{scope}[transform canvas={yshift=-.4em,xshift=.4em}]
      \draw[shorttwocell,shorten >=20pt, shorten <=20pt](m-2-1) to node [below,sloped] {\small$\mathcal{G}(X)$} (m-1-2);
    \end{scope}
    \draw[shorttwocell] (a) to node[above] {$\Gamma_1$} (b);
  \end{tikzpicture}\]
  for maps $F\cat{1}\xrightarrow{A\mapsto A_i}F\cat{n}$ for generators $A_i$ of $F\cat{n}$ with $n>1$. As $\alpha$ is invertible, these determine the components of $\Gamma_2,\Gamma_3$ uniquely at the generators, and hence fix unique monoidal natural transformations $\Gamma_2,\Gamma_3$ satisfying these equations. 

  It remains to show that $(\Gamma_{i})_{i=0}^3$ defines a modification $\Gamma\colon\alpha\to \mathcal{G}(X)$, \ie 
  that the desired equations are satisfied for all maps $F\cat{m}\to F\cat{n}$. However, this follows from our earlier observation: 
  \/given  $(\Gamma_{i})_{i=0}^3$, there is a unique pseudonatural transformation $\beta$\/ such that 
  $(\Gamma_{i})_{i=0}^3$ defines a modification $\Gamma\colon\alpha\to \beta$. By construction, $\beta$\/ agrees with 
  $\mathcal{G}(X)$\/ at its $1$-cell components and at the $2$-cell components corresponding to the monoidal functors 
  $\{F\cat{1}\xrightarrow{A\mapsto I}F\cat{0}\}\cup\{F\cat{1}\xrightarrow{A\mapsto A_i}F\cat{n}|n>1\}$, so that we must have 
  $\beta=\mathcal{G}(X)$ by Lemma~\ref{lem:psnat}. Thus $\mathcal{G}$ is essentially surjective as desired.
\end{proof}

\section{Monoidal fibrations}\label{sec:monfib}

In this final section we combine two of the earlier results by characterizing the 2-isotropy of the $2$-category of pseudofunctors 
$[\cat{C}\op,\cat{MonCat}]$, pseudonatural transformations and modifications. 
When $\cat{C}$\/ is locally discrete, this $2$-category is equivalent to the category of \emph{monoidal fibrations} over $\CC$\/ in the sense 
of~\cite{moeller:monoidalgrothendieck}. 
One might expect that the result is deducible directly from Theorems~\ref{thm:prestacks} and~\ref{thm:picard} 
given the similar but more general results proven in the one-dimensional case in~\cite{parker:presheaves}. 
However, as a full understanding of isotropy (both in one and two dimensions) in functor categories is outside the scope of this paper, 
we will proceed by a direct proof, highlighting the interesting components while merely sketching the repetitive parts of the arguments.

We begin by identifying a suitable dense sub-$2$-category of $[\cat{C}\op,\cat{MonCat}]$. Note that we have a diagram
\[
\begin{tikzpicture}
 \matrix (m) [matrix of math nodes,row sep=3em,column sep=4.5em,minimum width=2em]
    {
     {\cat{C} } & { [\cat{C}\op,\cat{Cat}] } & { [\cat{C}\op,\cat{MonCat}] } \\ };
    \path[->]
    (m-1-1) edge node [above] {$Y$} coordinate[midway] (a) (m-1-2)
     (m-1-2) edge node [above] {$F \circ -$} coordinate[midway] (b) (m-1-3);
\end{tikzpicture}
\]
where $Y$\/ is the Yoneda embedding and $F\colon\cat{Cat} \to \cat{MonCat}$\/ is the free monoidal category functor as before. We will write $\overline{A}$\/ for
the image of an object $A$\/ of $\CC$\/ under this composite functor. Thus $\overline{A}\colon\CC\op \to \cat{MonCat}$\/ acts on objects by $\overline{A}(C)=
F(\CC(C,A))$, the free monoidal category on the hom-category $\CC(C,A)$.

\begin{lemma}\label{lem:stuffisdenseinfunctorstomoncat} 
The full sub $2$-category $\DD$  of $[\cat{C}\op,\cat{MonCat}]$ spanned by the objects $\{\overline{A}, \overline{A}+\overline{A}, 
\overline{A}+\overline{A}+\overline{A}  \mid A\in\CC\}\cup\{\Delta (F\cat{1})\}$ is dense.
\end{lemma}

\begin{proof} The proof proceeds along the lines of Lemma~\ref{lem:stuffisdenseinmoncat}---indeed, 
when $\CC$ is the terminal category, we obtain Lemma~\ref{lem:stuffisdenseinmoncat} as a special case. 
Let us denote the inclusion  $\DD\hookrightarrow [\CC\op,\cat{MonCat}]$\/ by $\mathcal{F}$. 
Fixing pseudofunctors $G,H\colon\CC\to\DD$, it suffices to show that the functor 
\[ [\CC\op,\cat{MonCat}](G,H)\to [\DD\op,\cat{Cat}]([\CC\op,\cat{MonCat}](\mathcal{F}-,G),[\CC\op,\cat{MonCat}](\mathcal{F}-,H))\] 
is an equivalence. To see that it is faithful, consider two distinct modifications $\Gamma,\Gamma'$\/ between 
pseudonatural transformations $G\to H$. As $\Gamma\neq \Gamma'$, we must have $\Gamma_A\neq \Gamma'_A$\/
 for some $A\in\CC$. Consequently, the components of the induced modifications differ at $\overline{A}$, showing faithfulness. 
 To see that this functor is full, consider a modification $\Gamma$ between the pseudonatural transformations induced 
 by $\sigma,\tau\colon G\to H$. Now, the universal property of each $\overline{A}$ lets us recover a natural transformation 
 $\Gamma'_A\colon \sigma_A\to\tau_A$\/ for each $A$, and the universal property of $\overline{A}+\overline{A}$ lets us verify that 
 these natural transformations are monoidal. Consequently, $\Gamma'$ defines a modification $\sigma\to\tau$\/ and maps to $\Gamma$, showing fullness.

To see that the functor is essentially full, consider a pseudonatural transformation 
\[ \sigma\colon [\CC\op,\cat{MonCat}](\mathcal{F}-,G)\to [\CC\op,\cat{MonCat}](\mathcal{F}-,H)\]
 in $[\DD\op,\cat{Cat}]$. 
Now, we have a pseudonatural equivalence
\[ [\CC\op,\cat{Cat}](y-,UG)  \simeq  [\CC\op,\cat{MonCat}]((F \circ -)y-,G)\]
so that restricting to objects of the form $\overline{A}$, we obtain a pseudonatural transformation 
$[\CC\op,\cat{Cat}](y-,UG)\to [\CC\op,\cat{Cat}](y-,UH)$\/ which by Yoneda corresponds to a pseudonatural transformation 
$UG\to UH$. Hence it suffices to show that that the components of this pseudonatural transformation live in 
$\cat{MonCat}$. Considering objects of the form $\Delta (F\cat{1})$ and $\overline{A}+\overline{A}$ lets us endow the $1$-cell components with the data 
of a strong monoidal functor, and considering objects of the form  $\overline{A}+\overline{A}+\overline{A}$ lets us check that these are 
in fact strong monoidal. Considering objects of the form $\overline{A}+\overline{A}$\/ again lets us then verify that the $2$-cell components 
are monoidal natural transformations. Taken together, this lets us find a pseudonatural transformation $G\to H$\/ that induces $\sigma$ within isomorphism.
\end{proof}

The description of the isotropy $2$-group of a monoidal fibration uses the notion of a (conical) 
pseudolimit of a pseudofunctor $H\colon \CC\op \to \cat{2Grp}$. This notion is defined in the general case in~\cite[Section 5.1]{johnsonyau:2d} 
and in~\cite{kelly:2-dlimits}. Such pseudolimits are inherited from $\cat{MonCat}$, and may be explicitly described as follows.

\begin{description}
\item [Objects] of $\lim H$\/ consist of families 
$X=(X_A\in H(A))_{A\in\CC}$\/ together with isomorphisms $x_f\colon X_A \to H(f)X_B$\/ for 
$f\colon A\to B$ in $\CC$; these are subject to the unit and cocycle conditions $x_1=1$\/ and $x_{gf}=H(f)(x_g) \circ x_f$. 
\item [Morphisms] $m\colon (X,x)\to (Y,y)$\/ are 
families of morphisms $(m_A\colon X_A\to Y_A)_{A\in\CC}$\/ such that for all $f\colon A \to B$\/ we have $H(m_B) \circ x_f = y_f \circ m_A$. 
\item [Tensor] The monoidal structure and inverses of objects and of morphisms of $\lim H$\/ are computed pointwise, resulting in a $2$-group.
\end{description}

We can now state the main result of this section:

\begin{theorem}\label{thm:monoidalfib} For a small $2$-category $\CC$, the isotropy $2$-functor 
\[ \ZZ_{[\CC\op,\cat{MonCat}]}\colon [\CC\op,\cat{Moncat}]\to\cat{2Grp}\]
 is pseudonaturally equivalent to the functor 
sending $H\colon \CC\op\to\cat{MonCat}$ to $\Aut(\id[C])\times\lim (\mathrm{Pic}\circ H)$.
\end{theorem} 

  \begin{proof}  
  As the argument resembles the proof of Theorem~\ref{thm:picard}, we will restrict ourselves to a sketch. 
  The category $[\cat{C}\op,\cat{MonCat}]$\/ inherits binary coproducts from $\cat{MonCat}$, 
  so we may compute the isotropy $2$-group using Corollary~\ref{cor:2d-densitycontrolsisotropy} and the dense subcategory 
  $\mathcal{F}\colon \DD \to [\cat{C}\op,\cat{MonCat}]$\/ 
  provided by Lemma~\ref{lem:stuffisdenseinfunctorstomoncat}. Thus it suffices to provide equivalences 
  \[ \Aut(\id[C])\times\lim (\mathrm{Pic}\circ H)\to \LL_\mathcal{F}(H)\]
  that are pseudonatural in $H$.
  
  Given  an element $(\sigma,(X,x))$\/ of $\Aut(\id[C])\times\lim (\mathrm{Pic}\circ H)$, we define the corresponding 
  element 
   \[ (\LL_\mathcal{F}(H))_{\overline{A}}\colon H+\overline{A} \to H+\overline{A}\] 
  of $\LL_\mathcal{F}(H)$\/ at the object $\overline{A}$\/ as follows: given an object $C$\/ of $\CC$, 
  \[ (\LL_\mathcal{F}(H))_{\overline{A},C}\colon HC+\overline{A}C \to HC+\overline{A}C\] 
  is the identity on $HC$, and on $\overline{A}C$\/ sends a generator $l \in \CC(C,A)$\/ to 
 \[l \mapsto X_A\otimes (\sigma_A \circ l) \otimes X_A^{-1}.\]
  The action of $\LL_\mathcal{F}(H)$\/ at $\overline{A}+\overline{A}$ and $\overline{A}+\overline{A}+\overline{A}$\/ is
 defined similarly. Then the requisite pseudonaturality $2$-cells are induced by those of $\sigma$\/ and by the cups and caps as in the proof of 
  Theorem~\ref{thm:picard}. It is straightforward to check that this assignment on objects can be extend to a 
  ($\otimes$-reversing) monoidal functor that is fully faithful, and that this family of functors is pseudonatural in $H$.

  To see that each of these functors is essentially surjective on objects, 
  consider an arbitrary  pseudonatural autoequivalence of 
  \[ \DD\xrightarrow{\mathcal{F}}[\CC\op,\cat{MonCat}]\to H/[\CC\op,\cat{MonCat}].\] 
  The $1$-cells of such an autoequivalence consist of autoequivalences of $H+\overline{A}$, $H+\overline{A}+\overline{A}$ and of 
  $H+\overline{A}+\overline{A}+\overline{A}$ for each $A\in \CC$\/ that fix $H$\/ (within isomorphism). 
  By the universal property of $\overline{A}(A)=F(\CC(A,A))$\/ 
  and by Yoneda such autoequivalences correspond to elements of $H(A)+F\CC(A,A)$, which in turn correspond to words in $H(A)$\/ and $\CC(A,A)$. 
  As in the proof of Theorem~\ref{thm:picard}, the pseudonaturality constraints imply that these words have to be 
  (up to isomorphism) of the form $X_A\otimes \hat{f} \otimes X^{-1}_A$, with $\hat{f}\in F\CC(A,A)$. 
  For this to result in an equivalence, we must have $\hat{f}=f\in\CC(A,A)$ for some equivalence $f$. 
  Now, pseudonaturality in $A\in\CC$\/  implies that the family $X_A$ is in $\lim \mathrm{Pic}\circ H$, 
  and that the family $f\colon A\to A$\/ defines an autoequivalence of $\id[\CC]$. Thus we have demonstrated that 
  any element of $\LL_\mathcal{F}(H)$ is isomorphic to one whose $1$-cell components are exactly as determined 
  by an element of $\Aut(\id[C])\times\lim (\mathrm{Pic}\circ H)$. It remains to show that given such an element of 
  $\LL_\mathcal{F}(H)$, it is isomorphic to one whose $2$-cells are exactly as determined by an element of 
  $\Aut(\id[C])\times\lim (\mathrm{Pic}\circ H)$. This can be done by generalizing the similar argument in the 
  proof of Theorem~\ref{thm:picard} to show that, up to invertible modification, the $2$-cell components in $H$\/ are exactly as desired, 
  as the $2$-cell components in the other summand are as desired by Yoneda. 
\end{proof}

This result can be seen as a first step towards a two-dimensional generalization of~\cite{parker:presheaves}, 
where, under suitable conditions, one shows that there is a natural isomorphism between $\ZZ_{[\CC,T-Alg]}(H)$ and $\Aut(\id[\CC])\times \lim \ZZ_{T-Alg}\circ H$.

\bibliographystyle{plainurl}
\bibliography{2d-isotropy}

\end{document}

%% file: headeforpics.tex
\usepackage{tikz}
\usetikzlibrary{decorations.pathreplacing}
\usetikzlibrary{decorations.markings}
\usetikzlibrary{calc}
\usetikzlibrary{arrows}
\usetikzlibrary{arrows.meta}
\usetikzlibrary{shapes.geometric}
\tikzset{>={To[length=2.5pt,width=4pt]}}

\usetikzlibrary{intersections}

\newenvironment{pic}[1][]
{\begin{aligned}\begin{tikzpicture}[font=\tiny,#1]}
{\end{tikzpicture}\end{aligned}}

\pgfdeclarelayer{foreground}
\pgfdeclarelayer{background}
\pgfdeclarelayer{morphismlayer}
\pgfdeclarelayer{dashedmorphismlayer}
\pgfdeclarelayer{edgelayer}
\pgfdeclarelayer{nodelayer}
\pgfsetlayers{background,main,morphismlayer,dashedmorphismlayer,foreground,edgelayer,nodelayer}

\def\thickness{0.7pt}
\makeatletter
\pgfkeys{%
  /tikz/on layer/.code={
    \pgfonlayer{#1}\begingroup
    \aftergroup\endpgfonlayer
    \aftergroup\endgroup
  },
  /tikz/node on layer/.code={
    \gdef\node@@on@layer{%
      \setbox\tikz@tempbox=\hbox\bgroup\pgfonlayer{#1}\unhbox\tikz@tempbox\endpgfonlayer\pgfsetlinewidth{\thickness}\egroup}
    \aftergroup\node@on@layer
  },
  /tikz/end node on layer/.code={
    \endpgfonlayer\endgroup\endgroup
  }
}
\def\node@on@layer{\aftergroup\node@@on@layer}
\makeatother

\tikzstyle{braid}=[double=black, line width=3*\thickness, double distance=\thickness, draw=white, text=black]
\tikzset{every picture/.style={line width=\thickness, draw=black}}
\tikzstyle{pure}=[line width=.7pt]
\tikzstyle{string}=[line width=\thickness]
\tikzstyle{scalar}=[circle, inner sep=0pt, minimum width=15pt, draw, line width=\thickness, fill=white]
\tikzstyle{dot}=[circle, draw=black, fill=black!25, inner sep=.4ex, line width=\thickness, node on layer=foreground]
\tikzstyle{blackdot}=[circle, draw=black, fill=black!75, inner sep=.4ex, line width=\thickness, node on layer=foreground, text=white]
\tikzstyle{whitedot}=[circle, draw=black, fill=white, inner sep=.4ex, line width=\thickness, node on layer=foreground]
\tikzstyle{mixedmorphism}=[morphism, minimum width=30pt, minimum height=16pt, draw, font=\small, inner sep=0pt, fill=white, line width=\thickness,rounded corners=1ex]

\tikzstyle{triangle} = [regular polygon, regular polygon sides=3, draw=black, fill=black!20,scale=0.4, node on layer=foreground]

\tikzstyle{whitetriangle}=[triangle, fill=white]
\tikzstyle{greytriangle}=[triangle, fill=black!20]
\tikzstyle{darkgreytriangle}=[triangle, fill=black!50]
\tikzstyle{blacktriangle}=[triangle, fill=black]

\tikzstyle{invertedtriangle} = [triangle,scale=-1]
\tikzstyle{whiteinvertedtriangle}=[invertedtriangle, fill=white]
\tikzstyle{greyinvertedtriangle}=[invertedtriangle, fill=black!20]
\tikzstyle{darkgreyinvertedtriangle}=[invertedtriangle, fill=black!50]
\tikzstyle{blackinvertedtriangle}=[invertedtriangle, fill=black]

\tikzset{functor1/.style={gray!50}}
\tikzset{functor2/.style={gray!50}}
\tikzstyle{thick}=[line width=\thickness]
\tikzstyle{tiny}=[font=\tiny]
\tikzset{circlelabel/.style={draw, thick, circle, inner sep=-5pt,
 fill=white, minimum width=14pt, fill opacity=1, node on layer=foreground, font=\scriptsize}}
\tikzset{shade 1 transparent/.style={fill=black!20, fill opacity=0.8}}
\tikzset{shade 2 transparent/.style={fill=black!35, fill opacity=0.8}}
\tikzset{shade 3 transparent/.style={fill=black!50, fill opacity=0.8}}
\tikzset{shade 4 transparent/.style={fill=black!65, fill opacity=0.8}}
\tikzset{shade 1/.style={fill=black!16, fill opacity=1}}
\tikzset{shade 2/.style={fill=black!28, fill opacity=1}}
\tikzset{shade 3/.style={fill=black!40, fill opacity=1}}
\tikzset{shade 4/.style={fill=black!52, fill opacity=1}}

\def\strarr{line width=0.7pt, length=4pt, width=5pt, color=black}

\tikzset{arrow/.style={decoration={
    markings,
    mark=at position #1 with \arrow{>[\strarr]}},
    postaction=decorate},
    reverse arrow/.style={decoration={
    markings,
    mark=at position #1 with {{\arrow{<[\strarr]}}}},
    postaction=decorate}
}

\tikzset{overbrace/.style={
     decoration={brace},
     decorate}
}
\tikzset{underbrace/.style={
     decoration={brace, mirror},
     decorate}
}

\newif\ifblack\pgfkeys{/tikz/black/.is if=black}
\newif\ifwedge\pgfkeys{/tikz/wedge/.is if=wedge}
\newif\ifvflip\pgfkeys{/tikz/vflip/.is if=vflip}
\newif\ifhflip\pgfkeys{/tikz/hflip/.is if=hflip}
\newif\ifhvflip\pgfkeys{/tikz/hvflip/.is if=hvflip}
\newif\ifconnectsw\pgfkeys{/tikz/connect sw/.is if=connectsw}
\newif\ifconnectse\pgfkeys{/tikz/connect se/.is if=connectse}
\newif\ifconnectn\pgfkeys{/tikz/connect n/.is if=connectn}
\newif\ifconnects\pgfkeys{/tikz/connect s/.is if=connects}
\newif\ifconnectnw\pgfkeys{/tikz/connect nw/.is if=connectnw}
\newif\ifconnectne\pgfkeys{/tikz/connect ne/.is if=connectne}
\newif\ifconnectnwf\pgfkeys{/tikz/connect nw >/.is if=connectnwf}
\newif\ifconnectnef\pgfkeys{/tikz/connect ne >/.is if=connectnef}
\newif\ifconnectswf\pgfkeys{/tikz/connect sw >/.is if=connectswf}
\newif\ifconnectsef\pgfkeys{/tikz/connect se >/.is if=connectsef}
\newif\ifconnectnf\pgfkeys{/tikz/connect n >/.is if=connectnf}
\newif\ifconnectsf\pgfkeys{/tikz/connect s >/.is if=connectsf}
\newif\ifconnectnwr\pgfkeys{/tikz/connect nw </.is if=connectnwr}
\newif\ifconnectner\pgfkeys{/tikz/connect ne </.is if=connectner}
\newif\ifconnectswr\pgfkeys{/tikz/connect sw </.is if=connectswr}
\newif\ifconnectser\pgfkeys{/tikz/connect se </.is if=connectser}
\newif\ifconnectnr\pgfkeys{/tikz/connect n </.is if=connectnr}
\newif\ifconnectsr\pgfkeys{/tikz/connect s </.is if=connectsr}
\tikzset{keylengthnw/.initial=\connectheight}
\tikzset{keylengthn/.initial =\connectheight}
\tikzset{keylengthne/.initial=\connectheight}
\tikzset{keylengthsw/.initial=\connectheight}
\tikzset{keylengths/.initial =\connectheight}
\tikzset{keylengthse/.initial=\connectheight}
\tikzset{connect nw length/.style={connect nw=true, keylengthnw={#1}}}
\tikzset{connect n length/.style ={connect n =true, keylengthn ={#1}}}
\tikzset{connect ne length/.style={connect ne=true, keylengthne={#1}}}
\tikzset{connect sw length/.style={connect sw=true, keylengthsw={#1}}}
\tikzset{connect s length/.style ={connect s =true, keylengths ={#1}}}
\tikzset{connect se length/.style={connect se=true, keylengthse={#1}}}
\tikzset{connect nw < length/.style={connect nw <=true, keylengthnw={#1}}}
\tikzset{connect n < length/.style ={connect n <=true,  keylengthn ={#1}}}
\tikzset{connect ne < length/.style={connect ne <=true, keylengthne={#1}}}
\tikzset{connect sw < length/.style={connect sw <=true, keylengthnw={#1}}}
\tikzset{connect s < length/.style ={connect s <=true,  keylengths ={#1}}}
\tikzset{connect se < length/.style={connect se <=true, keylengthse={#1}}}
\tikzset{connect nw > length/.style={connect nw >=true, keylengthnw={#1}}}
\tikzset{connect n > length/.style ={connect n >=true,  keylengthn ={#1}}}
\tikzset{connect ne > length/.style={connect ne >=true, keylengthne={#1}}}
\tikzset{connect sw > length/.style={connect sw >=true, keylengthsw={#1}}}
\tikzset{connect s > length/.style ={connect s >=true,  keylengths ={#1}}}
\tikzset{connect se > length/.style={connect se >=true, keylengthse={#1}}}

\newlength\morphismheight
\newlength\wedgewidth
\newlength\minimummorphismwidth
\newlength\stateheight
\newlength\minimumstatewidth
\newlength\connectheight
\tikzset{width/.initial=\minimummorphismwidth}
\tikzset{colour/.initial=white}

\makeatletter
\pgfarrowsdeclare{thickarrow}{thickarrow}
{
  \pgfutil@tempdima=-0.84pt%
  \advance\pgfutil@tempdima by-1.3\pgflinewidth%
  \pgfutil@tempdimb=-1.7pt%
  \advance\pgfutil@tempdimb by.625\pgflinewidth%
  \pgfarrowsleftextend{+\pgfutil@tempdima}
  \pgfarrowsrightextend{+\pgfutil@tempdimb}
}
{
  \pgfmathparse{\pgfgetarrowoptions{thickarrow}}%
  \pgfsetlinewidth{1.25 pt}
  \pgfutil@tempdima=0.28pt%
  \advance\pgfutil@tempdima by.3\pgflinewidth%
  \pgfsetlinewidth{0.8\pgflinewidth}
  \pgfsetdash{}{+0pt}
  \pgfsetroundcap
  \pgfsetroundjoin
  \pgfpathmoveto{\pgfqpoint{-3\pgfutil@tempdima}{4\pgfutil@tempdima}}
  \pgfpathcurveto
  {\pgfqpoint{-2.75\pgfutil@tempdima}{2.5\pgfutil@tempdima}}
  {\pgfqpoint{0pt}{0.25\pgfutil@tempdima}}
  {\pgfqpoint{0.75\pgfutil@tempdima}{0pt}}
  \pgfpathcurveto
  {\pgfqpoint{0pt}{-0.25\pgfutil@tempdima}}
  {\pgfqpoint{-2.75\pgfutil@tempdima}{-2.5\pgfutil@tempdima}}
  {\pgfqpoint{-3\pgfutil@tempdima}{-4\pgfutil@tempdima}}
  \pgfusepathqstroke
}
\pgfarrowsdeclare{reversethickarrow}{reversethickarrow}
{
  \pgfutil@tempdima=-0.84pt%
  \advance\pgfutil@tempdima by-1.3\pgflinewidth%
  \pgfutil@tempdimb=0.2pt%
  \advance\pgfutil@tempdimb by.625\pgflinewidth%
  \pgfarrowsleftextend{+\pgfutil@tempdima}
  \pgfarrowsrightextend{+\pgfutil@tempdimb}
}
{
  \pgftransformxscale{-1}
  \pgfmathparse{\pgfgetarrowoptions{thickarrow}}%
  \ifpgfmathunitsdeclared%
    \pgfmathparse{\pgfmathresult pt}%
  \else%
    \pgfmathparse{\pgfmathresult*\pgflinewidth}%
  \fi%
  \let\thickness=\pgfmathresult
  \pgfsetlinewidth{1.25 pt}
  \pgfutil@tempdima=0.28pt%
  \advance\pgfutil@tempdima by.3\pgflinewidth%
  \pgfsetlinewidth{0.8\pgflinewidth}
  \pgfsetdash{}{+0pt}
  \pgfsetroundcap
  \pgfsetroundjoin
  \pgfpathmoveto{\pgfqpoint{-3\pgfutil@tempdima}{4\pgfutil@tempdima}}
  \pgfpathcurveto
  {\pgfqpoint{-2.75\pgfutil@tempdima}{2.5\pgfutil@tempdima}}
  {\pgfqpoint{0pt}{0.25\pgfutil@tempdima}}
  {\pgfqpoint{0.75\pgfutil@tempdima}{0pt}}
  \pgfpathcurveto
  {\pgfqpoint{0pt}{-0.25\pgfutil@tempdima}}
  {\pgfqpoint{-2.75\pgfutil@tempdima}{-2.5\pgfutil@tempdima}}
  {\pgfqpoint{-3\pgfutil@tempdima}{-4\pgfutil@tempdima}}
  \pgfusepathqstroke
}
\makeatother

\tikzset{diredge/.style={decoration={
  markings,
  mark=at position 0.525 with {\arrow{#1}}},postaction={decorate}}}
\tikzset{
    diredge/.default=>
}

\tikzset{diredgestart/.style={decoration={
  markings,
  mark=at position 4pt with {\arrow{#1}}},postaction={decorate}}}
\tikzset{
    diredgestart/.default=<
}

\tikzset{diredgeend/.style={decoration={
  markings,
  mark=at position 1 with {\arrow{#1}}},postaction={decorate}}}
\tikzset{
    diredgeend/.default=>
}

\makeatletter
\pgfdeclareshape{ground}
{
    \savedanchor\centerpoint
    {
        \pgf@x=0pt
        \pgf@y=0pt
    }
    \anchor{center}{\centerpoint}
    \anchorborder{\centerpoint}
    \anchor{north}
    {
        \pgf@x=0pt
        \pgf@y=0.5*0.33*\stateheight
    }
    \anchor{south}
    {
        \pgf@x=0pt
        \pgf@y=0pt
    }
    \saveddimen\overallwidth
    {
        \pgfkeysgetvalue{/pgf/minimum width}{\minwidth}
        \pgf@x=\minimumstatewidth
        \ifdim\pgf@x<\minwidth
            \pgf@x=\minwidth
        \fi
    }
    \backgroundpath
    {
        \begin{pgfonlayer}{foreground}
        \pgfsetstrokecolor{black}
        \pgfsetlinewidth{1.25pt}
        \ifhflip
            \pgftransformyscale{-1}
        \fi
        \pgftransformscale{0.5}
        \pgfpathmoveto{\pgfpoint{-0.5*\overallwidth}{0}}
        \pgfpathlineto{\pgfpoint{0.5*\overallwidth}{0}}
        \pgfpathmoveto{\pgfpoint{-0.33*\overallwidth}{0.33*\stateheight}}
        \pgfpathlineto{\pgfpoint{0.33*\overallwidth}{0.33*\stateheight}}
        \pgfpathmoveto{\pgfpoint{-0.16*\overallwidth}{0.66*\stateheight}}
        \pgfpathlineto{\pgfpoint{0.16*\overallwidth}{0.66*\stateheight}}
        \pgfpathmoveto{\pgfpoint{-0.02*\overallwidth}{\stateheight}}
        \pgfpathlineto{\pgfpoint{0.02*\overallwidth}{\stateheight}}
        \pgfusepath{stroke}
        \end{pgfonlayer}
    }
}
\tikzset{forward arrow style/.style={every to/.style, decoration={
    markings,
    mark=at position 0.5*\pgfdecoratedpathlength+2pt with \arrow{>[\strarr]}},
    postaction=decorate}}
\tikzset{reverse arrow style/.style={every to/.style, decoration={
    markings,
    mark=at position 0.5*\pgfdecoratedpathlength+2pt with \arrow{<[\strarr]}},
    postaction=decorate}}
\pgfdeclareshape{morphismshape}
{
    \savedanchor\centerpoint
    {
        \pgf@x=0pt
        \pgf@y=0pt
    }
    \anchor{center}{\centerpoint}
    \anchorborder{\centerpoint}
    \saveddimen\savedlengthnw
    {
        \pgfkeysgetvalue{/tikz/keylengthnw}{\len}
        \pgf@x=\len
    }
    \saveddimen\savedlengthn
    {
        \pgfkeysgetvalue{/tikz/keylengthn}{\len}
        \pgf@x=\len
    }
    \saveddimen\savedlengthne
    {
        \pgfkeysgetvalue{/tikz/keylengthne}{\len}
        \pgf@x=\len
    }
    \saveddimen\savedlengthsw
    {
        \pgfkeysgetvalue{/tikz/keylengthsw}{\len}
        \pgf@x=\len
    }
    \saveddimen\savedlengths
    {
        \pgfkeysgetvalue{/tikz/keylengths}{\len}
        \pgf@x=\len
    }
    \saveddimen\savedlengthse
    {
        \pgfkeysgetvalue{/tikz/keylengthse}{\len}
        \pgf@x=\len
    }
    \saveddimen\overallwidth
    {
        \pgfkeysgetvalue{/tikz/width}{\minwidth}
        \pgf@x=\wd\pgfnodeparttextbox
        \ifdim\pgf@x<\minwidth
            \pgf@x=\minwidth
        \fi
    }
    \savedanchor{\upperrightcorner}
    {
        \pgf@y=.5\ht\pgfnodeparttextbox
        \advance\pgf@y by -.5\dp\pgfnodeparttextbox
        \pgf@x=.5\wd\pgfnodeparttextbox
    }
    \anchor{north}
    {
        \pgf@x=0pt
        \pgf@y=0.5\morphismheight
    }
    \anchor{north east}
    {
        \pgf@x=\overallwidth
        \multiply \pgf@x by 2
        \divide \pgf@x by 5
        \pgf@y=0.5\morphismheight
    }
    \anchor{east}
    {
        \pgf@x=\overallwidth
        \divide \pgf@x by 2
        \advance \pgf@x by 5pt
        \pgf@y=0pt
    }
    \anchor{west}
    {
        \pgf@x=-\overallwidth
        \divide \pgf@x by 2
        \advance \pgf@x by -5pt
        \pgf@y=0pt
    }
    \anchor{north west}
    {
        \pgf@x=-\overallwidth
        \multiply \pgf@x by 2
        \divide \pgf@x by 5
        \pgf@y=0.5\morphismheight
    }
    \anchor{connect nw}
    {
        \pgf@x=-\overallwidth
        \multiply \pgf@x by 2
        \divide \pgf@x by 5
        \pgf@y=0.5\morphismheight
        \advance\pgf@y by \savedlengthnw
    }
    \anchor{connect ne}
    {
        \pgf@x=\overallwidth
        \multiply \pgf@x by 2
        \divide \pgf@x by 5
        \pgf@y=0.5\morphismheight
        \advance\pgf@y by \savedlengthne
    }
    \anchor{connect sw}
    {
        \pgf@x=-\overallwidth
        \multiply \pgf@x by 2
        \divide \pgf@x by 5
        \pgf@y=-0.5\morphismheight
        \advance\pgf@y by -\savedlengthsw
    }
    \anchor{connect se}
    {
        \pgf@x=\overallwidth
        \multiply \pgf@x by 2
        \divide \pgf@x by 5
        \pgf@y=-0.5\morphismheight
        \advance\pgf@y by -\savedlengthse
    }
    \anchor{connect n}
    {
        \pgf@x=0pt
        \pgf@y=0.5\morphismheight
        \advance\pgf@y by \savedlengthn
    }
    \anchor{connect s}
    {
        \pgf@x=0pt
        \pgf@y=-0.5\morphismheight
        \advance\pgf@y by -\savedlengths
    }
    \anchor{south east}
    {
        \pgf@x=\overallwidth
        \multiply \pgf@x by 2
        \divide \pgf@x by 5
        \pgf@y=-0.5\morphismheight
    }
    \anchor{south west}
    {
        \pgf@x=-\overallwidth
        \multiply \pgf@x by 2
        \divide \pgf@x by 5
        \pgf@y=-0.5\morphismheight
    }
    \anchor{south}
    {
        \pgf@x=0pt
        \pgf@y=-0.5\morphismheight
    }
    \anchor{text}
    {
        \upperrightcorner
        \pgf@x=-\pgf@x
        \pgf@y=-\pgf@y
    }
    \backgroundpath
    {
    \begin{scope}
        \pgfkeysgetvalue{/tikz/fill}{\morphismfill}
        \pgfsetstrokecolor{black}
        \pgfsetlinewidth{\thickness}
        \begin{scope}
            \begin{pgfonlayer}{morphismlayer}
        \pgfsetstrokecolor{black}
        \pgfsetfillcolor{\pgfkeysvalueof{/tikz/colour}}
        \pgfsetlinewidth{\thickness}
                \ifhflip
                    \pgftransformyscale{-1}
                \fi
                \ifvflip
                    \pgftransformxscale{-1}
                \fi
                \ifhvflip
                    \pgftransformxscale{-1}
                    \pgftransformyscale{-1}
                \fi
                \pgfpathmoveto{\pgfpoint
                    {-0.5*\overallwidth-5pt}
                    {0.5*\morphismheight}}
                \pgfpathlineto{\pgfpoint
                    {0.5*\overallwidth+5pt}
                    {0.5*\morphismheight}}
                \ifwedge
                    \pgfpathlineto{\pgfpoint
                        {0.5*\overallwidth + \wedgewidth}
                        {-0.5*\morphismheight}}
                \else
                    \pgfpathlineto{\pgfpoint
                        {0.5*\overallwidth + 5pt}
                        {-0.5*\morphismheight}}
                \fi
                \pgfpathlineto{\pgfpoint
                    {-0.5*\overallwidth-5pt}
                    {-0.5*\morphismheight}}
                \pgfpathclose
                \pgfusepath{fill,stroke}
            \end{pgfonlayer}
        \end{scope}
        \ifconnectnw
            \pgfpathmoveto{\pgfpoint
                {-0.4*\overallwidth}
                {0.5*\morphismheight}}
            \pgfpathlineto{\pgfpoint
                {-0.4*\overallwidth}
                {0.5*\morphismheight+\savedlengthnw}}
            \pgfusepath{stroke}
        \fi
        \ifconnectne
            \pgfpathmoveto{\pgfpoint
                {0.4*\overallwidth}
                {0.5*\morphismheight}}
            \pgfpathlineto{\pgfpoint
                {0.4*\overallwidth}
                {0.5*\morphismheight+\savedlengthne}}
            \pgfusepath{stroke}
        \fi
        \ifconnectsw
            \pgfpathmoveto{\pgfpoint
                {-0.4*\overallwidth}
                {-0.5*\morphismheight}}
            \pgfpathlineto{\pgfpoint
                {-0.4*\overallwidth}
                {-0.5*\morphismheight-\savedlengthsw}}
            \pgfusepath{stroke}
        \fi
        \ifconnectse
            \pgfpathmoveto{\pgfpoint
                {0.4*\overallwidth}
                {-0.5*\morphismheight}}
            \pgfpathlineto{\pgfpoint
                {0.4*\overallwidth}
                {-0.5*\morphismheight-\savedlengthse}}
            \pgfusepath{stroke}
        \fi
        \ifconnectn
            \pgfpathmoveto{\pgfpoint
                {0pt}
                {0.5*\morphismheight}}
            \pgfpathlineto{\pgfpoint
                {0pt}
                {0.5*\morphismheight+\savedlengthn}}
            \pgfusepath{stroke}
        \fi
        \ifconnects
            \pgfpathmoveto{\pgfpoint
                {0pt}
                {-0.5*\morphismheight}}
            \pgfpathlineto{\pgfpoint
                {0pt}
                {-0.5*\morphismheight-\savedlengths}}
            \pgfusepath{stroke}
        \fi
        \ifconnectnwf
            \draw [forward arrow style] (-0.4*\overallwidth,0.5*\morphismheight)
                to (-0.4*\overallwidth,0.5*\morphismheight+\savedlengthnw);
        \fi
        \ifconnectnef
            \draw [forward arrow style] (0.4*\overallwidth,0.5*\morphismheight)
                to (0.4*\overallwidth,0.5*\morphismheight+\savedlengthne);
        \fi
        \ifconnectswf
            \draw [forward arrow style] (-0.4*\overallwidth,-0.5*\morphismheight-\savedlengthsw)
                to (-0.4*\overallwidth,-0.5*\morphismheight);
        \fi
        \ifconnectsef
            \draw [forward arrow style] (0.4*\overallwidth,-0.5*\morphismheight-\savedlengthse)
                to (0.4*\overallwidth,-0.5*\morphismheight);
        \fi
        \ifconnectnf
            \draw [forward arrow style] (0,0.5*\morphismheight)
                to (0,0.5*\morphismheight+\savedlengthn);
        \fi
        \ifconnectsf
            \draw [forward arrow style] (0,-0.5*\morphismheight-\savedlengths)
                to (0,-0.5*\morphismheight);
        \fi
        \ifconnectnwr
            \draw [reverse arrow style] (-0.4*\overallwidth,0.5*\morphismheight)
                to (-0.4*\overallwidth,0.5*\morphismheight+\savedlengthnw);
        \fi
        \ifconnectner
            \draw [reverse arrow style] (0.4*\overallwidth,0.5*\morphismheight)
                to (0.4*\overallwidth,0.5*\morphismheight+\savedlengthne);
        \fi
        \ifconnectswr
            \draw [reverse arrow style] (-0.4*\overallwidth,-0.5*\morphismheight-\savedlengthsw)
                to (-0.4*\overallwidth,-0.5*\morphismheight);
        \fi
        \ifconnectser
            \draw [reverse arrow style] (0.4*\overallwidth,-0.5*\morphismheight-\savedlengthse)
                to (0.4*\overallwidth,-0.5*\morphismheight);
        \fi
        \ifconnectnr
            \draw [reverse arrow style] (0,0.5*\morphismheight)
                to (0,0.5*\morphismheight+\savedlengthn);
        \fi
        \ifconnectsr
            \draw [reverse arrow style] (0,-0.5*\morphismheight-\savedlengths)
                to (0,-0.5*\morphismheight);
        \fi
    \end{scope}
    }
}
\tikzset{morphism/.style={morphismshape, node on layer=foreground}}
\pgfdeclareshape{dashedmorphismshape}
{
    \savedanchor\centerpoint
    {
        \pgf@x=0pt
        \pgf@y=0pt
    }
    \anchor{center}{\centerpoint}
    \anchorborder{\centerpoint}
    \saveddimen\savedlengthnw
    {
        \pgfkeysgetvalue{/tikz/keylengthnw}{\len}
        \pgf@x=\len
    }
    \saveddimen\savedlengthn
    {
        \pgfkeysgetvalue{/tikz/keylengthn}{\len}
        \pgf@x=\len
    }
    \saveddimen\savedlengthne
    {
        \pgfkeysgetvalue{/tikz/keylengthne}{\len}
        \pgf@x=\len
    }
    \saveddimen\savedlengthsw
    {
        \pgfkeysgetvalue{/tikz/keylengthsw}{\len}
        \pgf@x=\len
    }
    \saveddimen\savedlengths
    {
        \pgfkeysgetvalue{/tikz/keylengths}{\len}
        \pgf@x=\len
    }
    \saveddimen\savedlengthse
    {
        \pgfkeysgetvalue{/tikz/keylengthse}{\len}
        \pgf@x=\len
    }
    \saveddimen\overallwidth
    {
        \pgfkeysgetvalue{/tikz/width}{\minwidth}
        \pgf@x=\wd\pgfnodeparttextbox
        \ifdim\pgf@x<\minwidth
            \pgf@x=\minwidth
        \fi
    }
    \savedanchor{\upperrightcorner}
    {
        \pgf@y=.5\ht\pgfnodeparttextbox
        \advance\pgf@y by -.5\dp\pgfnodeparttextbox
        \pgf@x=.5\wd\pgfnodeparttextbox
    }
    \anchor{north}
    {
        \pgf@x=0pt
        \pgf@y=0.5\morphismheight
    }
    \anchor{north east}
    {
        \pgf@x=\overallwidth
        \multiply \pgf@x by 2
        \divide \pgf@x by 5
        \pgf@y=0.5\morphismheight
    }
    \anchor{east}
    {
        \pgf@x=\overallwidth
        \divide \pgf@x by 2
        \advance \pgf@x by 5pt
        \pgf@y=0pt
    }
    \anchor{west}
    {
        \pgf@x=-\overallwidth
        \divide \pgf@x by 2
        \advance \pgf@x by -5pt
        \pgf@y=0pt
    }
    \anchor{north west}
    {
        \pgf@x=-\overallwidth
        \multiply \pgf@x by 2
        \divide \pgf@x by 5
        \pgf@y=0.5\morphismheight
    }
    \anchor{connect nw}
    {
        \pgf@x=-\overallwidth
        \multiply \pgf@x by 2
        \divide \pgf@x by 5
        \pgf@y=0.5\morphismheight
        \advance\pgf@y by \savedlengthnw
    }
    \anchor{connect ne}
    {
        \pgf@x=\overallwidth
        \multiply \pgf@x by 2
        \divide \pgf@x by 5
        \pgf@y=0.5\morphismheight
        \advance\pgf@y by \savedlengthne
    }
    \anchor{connect sw}
    {
        \pgf@x=-\overallwidth
        \multiply \pgf@x by 2
        \divide \pgf@x by 5
        \pgf@y=-0.5\morphismheight
        \advance\pgf@y by -\savedlengthsw
    }
    \anchor{connect se}
    {
        \pgf@x=\overallwidth
        \multiply \pgf@x by 2
        \divide \pgf@x by 5
        \pgf@y=-0.5\morphismheight
        \advance\pgf@y by -\savedlengthse
    }
    \anchor{connect n}
    {
        \pgf@x=0pt
        \pgf@y=0.5\morphismheight
        \advance\pgf@y by \savedlengthn
    }
    \anchor{connect s}
    {
        \pgf@x=0pt
        \pgf@y=-0.5\morphismheight
        \advance\pgf@y by -\savedlengths
    }
    \anchor{south east}
    {
        \pgf@x=\overallwidth
        \multiply \pgf@x by 2
        \divide \pgf@x by 5
        \pgf@y=-0.5\morphismheight
    }
    \anchor{south west}
    {
        \pgf@x=-\overallwidth
        \multiply \pgf@x by 2
        \divide \pgf@x by 5
        \pgf@y=-0.5\morphismheight
    }
    \anchor{south}
    {
        \pgf@x=0pt
        \pgf@y=-0.5\morphismheight
    }
    \anchor{text}
    {
        \upperrightcorner
        \pgf@x=-\pgf@x
        \pgf@y=-\pgf@y
    }
    \backgroundpath
    {
    \begin{scope}
        \pgfkeysgetvalue{/tikz/fill}{\morphismfill}
        \pgfsetstrokecolor{black}
        \pgfsetlinewidth{\thickness}
        \begin{scope}
            \begin{pgfonlayer}{dashedmorphismlayer}
        \pgfsetstrokecolor{black}
        \pgfsetdash{{3pt}{3pt}}{0pt}
        \pgfsetfillcolor{\pgfkeysvalueof{/tikz/colour}}
        \pgfsetlinewidth{\thickness}
                \ifhflip
                    \pgftransformyscale{-1}
                \fi
                \ifvflip
                    \pgftransformxscale{-1}
                \fi
                \ifhvflip
                    \pgftransformxscale{-1}
                    \pgftransformyscale{-1}
                \fi
                \pgfpathmoveto{\pgfpoint
                    {-0.5*\overallwidth-5pt}
                    {0.5*\morphismheight}}
                \pgfpathlineto{\pgfpoint
                    {0.5*\overallwidth+5pt}
                    {0.5*\morphismheight}}
                \ifwedge
                    \pgfpathlineto{\pgfpoint
                        {0.5*\overallwidth + \wedgewidth}
                        {-0.5*\morphismheight}}
                \else
                    \pgfpathlineto{\pgfpoint
                        {0.5*\overallwidth + 5pt}
                        {-0.5*\morphismheight}}
                \fi
                \pgfpathlineto{\pgfpoint
                    {-0.5*\overallwidth-5pt}
                    {-0.5*\morphismheight}}
                \pgfpathclose
                \pgfusepath{fill,stroke}
            \end{pgfonlayer}
        \end{scope}
        \ifconnectnw
            \pgfpathmoveto{\pgfpoint
                {-0.4*\overallwidth}
                {0.5*\morphismheight}}
            \pgfpathlineto{\pgfpoint
                {-0.4*\overallwidth}
                {0.5*\morphismheight+\savedlengthnw}}
            \pgfusepath{stroke}
        \fi
        \ifconnectne
            \pgfpathmoveto{\pgfpoint
                {0.4*\overallwidth}
                {0.5*\morphismheight}}
            \pgfpathlineto{\pgfpoint
                {0.4*\overallwidth}
                {0.5*\morphismheight+\savedlengthne}}
            \pgfusepath{stroke}
        \fi
        \ifconnectsw
            \pgfpathmoveto{\pgfpoint
                {-0.4*\overallwidth}
                {-0.5*\morphismheight}}
            \pgfpathlineto{\pgfpoint
                {-0.4*\overallwidth}
                {-0.5*\morphismheight-\savedlengthsw}}
            \pgfusepath{stroke}
        \fi
        \ifconnectse
            \pgfpathmoveto{\pgfpoint
                {0.4*\overallwidth}
                {-0.5*\morphismheight}}
            \pgfpathlineto{\pgfpoint
                {0.4*\overallwidth}
                {-0.5*\morphismheight-\savedlengthse}}
            \pgfusepath{stroke}
        \fi
        \ifconnectn
            \pgfpathmoveto{\pgfpoint
                {0pt}
                {0.5*\morphismheight}}
            \pgfpathlineto{\pgfpoint
                {0pt}
                {0.5*\morphismheight+\savedlengthn}}
            \pgfusepath{stroke}
        \fi
        \ifconnects
            \pgfpathmoveto{\pgfpoint
                {0pt}
                {-0.5*\morphismheight}}
            \pgfpathlineto{\pgfpoint
                {0pt}
                {-0.5*\morphismheight-\savedlengths}}
            \pgfusepath{stroke}
        \fi
        \ifconnectnwf
            \draw [forward arrow style] (-0.4*\overallwidth,0.5*\morphismheight)
                to (-0.4*\overallwidth,0.5*\morphismheight+\savedlengthnw);
        \fi
        \ifconnectnef
            \draw [forward arrow style] (0.4*\overallwidth,0.5*\morphismheight)
                to (0.4*\overallwidth,0.5*\morphismheight+\savedlengthne);
        \fi
        \ifconnectswf
            \draw [forward arrow style] (-0.4*\overallwidth,-0.5*\morphismheight-\savedlengthsw)
                to (-0.4*\overallwidth,-0.5*\morphismheight);
        \fi
        \ifconnectsef
            \draw [forward arrow style] (0.4*\overallwidth,-0.5*\morphismheight-\savedlengthse)
                to (0.4*\overallwidth,-0.5*\morphismheight);
        \fi
        \ifconnectnf
            \draw [forward arrow style] (0,0.5*\morphismheight)
                to (0,0.5*\morphismheight+\savedlengthn);
        \fi
        \ifconnectsf
            \draw [forward arrow style] (0,-0.5*\morphismheight-\savedlengths)
                to (0,-0.5*\morphismheight);
        \fi
        \ifconnectnwr
            \draw [reverse arrow style] (-0.4*\overallwidth,0.5*\morphismheight)
                to (-0.4*\overallwidth,0.5*\morphismheight+\savedlengthnw);
        \fi
        \ifconnectner
            \draw [reverse arrow style] (0.4*\overallwidth,0.5*\morphismheight)
                to (0.4*\overallwidth,0.5*\morphismheight+\savedlengthne);
        \fi
        \ifconnectswr
            \draw [reverse arrow style] (-0.4*\overallwidth,-0.5*\morphismheight-\savedlengthsw)
                to (-0.4*\overallwidth,-0.5*\morphismheight);
        \fi
        \ifconnectser
            \draw [reverse arrow style] (0.4*\overallwidth,-0.5*\morphismheight-\savedlengthse)
                to (0.4*\overallwidth,-0.5*\morphismheight);
        \fi
        \ifconnectnr
            \draw [reverse arrow style] (0,0.5*\morphismheight)
                to (0,0.5*\morphismheight+\savedlengthn);
        \fi
        \ifconnectsr
            \draw [reverse arrow style] (0,-0.5*\morphismheight-\savedlengths)
                to (0,-0.5*\morphismheight);
        \fi
    \end{scope}
    }
}
\tikzset{dashedmorphism/.style={dashedmorphismshape, node on layer=foreground}}
\pgfdeclareshape{swish right}
{
    \savedanchor\centerpoint
    {
        \pgf@x=0pt
        \pgf@y=0pt
    }
    \anchor{center}{\centerpoint}
    \anchorborder{\centerpoint}
    \anchor{north}
    {
        \pgf@x=\minimummorphismwidth
        \divide\pgf@x by 5
        \pgf@y=\morphismheight
        \divide\pgf@y by 2
        \advance\pgf@y by \connectheight
    }
    \anchor{south}
    {
        \pgf@x=-\minimummorphismwidth
        \divide\pgf@x by 5
        \pgf@y=-\morphismheight
        \divide\pgf@y by 2
        \advance\pgf@y by -\connectheight
    }
    \backgroundpath
    {
        \pgfsetstrokecolor{black}
        \pgfsetlinewidth{\thickness}
        \pgfpathmoveto{\pgfpoint
            {-0.2*\minimummorphismwidth}
            {-0.5*\morphismheight-\connectheight}}
        \pgfpathcurveto
            {\pgfpoint{-0.2*\minimummorphismwidth}{0pt}}
            {\pgfpoint{0.2*\minimummorphismwidth}{0pt}}
            {\pgfpoint
                {0.2*\minimummorphismwidth}
                {0.5*\morphismheight+\connectheight}}
        \pgfusepath{stroke}
    }
}
\pgfdeclareshape{swish left}
{
    \savedanchor\centerpoint
    {
        \pgf@x=0pt
        \pgf@y=0pt
    }
    \anchor{center}{\centerpoint}
    \anchorborder{\centerpoint}
    \anchor{north}
    {
        \pgf@x=-\minimummorphismwidth
        \divide\pgf@x by 5
        \pgf@y=\morphismheight
        \divide\pgf@y by 2
        \advance\pgf@y by \connectheight
    }
    \anchor{south}
    {
        \pgf@x=\minimummorphismwidth
        \divide\pgf@x by 5
        \pgf@y=-\morphismheight
        \divide\pgf@y by 2
        \advance\pgf@y by -\connectheight
    }
    \backgroundpath
    {
        \pgfsetstrokecolor{black}
        \pgfsetlinewidth{\thickness}
        \pgfpathmoveto{\pgfpoint
            {0.2*\minimummorphismwidth}
            {-0.5*\morphismheight-\connectheight}}
        \pgfpathcurveto
            {\pgfpoint{0.2*\minimummorphismwidth}{0pt}}
            {\pgfpoint{-0.2*\minimummorphismwidth}{0pt}}
            {\pgfpoint
                {-0.2*\minimummorphismwidth}
                {0.5*\morphismheight+\connectheight}}
        \pgfusepath{stroke}
    }
}
\pgfdeclareshape{state}
{
    \savedanchor\centerpoint
    {
        \pgf@x=0pt
        \pgf@y=0pt
    }
    \anchor{center}{\centerpoint}
    \anchorborder{\centerpoint}
    \saveddimen\overallwidth
    {
        \pgf@x=3\wd\pgfnodeparttextbox
        \ifdim\pgf@x<\minimumstatewidth
            \pgf@x=\minimumstatewidth
        \fi
    }
    \savedanchor{\upperrightcorner}
    {
        \pgf@x=.5\wd\pgfnodeparttextbox
        \pgf@y=.5\ht\pgfnodeparttextbox
        \advance\pgf@y by -.5\dp\pgfnodeparttextbox
    }
    \anchor{north}
    {
        \pgf@x=0pt
        \pgf@y=0pt
    }
    \anchor{south}
    {
        \pgf@x=0pt
        \pgf@y=0pt
    }
    \anchor{A}
    {
        \pgf@x=-\overallwidth
        \divide\pgf@x by 4
        \pgf@y=0pt
    }
    \anchor{B}
    {
        \pgf@x=\overallwidth
        \divide\pgf@x by 4
        \pgf@y=0pt
    }
    \anchor{text}
    {
        \upperrightcorner
        \pgf@x=-\pgf@x
        \ifhflip
            \pgf@y=-\pgf@y
            \advance\pgf@y by 0.4\stateheight
        \else
            \pgf@y=-\pgf@y
            \advance\pgf@y by -0.4\stateheight
        \fi
    }
    \backgroundpath
    {
       \begin{pgfonlayer}{foreground}
        \pgfsetstrokecolor{black}
        \pgfsetlinewidth{\thickness}
        \pgfpathmoveto{\pgfpoint{-0.5*\overallwidth}{0}}
        \pgfpathlineto{\pgfpoint{0.5*\overallwidth}{0}}
        \ifhflip
            \pgfpathlineto{\pgfpoint{0}{\stateheight}}
        \else
            \pgfpathlineto{\pgfpoint{0}{-\stateheight}}
        \fi
        \pgfpathclose
        \ifblack
            \pgfsetfillcolor{black!50}
            \pgfusepath{fill,stroke}
        \else
            \pgfusepath{stroke}
        \fi
       \end{pgfonlayer}
    }
}
\makeatother

\newcommand{\tinyhandle}[1][dot]{\ensuremath{\smash{\raisebox{-1pt}{\ensuremath{\hspace{-2pt}\begin{pic}[scale=0.33]
        \node (0) at (0,0) {};
        \node[dot, inner sep=1.0pt] (1) at (0,0.3) {};
        \node[dot, inner sep=1.0pt] (2) at (0,0.7) {};
        \node (3) at (0,1) {};
        \draw (0.center) to (1.center);
        \draw (2.center) to (3.center);
        \draw[in=180, out=180, looseness=2] (1.center) to (2.center);
        \draw[in=0, out=0, looseness=2] (1.center) to (2.center);
\end{pic}\hspace{-1pt}}}}}}

\makeatletter
\pgfdeclareshape{arrow shape}
{
    \savedanchor\centerpoint
    {
        \pgf@x=0pt
        \pgf@y=0pt
    }
    \anchor{center}{\centerpoint}
    \anchor{south}{\centerpoint}
    \anchor{north}{\centerpoint}
    \anchorborder{\centerpoint}
    \backgroundpath
    {
        \begin{pgfonlayer}{foreground}
        \pgfsetstrokecolor{black}
        \pgfsetarrowsstart{>[\strarr]}
        \pgfsetlinewidth{\thickness}
        \pgfpathmoveto{\pgfpoint{0}{-2pt}}
        \pgfpathlineto{\pgfpoint{0}{1pt}}
        \pgfusepath{stroke}
        \end{pgfonlayer}
    }
}
\pgfdeclareshape{double arrow shape}
{
    \savedanchor\centerpoint
    {
        \pgf@x=0pt
        \pgf@y=0pt
    }
    \anchor{center}{\centerpoint}
    \anchor{south}{\centerpoint}
    \anchor{north}{\centerpoint}
    \anchorborder{\centerpoint}
    \backgroundpath
    {
        \begin{pgfonlayer}{foreground}
        \pgfsetstrokecolor{black}
        \pgfsetarrowsstart{>[\strarr]}
        \pgfsetlinewidth{\thickness}
        \pgfpathmoveto{\pgfpoint{0}{-3.5pt}}
        \pgfpathlineto{\pgfpoint{0}{-0.5pt}}
        \pgfusepath{stroke}
        \pgfpathmoveto{\pgfpoint{0}{-0.5pt}}
        \pgfpathlineto{\pgfpoint{0}{2.5pt}}
        \pgfusepath{stroke}
        \end{pgfonlayer}
    }
}
\pgfdeclareshape{reverse arrow shape}
{
    \savedanchor\centerpoint
    {
        \pgf@x=0pt
        \pgf@y=0pt
    }
    \anchor{center}{\centerpoint}
    \anchor{south}{\centerpoint}
    \anchor{north}{\centerpoint}
    \anchorborder{\centerpoint}
    \backgroundpath
    {
        \begin{pgfonlayer}{foreground}
        \pgfsetstrokecolor{black}
        \pgfsetarrowsstart{<[\strarr]}
        \pgfsetlinewidth{\thickness}
        \pgfsetlinewidth{\thickness}
        \pgfpathmoveto{\pgfpoint{0}{-2pt}}
        \pgfpathlineto{\pgfpoint{0}{-1pt}}
        \pgfusepath{stroke}
        \end{pgfonlayer}
    }
}
\pgfdeclareshape{reverse double arrow shape}
{
    \savedanchor\centerpoint
    {
        \pgf@x=0pt
        \pgf@y=0pt
    }
    \anchor{center}{\centerpoint}
    \anchor{south}{\centerpoint}
    \anchor{north}{\centerpoint}
    \anchorborder{\centerpoint}
    \backgroundpath
    {
        \begin{pgfonlayer}{foreground}
        \pgfsetstrokecolor{black}
        \pgfsetarrowsstart{<[\strarr]}
        \pgfsetlinewidth{\thickness}
        \pgfpathmoveto{\pgfpoint{0}{-3.5pt}}
        \pgfpathlineto{\pgfpoint{0}{-0.5pt}}
        \pgfusepath{stroke}
        \pgfpathmoveto{\pgfpoint{0}{-0.5pt}}
        \pgfpathlineto{\pgfpoint{0}{2.5pt}}
        \pgfusepath{stroke}
        \end{pgfonlayer}
    }
}
\makeatother

\tikzset{arrow/.pic={
    \path[pic actions, decoration={ markings,
      mark=at position 1 with {\arrow{>[\strarr]}}
    },
    postaction={decorate}] (0,1pt) -- (0,2pt);
},double arrow/.pic={
    \path[pic actions, decoration={ markings,
      mark=at position \pgfdecoratedpathlength-3pt with {\arrow{>[\strarr]}},
      mark=at position \pgfdecoratedpathlength with {\arrow{>[\strarr]}}
    },
    postaction={decorate}] (0,-6pt) -- (0,4pt);
},
  reverse arrow/.pic={
    \path[pic actions, decoration={ markings,
      mark=at position 1 with {\arrow{<[\strarr]}}
    },
    postaction={decorate}] (0,1pt) -- (0,2pt);
}
}
